\documentclass[reqno,11pt]{amsart}
\usepackage{etex}
\usepackage{mathrsfs}
\usepackage{verbatim}
\usepackage{upgreek}
\usepackage{amsmath}
\usepackage{amsxtra}
\usepackage{amscd}
\usepackage{amsthm}
\usepackage{amsfonts}
\usepackage{amssymb}
\usepackage{eucal}
\usepackage{tikz-cd}
\usepackage[matrix,arrow,curve]{xy}
\usepackage{young}
\usepackage[vcentermath]{youngtab}

\usepackage{amsmath}

\usepackage{hyperref}
\usepackage[capitalise]{cleveref}
\makeatletter
\newtheoremstyle{break}
  {8pt}{8pt}%
  {\normalfont}{}%
  {\bfseries}{}%
  {\newline}%
  {\thmname{#1}\thmnumber{\@ifnotempty{#1}{ }\@upn{#2}}%
    \thmnote{ {\normalfont\itshape(#3)}}}
\theoremstyle{break}
\makeatother

\newcounter{theorem}

\theoremstyle{break}
\newtheorem{Thm}{Theorem}[section]
\newtheorem{Cor}[Thm]{Corollary}
\newtheorem{Lem}[Thm]{Lemma}
\newtheorem{Prop}[Thm]{Proposition}
\newtheorem{Conj}[Thm]{Conjecture}

\theoremstyle{break}
\newtheorem{Def}[Thm]{Definition}
\newtheorem{Rem}[Thm]{Remark}

\theoremstyle{definition}

\theoremstyle{remark}

\newtheorem{conj}{Conjecture}

\numberwithin{equation}{section}
\renewcommand{\rm}{\normalshape}

\newif\ifShowLabels
\ShowLabelstrue
\newdimen\theight
\def\TeXref#1{%
	\leavevmode\vadjust{\setbox0=\hbox{{\tt
				\quad\quad  {\small \rm #1}}}%
		\theight=\ht0
		\advance\theight by \lineskip
		\kern -\theight \vbox to
		\theight{\rightline{\rlap{\box0}}%
			\vss}%
}}%

\ShowLabelsfalse

\renewcommand{\sec}[2]{\section{#2}\label{S:#1}%
	\ifShowLabels \TeXref{{S:#1}} \fi}
\newcommand{\ssec}[2]{\subsection{#2}\label{SS:#1}%
	\ifShowLabels \TeXref{{SS:#1}} \fi}

\newcommand{\sssec}[2]{\subsubsection{#2}\label{SSS:#1}%
	\ifShowLabels \TeXref{{SSS:#1}} \fi}

\newenvironment{thm}[1]%
{ \begin{Thm} \label{T:#1}  \ifShowLabels \TeXref{T:#1} \fi }%
	{ \end{Thm} }

\renewcommand{\th}[1]{\begin{thm}{#1} \sl }
	\renewcommand{\eth}{\end{thm} }

\newenvironment{lemma}[2][]%
{ \begin{Lem}[#1] \label{L:#2}  \ifShowLabels \TeXref{L:#2} \fi }%
	{ \end{Lem} }
\newcommand{\lem}[2][]{\begin{lemma}[#1]{#2} \sl}
	\newcommand{\elem}{\end{lemma}}

\newenvironment{propos}[2][]%
{ \begin{Prop}[#1] \label{P:#2}  \ifShowLabels \TeXref{P:#2} \fi }%
	{ \end{Prop} }
\newcommand{\prop}[2][]{\begin{propos}[#1]{#2}\sl }
	\newcommand{\eprop}{\end{propos}}

\newenvironment{corol}[1]%
{ \begin{Cor} \label{C:#1}  \ifShowLabels \TeXref{C:#1} \fi }%
	{ \end{Cor} }
\newcommand{\cor}[1]{\begin{corol}{#1} \sl }
	\newcommand{\ecor}{\end{corol}}

\newenvironment{defeni}[2][]%
{ \begin{Def}[#1] \label{D:#2}  \ifShowLabels \TeXref{D:#2} \fi }%
	{ \end{Def} }
\newcommand{\defe}[2][]{\begin{defeni}[#1]{#2} \sl }
	\newcommand{\edefe}{\end{defeni}}

\newenvironment{remark}[1]%
{ \begin{Rem} \label{R:#1}  \ifShowLabels \TeXref{R:#1} \fi }%
	{ \end{Rem} }
\newcommand{\rem}[1]{\begin{remark}{#1}}
	\newcommand{\erem}{\end{remark}}

\newenvironment{conjec}[1]%
{ \begin{Conj} \label{Co:#1}  \ifShowLabels \TeXref{Co:#1} \fi }%
	{ \end{Conj} }
\renewcommand{\conj}[1]{\begin{conjec}{#1} \sl }
	\newcommand{\econj}{\end{conjec}}

\newenvironment{First proof}[1]%
{ \begin{First proof} \label{Co:#1}  \ifShowLabels \TeXref{Co:#1} \fi }%
	{ \end{First proof} }

\newenvironment{Second proof}[1]%
{ \begin{Second proof} \label{Co:#1}  \ifShowLabels \TeXref{Co:#1} \fi }%
	{ \end{Second proof} }

\newcommand{\eq}[1]%
{ \ifShowLabels \TeXref{E:#1} \fi
	\begin{equation} \label{E:#1} }
\newcommand{\eeq}{ \end{equation} }

\newcommand{\prf}{ \begin{proof} }
	\newcommand{\epr}{ \end{proof} }

\newcommand{\prft}{ \begin{proof} }
	\newcommand{\eprt}{ \end{proof} }

\newcommand\nc{\newcommand}

\nc{\on}{\operatorname}
\nc{\BA}{{\mathbb{A}}}
\nc{\C}{{\mathbb{C}}}
\nc{\BN}{{\mathbb{N}}}
\nc{\BZ}{{\mathbb{Z}}}
\nc{\BP}{{\mathbb{P}}}
\nc{\CA}{{\mathcal{A}}}
\nc{\CE}{{\mathcal{E}}}
\nc{\CF}{{\mathcal{F}}}
\nc{\CG}{{\mathcal{G}}}
\nc{\CH}{{\mathcal{H}}}
\nc{\CK}{{\mathcal{K}}}
\nc{\CL}{{\mathcal{L}}}
\nc{\CM}{{\mathcal{M}}}
\nc{\CN}{{\mathcal{N}}}
\nc{\CO}{{\mathcal{O}}}
\nc{\CP}{{\mathcal{P}}}
\nc{\CR}{{\mathcal{R}}}
\nc{\CT}{{\mathcal{T}}}
\nc{\CU}{{\mathcal{U}}}
\nc{\CV}{{\mathcal{V}}}

\nc{\Cx}{{\mathbb{C}^\times}}

\nc{\fm}{{\mathfrak{m}}}
\nc{\unl}{\underline}
\nc{\ol}{\overline}
\newcommand\iso{\,\vphantom{j^{X^2}}\smash{\overset{\sim}{\vphantom{\rule{0pt}{0.20em}}\smash{\longrightarrow}}}\,}
\nc{\ul}{\underline}
\nc{\Mvd}{\mathfrak{M}(\underline{v},\underline{d})}
\nc{\MvdT}{\mathfrak{M}(\underline{v}^{\dagger},\underline{d}^{\dagger})}
\nc{\MVD}{\mathfrak{M}(V,D)}
\nc{\mt}{\mapsto}
\nc{\sm}{\setminus}
\nc{\ra}{\rightarrow}
\nc{\lar}{\leftarrow}
\nc{\hr}{\hookrightarrow}
\nc{\La}{\Lambda}
\nc{\Lap}{\Lambda^{+}}
\nc{\oZal}{\overset{\circ}{Z^{\alpha}}}
\nc{\sig}{\sigma}
\nc{\al}{\alpha}
\nc{\la}{\lambda}
\nc{\is}{\simeq}
\nc{\ip}{\iota^{+}_{\la, \mu}}
\nc{\im}{\iota^{-}_{\la, \mu}}
\nc{\jp}{j^{+}_{\la, \mu}}
\nc{\jm}{j^{-}_{\la, \mu}}
\nc{\pip}{\pi^{+}_{\la, \mu}}
\nc{\pim}{\pi^{-}_{\la, \mu}}
\nc{\s}{\star}
\nc{\fpt}{[A^{\la},B^{\la},\gamma^{\la},\delta^{\la}]}
\nc{\ulfpt}{[\ul{A}^{\la},\ul{B}^{\la},\ul{\gamma}^{\la},\ul{\delta}^{\la}]}
\nc{\Hom}{\on{Hom}}
\nc{\Fl}{\on{\mathbf{Fl}}}
\nc{\Gr}{\on{\mathbf{Gr}}}

\newcommand{\BC}{\mathbb{C}}
\newcommand{\CW}{\mathcal{W}}

\author{Vasily Krylov}
\address{Department of Mathematics
Massachusetts Institute of Technology
\newline
77 Massachusetts Avenue,
Cambridge, MA 02139,
USA;
\newline National Research University Higher School of Economics, Russian Federation\newline
Department of Mathematics, 6 Usacheva st., Moscow 119048;
}
\email{krvas@mit.edu, kr-vas57@yandex.ru}
\author{Ivan Perunov}
\address{National Research University Higher School of Economics, Russian Federation\newline
Department of Mathematics, 6 Usacheva st., Moscow 119048;\newline
Skolkovo Institute of Science and Technology
}
\email{ivan-perunov@ya.ru}
\title{Almost dominant generalized slices and convolution diagrams over them}
\begin{document}

     \begin{abstract}
     Let $G$ be a connected reductive complex algebraic group with a maximal torus $T$. We denote by $\La$ the coweight lattice of $T$.
     Let $\La^+ \subset \La$ be the submonoid of dominant coweights. For $\la \in \La^+,\,\mu \in \La,\,\mu \leqslant \la$, in~\cite{BFN}, authors defined a generalized transversal slice $\ol{\CW}^\la_\mu$. This is an algebraic variety of the dimension $\langle 2\rho^{\vee}, \la-\mu \rangle$,
     where $2\rho^{\vee}$
    is the sum of positive roots of $G$. 
    In this paper, 
    we construct an isomorphism $\ol{\CW}^\la_\mu \simeq \ol{\CW}^\la_{\mu^+} \times \BA^{\langle2\rho^{\vee},\, \mu^+-\mu\rangle}$ for $\mu \in \La$ such that $\langle \al^{\vee},\mu\rangle \geqslant -1$ for any positive root $\al^{\vee}$,   
    here $\mu^+ \in W\mu$ is the dominant representative in the Weyl group orbit of $\mu$.                           
    We consider the example when $\la$ is minuscule, $\mu \in W\la$ and describe natural coordinates, Poisson structure on $\ol{\CW}^\la_\mu \simeq \BA^{\langle 2\rho^\vee,\,\la-\mu \rangle}$ and its $T\times \BC^\times$-character. We apply these results to compute $T \times \BC^\times$-characters of tangent spaces at fixed points of convolution diagrams $\widetilde{\CW}^{\ul{\la}}_\mu$ with minuscule $\la_i$.
    We also apply our results to construct open coverings by affine spaces of convolution diagrams $\widetilde{\CW}^{\ul{\la}}_\mu$ over slices with $\mu$ such that $\langle \al^{\vee},\mu\rangle \geqslant -1$ for any positive root $\al^{\vee}$ and minuscule $\la_i$ and to compute Poincar\'e polynomials of such convolution diagrams       $\widetilde{\CW}^{\ul{\la}}_{\mu}$.
    \end{abstract}    
     \maketitle
   
\sec{Int}{Introduction}

\ssec{Int}{Generalized transversal slices}
$G$ is a connected reductive complex algebraic group with a maximal torus $T \subset G$. In~\cite{BFN}, the authors constructed a generalized transversal slice $\ol{\CW}^\la_\mu$, which depends on a pair of coweights $\la,\,\mu$ of $T$ such that $\la$ is dominant and $\mu \leqslant \la$ with respect to the dominance order on $\La$. 
Let $\on{Gr}_G$ be the affine Grassmannian of $G$. Any coweight $\mu$ gives rise to a point of $\on{Gr}_G$ to be denoted $z^\mu$. We set $\CO:=\BC[[z]]$. It is known (see~\cite[Section~2]{BF1}) that if a coweight $\mu$ is dominant, then the variety $\ol{\CW}^\la_\mu$ coincides with the transversal slice to the $G(\CO)$-orbit $G(\CO)\cdot z^\mu$ inside $\ol{G(\CO) \cdot z^\la}$. 
It is also known that for $\la=0$ the variety $\ol{\CW}^\la_\mu$ parametrizes  the based maps of degree $w_0(\mu)$ from $\BP^1$ to the flag variety $\mathcal{B}:=G/B$ (so-called open zastava space). In general, we have a locally closed embedding $\imath\colon \ol{\CW}^\la_\mu \hookrightarrow \ol{\on{Gr}}{}^\la_G \times Z^{-w_0(\la-\mu)}$, where $\ol{\on{Gr}}{}^\la_G:=\ol{G(\CO)\cdot z^\la} \subset \on{Gr}_G$ and $Z^{-w_0(\la-\mu)}$ is the space of based quasi-maps of degree $-w_0(\la-\mu)$ from $\BP^1$ to $\mathcal{B}$. 

\ssec{Int}{Main results and structure of the paper}
The paper is organized as follows. In Section~\ref{defnot}, we give the definitions of the main geometric objects of our study and formulate their basic properties. 
In Section~\ref{main_result_fibration}, we 
construct an isomorphism $A^\la_\mu \simeq p^{-1}(z^\mu) \times U^\mu_- \times \ol{\CW}^\la_{\mu^+}$ (see Proposition~\ref{fibration}), here 
$A^\la_\mu \subset \ol{\CW}^\la_\mu$ is the attractor with respect to the loop rotation action, $p\colon \ol{\CW}^{\la}_\mu \ra \ol{\on{Gr}}{}^{\la}_G$ is the natural morphism and $U^\mu_-$ is a certain subgroup of $U_-$. In the case when $\langle \al^{\vee},\mu \rangle \geqslant -1$ for every positive root $\al^{\vee}$ we have $A^\la_\mu=\ol{\CW}^\la_\mu$ (see Proposition~\ref{loop_contracts}) and we obtain the isomorphism $\ol{\CW}^\la_\mu \simeq \ol{\CW}^{\la}_{\mu^+}\times \BA^{\langle2\rho^{\vee},\, \mu^+-\mu\rangle}$ (see Theorem~\ref{affine_fibr}).
In Section~\ref{minuscule_case} we consider the example when $\la$ is minuscule and $\mu \in W\la$. We construct natural coordinates on $\ol{\CW}^\la_\mu$ (see Theorem~\ref{coordinates_and_isomorphisms}) and describe Poisson structure on $\ol{\CW}^\la_\mu$ (see Theorem~\ref{poiss_structure}).
In Section~\ref{application}, we study convolution diagrams $\widetilde{\CW}^{\ul{\la}}_\mu$ over slices $\ol{\CW}^\la_\mu$ and prove that for $\mu$ such that $\langle \al^{\vee},\mu \rangle \geqslant -1$ for every positive root $\al^{\vee}$ they can be covered by the images of so-called multiplication morphisms $\tilde{m}^{\ul{\la}}_{\ul{\mu}}\colon \ol{\CW}^{\la_1}_{\mu_1} \times \ldots \times \ol{\CW}^{\la_N}_{\mu_N} \ra \widetilde{\CW}^{\ul{\la}}_\mu$ with $\mu_i$ being weights of $V^{\la_i}$ (see Theorem~\ref{covering_figovering}). 
As a result we obtain open coverings by affine spaces of convolution diagrams $\widetilde{\CW}^{\ul{\la}}_\mu$ over slices for minuscule $\la_i$ (Theorem~\ref{theo_covering_min}).
We compute characters of tangent spaces at fixed points of convolution diagrams $\widetilde{\CW}^{\ul{\la}}_{\mu}$ for minuscule $\la_i$ (Remark~\ref{characters_resolv}) and use this computation to find Poincar\'e polynomials (i.e. to compute cohomology groups with compact support) of convolution diagrams $\widetilde{\CW}^{\ul{\la}}_{\mu}$ such that $\langle \al^{\vee},\mu \rangle \geqslant -1$ for $\al^{\vee} \in \Delta^{\vee}_+$ (Remark~\ref{Poincare_pol}).
We also prove that without any restrictions on $\mu$
convolution diagrams $\widetilde{\CW}^{\ul{\la}}_\mu$ are covered by the images of $\tilde{m}^{\ul{\la}}_{\ul{\mu}}$ with $\mu_i$ being a weight of $V^{\la_i}$ for all $\mu_i$ possibly except the last one (see Theorem~\ref{covering_general}). 

\ssec{Int}{Acknowledgements} 
We would like to thank Hiraku Nakajima for 
lots of explanations especially on the material of Section~\ref{application}. The main idea in the proof of Theorem~\ref{theo_covering_min} belongs to Hiraku Nakajima. We would also like to thank our advisor Michael Finkelberg for many helpful discussions and numerous explanations.
We would also like to thank Alexander Braverman,
Dinakar Muthiah and Alex Weekes for helpful discussions, suggestions and comments.
We are gratefull to Dinakar Muthiah for explaining us the modification of the proof of Theorem~\ref{affine_fibr} that (potentially) works in any charactersistic (our original proof used the fact that the characteristic of the base field is $0$). 
We would like to thank the anonymous referee for carefull proofreading of the text, very usefull comments and suggestions, in particular for strengthening our original version of Theorem~\ref{poiss_structure}. 
V.K. was partially supported by the grant RSF 19-11-00056.
    
\sec{defnot}{Definitions and Basic Properties}\label{defnot}
\ssec{defnot}{Main objects}   
We fix a triple $G\supset B\supset T$, consisting of a connected reductive algebraic group over $\BC$, a Borel subgroup $B$ and a maximal torus $T$; $\mathfrak{g}\supset\mathfrak{b}\supset\mathfrak{t}$ are their Lie algebras. We denote by $B_- \supset T$ the opposite Borel subgroup of $G$.
Also we denote by $U$ and $U_-$ the unipotent radicals of $B$ and $B_-$ respectively.
    We denote by $\Lambda$ the coweight lattice of $T$ that is by definition the image of the coharacter lattice $\on{Hom}(\BC^\times,T)$ via the map $\on{Hom}(\BC^\times,T) \ra \mathfrak{t}$ sending a cocharacter $\eta$ to $(d_1\eta)(1)$. We denote by $\Lambda^+ \subset \Lambda$ the submonoid of dominant coweights. We similarly denote the lattice $\La^\vee$ of characters of $T$ and the submonoid $\La^{\vee+} \subset \La^\vee$ of dominant characters.
    We denote by $\Delta^{\vee}$ (resp. $\Delta$) the set of roots (resp. coroots) of $(T,G)$ and by $\Delta^{\vee}_+$ (resp. $\Delta_+$) the set of positive roots (resp. coroots) with respect to the Borel subgroup $B \subset G$. 
    We also denote by $W$ the Weyl group of $(T,G)$ and by $w_0 \in W$ its longest element with respect to the  Borel $B$.
    Set $\rho^\vee=\frac{1}{2}\sum_{\al^\vee \in \Delta^\vee_+}\al^\vee$, $\rho=\frac{1}{2}\sum_{\al \in \Delta_+}\al$. For two coweights $\la,\, \mu \in \La$ we say that $\la \geqslant \mu$ if $\la-\mu$ can be written as the sum of positive coroots with integer nonnegative coefficients.

    \defe[Affine Grassmannian]{affgr}
        We define $\on{Gr}_G$ as the moduli space of the data $(\CP,\sigma)$, where
        \item
    (a) $\mathcal{P}$ is a G-bundle  on $\mathbb{P}^{1};$
    \item
    (b) $\sig\colon \mathcal{P}^{triv}|_{\mathbb{P}^{1} \sm \{0\}} \iso \mathcal{P}|_{\mathbb{P}^{1} \setminus \{0\}}$ is 
    a trivialization of $\CP$ restricted to $\BP^1 \setminus \{0\}$.
    \edefe

    The set of $\BC$-points of $\on{Gr}_G$ can be described as follows: set $\CK: =\BC((z)),\, \CO:=\BC[[z]]$, then $\on{Gr}_G(\BC)$ is the quotient $G(\CK)/G(\CO)$.
    Any coweight $\la\colon \BC^\times \ra T$ is an element of $G(\BC[z^{\pm 1}])$ so defines an element of $\on{Gr}_G$ to be
    denoted by $z^\la$, also for every $t \in \BC^\times$ we denote by $t^\mu \in T$ the element $\mu(t)$.
    The group $G(\CO)$ acts on $\on{Gr}_G$ via left multiplication. For $\la \in \Lambda^+$, denote by
    $\on{Gr}^\la_G$ the $G(\CO)$-orbit of $z^\la$.
   We have the following decompositions:
        \begin{equation}\label{decomp}
        \begin{aligned}
       &\on{Gr}_G=\bigsqcup_{\lambda \in \Lambda^+} \on{Gr}_G^\lambda, \\
       &\ol{\on{Gr}}{}^\lambda_G=\bigsqcup_{
       \mu \leqslant \la,\, \mu \in \La^+
       } \on{Gr}_G^\mu.
       \end{aligned}
       \end{equation}
    It is known that for any $\la \in \Lambda^+$ $\ol{\on{Gr}}{}^\la_G$ is a projective algebraic variety of dimension $\langle 2\rho^{\vee},\la \rangle$. It follows that $\on{Gr}_G= \underset{\longrightarrow}{\on{lim}}\,\ol{\on{Gr}}{}^\la_G$ is an ind-projective scheme.
    
    \defe[Generalized transversal slice]{slice}\label{def_genslice}
    Let $\lambda$ be a dominant coweight, let $\mu\leqslant\lambda$ be any coweight.
    Following~\cite{BFN}, we define the generalized transversal slice in the affine Grassmannian $\ol{\mathcal W}_\mu^\lambda$. 
    It is the moduli space of the data
    $(\mathcal{P},\sigma,\phi)$, where 
    \item
    (a) $\mathcal{P}$ is a G-bundle  on $\mathbb{P}^{1};$
    \item
    (b) $\sig\colon \mathcal{P}^{triv}|_{\mathbb{P}^{1} \sm \{0\}} \iso \mathcal{P}|_{\mathbb{P}^{1} \setminus \{0\}}$ --
    a trivialization, having a pole of degree $\leqslant \la$. This means that the point $(\CP, \sig) \in \on{Gr}_G$ lies in $\ol{\on{Gr}}{}^{\la}_G$;
    \item
    (c) $\phi$ is a $B$-structure on $\CP$ (i.e. a $B$-subbundle of  $\CP$) of degree $w_{0}(\mu)$, having no defect at $\infty$ and having fiber $B_{-}$ at $\infty$ (with respect to $\sig$). 
    
    \edefe
    It follows from~\cite[Lemmas~2.5,~2.7]{BFN} that $\ol{\CW}^\la_\mu$ is an affine algebraic variety of dimension $\langle 2\rho^{\vee}, \la-\mu \rangle$. We will denote by $\CW^\la_\mu \subset \ol{\CW}^\la_\mu$ the open subvariety consisting of $(\CP,\sigma,\phi)$ such that $\sigma$ has a pole of degree exactly $\la$. It follows from~\cite{smooth_slice} that $\CW^\la_\mu$ is a smooth variety. One can show (see Remark~\ref{smooth_transv_crit}) that $\CW^{\la}_{\mu}$ is the smooth locus of $\ol{\CW}^\la_\mu$.

     Let us denote by $J$ the set of simple coroots $\al_i,\, i \in J$ of $G$.

      \prop{irred}\label{irred}
    The variety $\ol{\CW}^\la_\mu$ is irreducible.
    \eprop
    \prf
    Set $\al=\la-\mu$ and $\al^*=-w_0(\al)$. Let us define the space $\BA^{\al^*}$ of colored effective divisors of multidegree $\al^*$.
    For any $n \in \BZ_{\geqslant 0}$ we denote the $n$-th symmetric power of the curve $\BA^1$ by $\BA^{(n)}$. A point $D \in \BA^{\al^*}$ is a collection of effective divisors $D_{\la^\vee} \in \BA^{(\langle\la^\vee ,\al^*\rangle)}$ for $\la^\vee \in \La^{\vee+}$ such that $D_{\la_1^\vee}+D_{\la_2^{\vee}}=D_{\la_1^\vee+\la_2^\vee}$.  
    We can write $\al^*=\sum_{i \in J} a_i\al_i$ for some $a_i \in \BZ_{\geqslant 0}$. 
    Note that when the derived group $G^{\mathrm{der}}=[G,G]$ is simply-connected then $\BA^{\al^*}=\prod_{i \in J}\BA^{(a_i)}$ and in general we have a closed (``diagonal") embedding $\prod_{i \in J}\BA^{(a_i)} \hookrightarrow \BA^{\al^*}$.
    
    Recall the factorization morphism $\pi\colon \ol{\CW}^\la_\mu \ra \mathbb{A}^{\al^*}$ (see~\cite[Lemma~2.7]{BFN}). The variety $\mathbb{A}^{\al^*}$ is irreducible, morphism $\pi$ is flat (\cite[Lemma~2.7]{BFN}) and there exists an open dense subset $U \subset \mathbb{A}^{\al^*}$, such that $\pi^{-1}(U)$ is irreducible (see~\cite[Section~2(ix)]{BFN}). So the closure $\ol{\pi^{-1}(U)}$ is an irreducible component of $\ol{\CW}^\la_\mu$. Suppose that there exists an other irreducible component $X\subset \ol{\CW}^\la_\mu$. 
    We denote by $\overset{\circ}{X}$ the open dense subvariety of $X$, consisting of points which do not lie in any other irreducible component.
    It follows from the flatness of $\pi$ and irreducibility of $\BA^{\al^*}$ that the restriction $\pi|_{\overset{\circ}{X}}\colon  \overset{\circ}{X} \ra \mathbb{A}^{\al^*}$ is dominant. It contradicts to the fact that ${\overset{\circ}{X}} \cap \ol{\pi^{-1}(U)} = \varnothing$.
    \epr

Let us denote by $\La^{\on{pos}} \subset \Lambda$ the submonoid of $\La$, spanned by the simple coroots $\alpha_{i},\,i \in J$.

    \defe[Zastava]{zast}
   For $\al \in \La^{\on{pos}}$, we define $Z^\al$ as a moduli space of degree $\al$ quasi-maps $f \in \on{Qmaps}^\al(\BP^1,\mathcal{B})$ to the full flag variety $\mathcal{B}$
   having no defect at $\infty$ and such that $f(\infty)=B_-$
   (see~\cite{b} for the definition of the notion of a quasi-map, see also~\cite{FM} for other equivalent definitions of $Z^\al$).
    This is an algebraic variety of dimension $\langle 2\rho^{\vee}, \al \rangle$. We define by $\oZal \subset Z^\al$ the open subvariety of $Z^\al$ consisting of actual maps $f\colon \BP^1 \ra \mathcal{B}$ of degree $\al$ and such that $f(\infty)=B_-$. 
   \edefe

\rem{}
{\em{Note that a map $f\colon \BP^1 \ra \mathcal{B}$ of degree $\al$ and such that $f(\infty)=B_-$ is the same as the degree $\al$ $B$-structure in the trivial $G$-bundle $\CP^{triv}$  having fiber $B_-$ at $\infty$. It follows that the variety $\oZal$ parametrizes degree $\al$ $B$-structures in the trivial $G$-bundle over $\BP^1$ having fiber $B_-$ at $\infty$.
So  we see that $\oZal=\ol{\CW}^{0}_{w_0(\al)}$.
}}
\erem

\ssec{}{Action of $B_- \times \BC^\times$ and matrix description of slices}\label{matrix}\label{actionss}    
    Group
    $B_-$ acts on $\ol{\mathcal W}_\mu^\lambda,\,\on{Gr}_G,\, Z^\al$ via changing the trivialization. 
    We also have the natural
    $\Cx$-action on $\ol{\mathcal W}_\mu^\lambda,\,\on{Gr}_G,\,Z^\al$ to be called loop rotation action. It is induced from the following action  on $\BP^1$: 
    $(x:y) \mapsto (tx:y)$, here $0=(1:0),\,\infty = (0:1)$.

    For a complex algebraic group $H$, we recall that $H[z]:=H(\BC[z]),\,H[[z^{-1}]]:=H(\BC[[z^{-1}]])$ and denote by $H[[z^{-1}]]_{1}$ the kernel (preimage of $1 \in H$) of the natural evaluation at $\infty$
morphism $H[[z^{-1}]] \ra H$.
    \rem{}
    {{
    Note that the action of $B_-$ on $\on{Gr}_G$ extends to the action of the whole group $G(\CK) \curvearrowright \on{Gr}_G$ via changing the trivialization. Clearly this action does not extend to $\ol{\CW}^\la_\mu$. On the other hand, one can consider the ind-scheme $\CW_\mu^{\mathrm{rat}}:=\underset{\longrightarrow}{\on{lim}}\,\ol{\CW}^\la_\mu$ parametrizing triples $(\CP,\sigma,\phi)$ as in Definition~\ref{def_genslice}, where we put no restrictions on the defect of $\sigma$. Then we have the action $G[z^{-1}]_{B_-} \curvearrowright \CW_\mu^{\mathrm{rat}}$ (and even of the group $G[[z^{-1}]]_{B_-} \cap G(z)$) via changing the trivialization and this action extends the action of $B_-$ on $\ol{\CW}^\la_\mu$, here $G[z^{-1}]_{B_-} \subset G[z^{-1}] \subset G(\CK)$ is the preimage of $B_-$ with respect to the natural evaluation at $\infty$ morphism $G[z^{-1}] \ra G$. For details, see the proof of Proposition~\ref{fibration}.
    }}
    \erem
    
    



    In~\cite[Section~2($\on{xi}$)]{BFN}, the following isomorphism was constructed:
    \begin{equation}\label{matrixmatrix}
    \Psi\colon\ol{\CW}^{\la}_{\mu} \simeq \Big( U[[z^{-1}]]_{1}z^{\mu}B_{-}[[z^{-1}]]_{1}\cap \ol{G[z]z^{\la}G[z]}\Big),
    \end{equation}
    where the right hand side is considered as a locally closed subvariety in the ind-scheme $G((z^{-1})):=G(\BC((z^{-1})))$. 
    
    \lem{action via matrix}\label{action_via_matrix}
    (a) The $\BC^\times$-action on $\ol{\CW}^\la_\mu$ via loop rotation comes from the following action $\BC^\times \curvearrowright G((z^{-1}))$:
    $
    g(z) \mapsto g(t^{-1}z) \cdot t^{\mu}
    $.
    
    \item
    (b) The $T$-action on $\ol{\CW}^\la_\mu$ 
    comes from the action on $G((z^{-1}))$ via conjugation.
    \elem
    
    \prf
To prove $(a)$, let us recall the construction of the isomorphism $\Psi$ of~\cite[Section~2(xi)]{BFN}.
    Take a point $(\CP,\sigma,\phi) \in \ol{\CW}^\la_\mu$. Let $\CP_B$ be the $B$-bundle of degree $w_0(\mu)$ that corresponds to $\phi$. We denote by $\CP_{B_-}$ the corresponding $B_-$-bundle of degree $\mu$. Fix a trivialization
    $\sigma_{B_-}\colon (\CP^{triv}_{B_-})|_{\mathbb{A}^1} \iso (\CP_{B_-})|_{\mathbb{A}^1}$ of $\CP_{B_-}$, restricted to $\mathbb{A}^1$.
    Fix also a trivialization $\sigma_{U}\colon \CP^{triv}_{U}|_{\mathbb{A}^1} \iso \CP^{triv}_{U}|_{\mathbb{A}^1}$ of the trivial $U$-bundle $\CP^{triv}_U$, restricted to $\mathbb{A}^1$.
    Now $\sigma\colon \CP^{triv}|_{\BP^1 \setminus \{0\}} \iso \CP|_{\BP^{1}\setminus \{0\}}$ defines an element of $\overline{G[z]z^{\la}G[z]} \subset  G(z)$ (given by $\on{Ind}^G_{U}(\sigma^{-1}_{U}) \circ \sigma^{-1} \circ \on{Ind}^G_{B_-}(\sigma_{B_-})$), well-defined up to the right multiplication by $B_-[z]$ and the left multiplication by $U[z]$. After that we consider an embedding
    $G(z) \hookrightarrow G((z^{-1}))$. Due to the condition on $\phi$ at $\infty$ and the fact that the degree of $\CP_{B_-}$ equals to $\mu$ we obtain an element of
    $U[z] \backslash U((z^{-1})) z^{\mu}T[[z^{-1}]] U_{-}((z^{-1}))/B_-[z]$. 
    It lifts uniquely to an element 
    $g\in U[[z^{-1}]]_1z^{\mu}B_-[[z^{-1}]]_1$.
    We set 
    \begin{equation*}
    \Psi(\CP,\sigma,\phi):=g \in U[[z^{-1}]]_1z^{\mu}B_{-}[[z^{-1}]]_1 \cap \overline{G[z]z^{\la}G[z]}.
    \end{equation*}

    Let us also describe the inverse morphism $\Psi^{-1}$. 
    We take
    \begin{equation*}
    g \in B[[z^{-1}]]_1z^{\mu}B_-[[z^{-1}]]_1 \cap \overline{ G[z] z^\lambda G[z]},
    \end{equation*}
    it defines a transition function for a $G$-bundle $\CP$ together with trivializations $\sigma_{0}\colon \CP^{triv}|_{\BP^1\setminus \{0\}} \iso \CP|_{\BP^1\setminus \{0\}}$ and 
    $\sigma_{\infty}\colon \CP^{triv}|_{\mathbb{A}^1} \iso \CP|_{\mathbb{A}^1}$ ($g$ is a rational function which corresponds to the composition $\sigma^{-1}_0 \circ \sigma_{\infty}$). Let $\phi_{-}$ be the image of the standard $B_-$-structure in 
    $\CP^{triv}$ under the morphism $\sigma_{\infty}$.
    Let $\phi$ be the corresponding $B$-structure. Then $\Psi^{-1}(g)=(\CP,\sigma_{0},\phi)$.
    
    Now let us prove $(a)$.
    Fix an element $t^{-1} \in \BC^\times$, $(\CP,\sigma,\phi) \in \ol{\CW}^\la_\mu$ and set $g:=\Psi(\CP,\sigma,\phi)$, we also denote by $\phi_-$ the $B_-$-structure which corresponds to $\phi$. Our goal is to compute $t^{-1} \cdot g:= \Psi(t^{-1} \cdot (\CP,\sigma,\phi))$. Note that $t^{-1} \cdot (\CP,\sigma,\phi)=(t^*\CP,t^*\sigma,t^*\phi)$ where 
    $t^{*}\CP,\, t^*\phi$ are the pullbacks of $\CP,\, \phi$, and $t^{*}\sigma\colon (t^{*}\CP^{triv})|_{\BP^1\setminus \{0\}} \iso (t^{*}\CP)|_{\BP^1\setminus \{0\}}$ is the corresponding trivialization.
    Recall also that the point $g \in G(z)$ gives us the trivialization $\sigma_{\infty}\colon (\CP^{triv})|_{\mathbb{A}^1} \iso (\CP)|_{\mathbb{A}^1}$. We denote by $t^*\sigma_{\infty}$ the corresponding pullback.
    
    Consider now the automorphism $t^{\mu}\colon \CP^{triv} \iso \CP^{triv}$. Note that the point $(t^*\CP,t^{*}\sigma,t^{*}\sigma_{\infty} \circ t^\mu)=g(t^{-1}z)t^{\mu}$ lies in $U[[z^{-1}]]_{1}z^{\mu}B_{-}[[z^{-1}]]_{1}\cap \ol{G[z]z^{\la}G[z]}$ (directly follows from the fact that $g(z)$ lies in this intersection). It remains to show that $\Psi^{-1}(g(t^{-1}z)t^{\mu})=(t^*\CP,t^*\sigma,t^*\phi)$.
    Recall that the $B_-$-structure $\phi_-$ is the image of the standard $B_-$ structure in $\CP^{triv}$ under the morphism $\sigma_{\infty}$.
    Let us denote this standard $B_-$-structure by $\phi^{stand}_-$. It follows that $t^{*}\phi_-$ is the image of $t^*\phi^{stand}_-=\phi^{stand}_-$ under $t^*\sigma_{\infty}$. The automorphism $t^{\mu}$ preserves $\phi^{stand}$ so $t^{*}\phi_-$ is the image of $\phi^{stand}_-$ under $t^*\sigma_{\infty} \circ t^{\mu}$. It now directly follows from the definitions that $\Psi^{-1}(g(t^{-1}z)t^{\mu})=(t^*\CP,t^*\sigma,t^*\phi)$.

    Part $(b)$ can be proved analogously.
    \epr
    
   \rem{}{}\label{act_U-_matrix}
{{    
Let us point out that we have the following description of the $B_-$-action (see Definition~\ref{actionss}) in matrix terms (we are grateful to Alex Weekes and Dinakar Muthiah for explaining this to us). Set 
\begin{equation*}
\CW_\mu := U[[z^{-1}]]_1 z^\mu T[[z^{-1}]]_1U_{-}[[z^{-1}]]_1,\, X_\mu:=U[[z^{-1}]]_1 z^\mu T[[z^{-1}]]_1U_{-} ((z^{-1})).
\end{equation*} 
We have the natural projection morphism $p\colon X_\mu \twoheadrightarrow \CW_\mu$ and the natural section $i\colon \CW_\mu \hookrightarrow X_\mu$ of the morphism $p$. We also have the Gauss decomposition
\begin{equation}\label{Gauss}
G[[z^{-1}]]_1=U[[z^{-1}]]_1 T[[z^{-1}]]_1 U_-[[z^{-1}]]_1.
\end{equation} We claim that the conjugation action of $B_- \curvearrowright G((z^{-1}))$ restricts to the $B_-$-action on $X_\mu$. To see this let us fix a point $u z^{\mu} t u_- \in X_\mu,\, u \in U[[z^{-1}]]_1, t \in T[[z^{-1}]]_1, u_- \in U_-((z^{-1}))$ and an element $g \in B_-$.
Using~(\ref{Gauss}) we decompose $gug^{-1}=\tilde{u}\tilde{t}\tilde{u}_-,
\, \tilde{u} \in U[[z^{-1}]]_1, \tilde{t} \in T[[z^{-1}]]_1, \tilde{u}_- \in U_-[[z^{-1}]]_1$. We obtain 
\begin{equation*}
guz^{\mu}tu_-g^{-1}=(gug^{-1})gz^{\mu}tu_-g^{-1}=
(\tilde{u}\tilde{t}\tilde{u}_-)gz^{\mu}tu_-g^{-1}=
\tilde{u}\cdot (z^{\mu} \cdot \tilde{t}t \cdot t^{-1}z^{-\mu})\tilde{u}_-g z^{\mu}tu_-g^{-1}
\end{equation*}
and see that $\tilde{u} \in U[[z^{-1}]],\, \tilde{t}t \in T[[z^{-1}]]_1,\, t^{-1}z^{-\mu}\tilde{u}_-gz^{\mu}tu_-g^{-1} \in U_-((z^{-1}))$, so $guz^{\mu}tu_-g^{-1} \in X_\mu$.
It follows that the conjugation action $B_- \curvearrowright G((z^{-1}))$ restricts to the $B_-$-action on $X_\mu$. Therefore, we obtain the action of $B_-$ on $\CW_\mu$ given by $g \cdot x = p(g i(x)g^{-1})=p(g i(x))$, here $g \in B_-,\, x \in \CW_\mu$.
This action restricts to the desired $B_-$-action on $\ol{\CW}^{\la}_{\mu} \subset \CW_\mu$. }}
\erem

\rem{} 
Another similar way to describe the action $U_- \curvearrowright \ol{\CW}^\la_\mu$ and more generally $G[[z^{-1}]]_{U_-} \curvearrowright \CW_\mu$ is the following. Recall that 
\begin{equation*}
\CW_\mu=U[[z^{-1}]]_1 z^\mu B_{-}[[z^{-1}]]_1   = U[[z^{-1}]]_1 z^{\mu} T[[z^{-1}]]_1 U((z^{-1}))/U[z].
\end{equation*}
Recall also the Gauss decomposition $G[[z^{-1}]]_1=U[[z^{-1}]]_1 T[[z^{-1}]]_1 U_-[[z^{-1}]]_1$ which implies the decomposition $G[[z^{-1}]]_{U_-}=U[[z^{-1}]]_1 T[[z^{-1}]]_1 U_-[[z^{-1}]]$. It follows that $U[[z^{-1}]]_1 z^{\mu} T[[z^{-1}]]_1 U_-((z^{-1}))/U_-[z]=G[[z^{-1}]]_{U_-} z^\mu T[[z^{-1}]]_1 U_-((z^{-1}))/U_-[z]$ so we conclude that 
\begin{equation*}
\CW_\mu=G[[z^{-1}]]_{U_-} z^\mu T[[z^{-1}]]_1 U_-((z^{-1}))/U_-[z]    
\end{equation*}
and the action of $G[[z^{-1}]]_{U_-}$ in these terms simply corresponds to the action via left multiplication (compare with~\cite[Section~$3.1$]{y}). The action of $G[[z^{-1}]]_{U_-} \curvearrowright \CW_\mu$ restricts to the desired $U_-$-action on $\ol{\CW}^\la_\mu \subset \CW_\mu$.
\erem
    
\ssec{}{Cartan involution via matrix description}

In~\cite{BFN} the authors defined a certain  involution $\iota\colon \ol{\CW}^{\la}_{\mu} \iso \ol{\CW}^{\la}_{\mu}$ (see~\cite[Section~2(vii)]{BFN}) that they call the Cartan involution of $\ol{\CW}^{\la}_\mu$.
This involution is defined using the moduli description of $\ol{\CW}^{\la}_\mu$. The goal of this Section is to describe $\iota$ using the matrix description of slices (see Proposition~\ref{cartan_via_matrix}) and to show that $\iota$ is anti-$T$ equivariant and ``$\mu$-twisted equivariant" with respect to the loop rotation action of $\BC^\times$ (see Corollary~\ref{Cartan_invol_via_loop_rot}). 

Let $\mathfrak{C}\colon G \iso G$ be the Cartan involution of the group $G$ (it interchanges $B$ and $B_-$, and acts as $t \mapsto t^{-1}$ on $T$). The authomorphism $\mathfrak{C}$ induces the authomorphism of $\BC((z^{-1}))$-points of $G$ that we denote by the same symbol $\mathfrak{C} \colon G((z^{-1})) \iso G((z^{-1}))$.
\prop{}\label{cartan_via_matrix}
The Cartan involution $\iota$ in matrix description~(\ref{matrixmatrix}) is given by 
\begin{equation*}
 g \mapsto \mathfrak{C}(g)^{-1},\, g \in U[[z^{-1}]]_1z^{\mu}B_{-}[[z^{-1}]]_1 \cap \overline{G[z]z^{\la}G[z]}.
\end{equation*}
\eprop
\prf
Let us first of all recall that for $\mu_-,\,\mu_+ \in \La$ with $\mu=\mu_++\mu_-$ we can define the variety $\ol{\CW}^\la_{\mu_{-},\mu_+}$ parametrizing tuples $(\CP_{-},\CP_{+},\sigma,\phi_-,\phi_+,s)$, where 
$\CP_-,\,\CP_+$ are $G$-bundles on $\BP^1$,      
$\sigma\colon {\CP_-}|_{\BP^1 \setminus \{0\}} \iso {\CP_+}|_{\BP^1 \setminus \{0\}}$ is an isomorphism having a pole of degree $\leqslant \la$ at $0 \in \BP^1$, $s\colon G \iso {\CP_-}|_{\infty}$ is a trivialization, $\phi_-$ 
is a $B_-$-structure in $\CP_-$ of degree $-w_0\mu_-$ equal to $B$ at $\infty$ w.r.t. $s$, $\phi_+$ is a $B$-structure in $\CP_+$ of degree $w_0\mu_+$ equal to $B_-$ at $\infty$ with respect to $s$, see Section~\ref{symmetric_slices} and~\cite[Section~$2(\text{v})$]{BFN} for details. We will  omit $s$ from the notations when one of the bundles is canonically trivial. The variety $\ol{\CW}^{\la}_{\mu_-,\mu_+}$ is isomorphic to $\ol{\CW}^\la_\mu$ (see Proposition~\ref{iso_symm_ord} and~\cite[Section~$2(\text{v})$]{BFN}). In the same way as for $\ol{\CW}^\la_\mu$ the isomorphism 
\begin{equation*}
\Psi_{\mu_-,\mu_+}\colon \ol{\CW}^\la_{\mu_-,\mu_+} \iso 	U[[z^{-1}]]_1z^{\mu}B_{-}[[z^{-1}]]_1 \cap \overline{G[z]z^{\la}G[z]}
\end{equation*}
can be constructed as follows (see~\cite[Section~$2(\text{xi})$]{BFN} and the proof of Lemma~\ref{action_via_matrix}). Let $\CP_{B}$ be the $B$-bundle of degree $-\mu_-$ corresponding to $\phi_-$ and let $\CP_{B_-}$ be the $B_-$-bundle of degree $\mu_+$ corresponding to $\phi_+$. 
We denote by $\CP^{triv}_{B_-}$ the trivial $B_-$-bundle on $\BP^1$ and by $\CP^{triv}_{B}$ the trivial $B$-bundle on $\BP^1$. We restrict $\CP_{B_-},\,\CP_{B}$ to $\BA^1$ and choose the trivializations of these bundles $\sigma_{B_-}\colon {\CP^{triv}_{B_-}}|_{\BA^1} \iso {\CP^{triv}_{B_-}}|_{\BA^1}$, $\sigma_{B}\colon {\CP^{triv}_{B}}|_{\BA^1} \iso {\CP^{triv}_{B}}|_{\BA^1}$ such that after these identifications $\sigma$ becomes an element of $U[[z^{-1}]]_1 z^\mu B_-[[z^{-1}]]_1$. 
Then we define
\begin{equation*}
\Psi_{\mu_-,\mu_+}((\CP_-,\CP_+,\sigma,\phi_-,\phi_+,s))=\on{Ind}_{B}^G(\sigma_{B}^{-1}) \circ \sigma^{-1} \circ \on{Ind}_{B_-}^G(\sigma_{B_-}) \in G(z).
\end{equation*}

Let us now recall the definition of $\iota$. 
Recall that $\phi^{stand}_-$ is the standard $B_-$-structure in $\CP^{triv}$ 
.
Pick a point 
$
(\CP^{triv},\CP,\sigma,\phi^{stand}_-,\phi) \in \ol{\CW}^{\la}_{0,\mu}=\ol{\CW}^\la_\mu
$
and consider the point 
$
(\mathfrak{C}\CP,\CP^{triv},\mathfrak{C}(\sigma^{-1}),\mathfrak{C}\phi,\phi^{stand}_+) \in \ol{\CW}^{\la}_{0,\mu}.
$
Then use the identification $\Psi_{\mu,0}$ and define 
\begin{equation*}
\iota(x)=\Psi_{\mu,0}((\mathfrak{C}\CP,\CP^{triv},\mathfrak{C}\sigma^{-1},\mathfrak{C}\phi,\phi^{stand}_+)).
\end{equation*}

Let us now finally prove that the involution $\iota$ in matrix description is given by $g \mapsto \mathfrak{C}(g)^{-1}$. We pick a point $(\CP^{triv},\CP,\sigma,\phi_-^{stand},\phi_+,s) \in \ol{\CW}^{\la}_{0,\mu}$ corresponding to $g$ i.e. such that 
\begin{equation*}
g=\on{Ind}_{B}^{G}(\sigma_B^{-1}) \circ	\sigma^{-1} \circ \on{Ind}_{B_-}^G(\sigma_{B_-})
\end{equation*}
for certain trivializations $\sigma_{B},\, \sigma_{B_-}$ of bundles ${\CP_{B}^{triv}}|_{\BA^1},\, {\CP_{B_-}}|_{\BA^1}$ (here $\CP_{B_-}$ is the $B_-$-bundle of degree $\mu$ corresponding to $\phi$).

It remains to note that $\Psi_{\mu,0}((\mathfrak{C}\CP,\CP^{triv},\mathfrak{C}\sigma^{-1},\mathfrak{C}\phi,\phi^{stand}_+))$ is equal to
\begin{equation*}
\on{Ind}^G_{B_-}(\mathfrak{C}\sigma_{B_-}^{-1}) \circ \mathfrak{C}\sigma \circ \on{Ind}_{B_-}^G(\mathfrak{C}\sigma_B)	
\end{equation*}
that is exactly $\mathfrak{C}(g)^{-1}$.
\epr

\cor{}\label{Cartan_invol_via_loop_rot}
$(i)$ The isomorphism $\iota$ is anti-$T$-equivariant: for every $t \in T$ we have 
\begin{equation*}
\iota(tx)=t^{-1}\iota(x)~\forall~x\in \ol{\CW}^\la_\mu.\end{equation*}

$(ii)$ 
Recall the action of $\BC^\times$ on $\ol{\CW}^\la_\mu$ via the loop rotation. Recall also the $T$-action on  $\ol{\CW}^\la_\mu$.
For every $s \in \BC^\times$ we have 
\begin{equation*}
\iota(sx)= \mu(s) s \iota(x). 
\end{equation*}
\ecor
\prf
Let $g \in G(z)$ be an element corresponding to $x$ via the matrix decription.
Using Lemma~\ref{action_via_matrix} and Proposition~\ref{cartan_via_matrix} we get
\begin{equation*}
\iota(tx)=\mathfrak{C}(tgt^{-1})^{-1}=(t^{-1}\mathfrak{C}(g)t)^{-1}=t^{-1}\mathfrak{C}(g)^{-1}t=t^{-1}\iota(x).
\end{equation*}

\begin{equation*}
\iota(sx)=\mathfrak{C}(g(s^{-1}z)s^{\mu})^{-1}=s^{\mu}\mathfrak{C}(g(s^{-1}z))^{-1}=s^{\mu} \big(\mathfrak{C}(g(s^{-1}z))^{-1}s^{\mu}\big) s^{-\mu}=\mu(s)s \iota(x).\end{equation*}

\epr

\ssec{fixedpoints}{Torus fixed points}\label{fixedpoints}
    \prop{tfixed}
    The set of
    $T$-fixed points $(\ol{\mathcal W}_\mu^\lambda)^T$ consists of one element if $\mu$ is a weight of $V^\la$ (the irreducible representation of the Langlands dual group $G^{\vee}$ with the highest weight $\la$) and is empty otherwise.
    We denote the corresponding fixed point by $z^\mu \in (\ol{\mathcal W}_\mu^\lambda)^T$.
    \eprop
    
    \prf
    \cite[Lemma~2.8]{kr}.
    \epr
    
    \rem{}
    {\em{Note that the notation $z^\mu$ for the $T$-fixed point of $\ol{\CW}^\la_\mu$ is consistent with the matrix description of $\ol{\CW}^\la_\mu$ i.e. after the identification~(\ref{matrixmatrix}) point $z^\mu \in \ol{\CW}^\la_\mu$ corresponds to the point $z^\mu \in U[[z^{-1}]]_1z^{\mu}B_{-}[[z^{-1}]]_1 \cap \overline{G[z]z^{\la}G[z]}$, note also that the image of the point $z^\mu \in \ol{\CW}^\la_\mu$ under the natural  forgetful morphism $p\colon \ol{\CW}^\la_\mu \ra \ol{\on{Gr}}{}^{\la}_G$
    is exactly $z^\mu \in \on{Gr}_G$.}}
    \erem
    
    \prop{loopcontrZ}\cite{BF1}\label{zastava_contract}
    Fix $\al \in \La^{\on{pos}}$, the loop rotation $\Cx\curvearrowright Z^\alpha$ contracts $Z^\alpha$ to a point.
    \eprop
    \prf
    Recall that there is a $\Cx$-equivariant locally closed embedding $Z^\alpha \hookrightarrow \on{Qmaps}^\alpha(\BP^1,\mathcal{B})$ of zastava space to the projective variety parametrizing quasi-maps of degree $\al$ from $\BP^1$ to the flag variety $\mathcal{B}$. It follows that any point $x \in Z^\al$ has a limit in $\on{Qmaps}^\alpha(\BP^1,\mathcal{B})$. Let us denote this limit by $x_0 \in \on{Qmaps}^\alpha(\BP^1,\mathcal{B})$. It follows that the 
    quasi-map $x_0$ equals to $B_-$ at $\infty$. On the other hand, all $\Cx$-fixed points in $\on{Qmaps}^\alpha(\BP^1,\mathcal{B})$ are ``constant maps" and are uniquely determined by their value at $\infty$, so the only available limit for points of $Z^\alpha$ is its fixed point. 
    \epr
    
\defe[Repellents and Attractors]{rep}\label{rep_attr}
    The repellent (resp. attractor) to the (unique) $T$-fixed point $z^{\mu} \in \ol{\CW}^\la_\mu$ is defined as 
    \begin{equation*}
    \mathcal{R}^\la_\mu:=
    \{x \in \ol{\CW}^\la_\mu \,|\,  \underset{t \ra \infty}{\on{lim}}\,2{\rho}(t)\cdot x=z^{\mu} \}~ (\on{resp}.~\mathcal{A}^\la_\mu:=
    \{x \in \ol{\CW}^\la_\mu \,|\,  \underset{t \ra 0}{\on{lim}}\,2{\rho}(t)\cdot x=z^{\mu} \}).
    \end{equation*}
\edefe

\rem{aff_attr_rep}
{{
Note that $\ol{\CW}^{\la}_{\mu}$ is an affine variety, hence, $\CR^{\la}_{\mu},\, 
\CA^{\la}_{\mu}$ are closed affine subvarieties of $\ol{\CW}^\la_\mu$ (see~\cite[Section~1.4.7]{DrGa}).}}
\erem

    \prop{reptogr}\label{reptogr}
    The natural morphism $p\colon\ol{\mathcal W}_\mu^\lambda\to\ol{\on{Gr}}{}^\lambda_G$, being restricted to ${\mathcal R}^\lambda_\mu$, is an isomorphism onto its image, $\on{dim}(\CR^\la_\mu)=\langle\rho^{\vee},\la-\mu\rangle$. 
    \eprop
    
    \prf
    Follows from~\cite[Theorem~3.1(1)]{kr}.
    \epr

        \prop{gen}\label{fixed_gen}
    Recall the $\BC^\times$-action on $\ol{\CW}^\la_\mu$ via the loop rotation. Then $(\ol{\CW}^\la_\mu)^{\BC^\times}=U_-\cdot z^\mu$ if $\mu$ is the weight of $V^\la$ and is empty otherwise.
    \eprop
    
\prf
 Fix a point $x \in (\ol{\CW}^\la_\mu)^{\BC^\times}$. Let us prove that $x \in \CR^\la_\mu$.
    Using matrix representation $\Psi$ we decompose $x=uz^{\mu}b_-$ where $u \in U[[z^{-1}]]_1,\, b_- \in B_-[[z^{-1}]]_1$. 
    By Lemma~\ref{action_via_matrix} $(a)$ an element $t \in \BC^\times$ sends $x$ to $u(tz^{-1})z^{\mu}t^{-\mu}b_-(tz^{-1})t^{\mu}$.
    Recall that $t\cdot x = x$, hence,
    $u(tz^{-1})=u(z)$ for any $t \in \BC^\times$. It follows that $u \in U \cap (U[[z^{-1}]]_1)=\{1\}$ so $x \in \CR^\la_\mu$.
    It follows from~\cite[Theorem~3.1(1)]{kr} that we have a $\BC^\times$-equivariant isomorphism $\CR^\la_\mu \iso T_\mu \cap \ol{\on{Gr}}^\la_G$ where $T_\mu:=U_-(\CK)\cdot z^\mu \subset \on{Gr}_G$. Note that
    \begin{equation}\label{fixed_inter}
    (T_\mu \cap \ol{\on{Gr}}^\la_G)^{\BC^\times}=\bigsqcup_{\nu \leqslant \lambda,\,\nu\in \Lambda^+}T_\mu \cap (G \cdot z^{\nu}).
    \end{equation}
    I claim that $T_\mu \cap (G \cdot z^{\nu})$ is empty if $\nu \notin W\mu$. Indeed recall that we can consider the $\BC^\times$-action on $\on{Gr}_G$ induced by the cocharacter $-2\rho\colon \BC^\times \ra T$. This action contracts $T_\mu$ to the point $z^\mu \in \on{Gr}_G$. Pick now a point $x \in T_\mu \cap (G \cdot z^{\nu})$. Note that $G \cdot z^\nu \subset \on{Gr}_G$ is closed so we must have $z^\mu=\underset{t \ra 0}{\on{lim}}(-2\rho)(t)x \in G \cdot z^\nu$. It remains to note that $z^\mu \in G \cdot z^\nu$ iff $\nu \in W\mu$.
    We then conclude from~(\ref{fixed_inter}) that $(T_\mu \cap \ol{\on{Gr}}^\la_G)^{\BC^\times}=U_- \cdot z^{\mu}$.
\epr
    
    \prop{attr_contr}\label{attr_contr}
    The $\BC^\times$-action on $\CA^\la_\mu$ via 
    loop rotation contracts $\CA^\la_\mu$
    to $z^\mu \in \CA^\la_\mu$.
    \eprop
    
    \prf
    Recall the isomorphism $\Psi$ from~(\ref{matrixmatrix}). It restricts to the isomorphism 
    $\CA^\la_\mu \simeq U[[z^{-1}]]_{1}z^{\mu}\cap \ol{G[z]z^{\la}G[z]}$. By Lemma~\ref{action_via_matrix}, the $\BC^\times$-action via loop rotation sends $u(z^{-1})z^{\mu} \in U[[z^{-1}]]_{1}z^{\mu}$ 
    to $u(tz^{-1})z^{\mu}$, hence, contracts it to $z^\mu$. The desired follows.
    \epr

\prop{}\label{loop_contracts}
The loop rotation action $\BC^\times \curvearrowright \ol{\CW}^\la_\mu$ contracts it to its fixed points
if $\langle\al^{\vee},\mu\rangle \geqslant -1$ for any positive root $\al^{\vee} \in \Delta_+^{\vee}$.  
\eprop
\prf
Recall that we have an identification (see Section~\ref{matrix})
    \begin{equation*}
    \ol{\CW}^{\la}_{\mu} \simeq \Big(U[[z^{-1}]]_{1}z^{\mu}B_{-}[[z^{-1}]]_{1}\cap \ol{G[z]z^{\la}G[z]}\Big)
    \end{equation*}
    and by Lemma~\ref{action_via_matrix} $(a)$ the $\BC^\times$-action in these terms is given by 
    \begin{equation}\label{actact}
    g(z) \mapsto g(t^{-1}z)t^\mu.
    \end{equation}
    Since $\ol{G[z]z^{\la}G[z]} \subset G[z^{\pm 1}]$ is closed it is enough to show that any point $g \in U[[z^{-1}]]_{1}z^{\mu}B_{-}[[z^{-1}]]_{1}$ flows to some point of $U[[z^{-1}]]_{1}z^{\mu}B_{-}[[z^{-1}]]_{1}$ via the action~(\ref{actact}). 
    Indeed we can write $g=uz^\mu hu_-$, where $u \in U[[z^{-1}]]_1,\, h \in T[[z^{-1}]]_1,\, u_- \in U_-[[z^{-1}]]_1$. The action~(\ref{actact}) sends $g$ to the point 
    $u(t^{-1}z)z^{\mu}h(t^{-1}z)t^{-\mu}u_-(t^{-1}z)t^\mu$. It remains to note that from the fact that $u_- \in U_-[[z^{-1}]]_1,\, \langle \al,\mu \rangle \geqslant -1$ for $\al \in \Delta_+$ it immediately follows that the limit $\underset{t \ra 0}{\on{lim}}\,t^{-\mu}u_-(t^{-1}z)t^\mu$ exists and lies in $U_-[[z^{-1}]]_1$ so $\underset{t \ra 0}{\on{lim}}\,t \cdot g \in z^\mu U_-[[z^{-1}]]_1 \subset U[[z^{-1}]]_{1}z^{\mu}B_{-}[[z^{-1}]]_{1}$. This observation finishes the proof.
    
    
\epr

\rem{}\label{loop_opp}
{{
The converse to the claim in the Proposition~\ref{loop_contracts} is also true. One can deduce it from~\cite[Remark~B.20]{BFN} together with~\cite[Theorem~B.18]{BFN}.
}
}
\erem
    
\cor{}\label{contr_to_pt}
Assume that 
$\langle\al^{\vee},\mu\rangle \geqslant -1$ for any $\al^{\vee} \in \Delta_+^{\vee}$. Then there exists a $\BC^\times$-action on $\ol{\CW}^\la_\mu$ which contracts it to the point $z^\mu$.
\ecor
\prf
It follows from Proposition~\ref{loop_contracts} and Proposition~\ref{fixed_gen} that the loop rotation action contracts $\ol{\CW}^\la_\mu$ to $U_- \cdot z^\mu$. 
It remains to construct a $\BC^\times$-action on $\ol{\CW}^\la_\mu$ which contracts $U_- \cdot z^\mu \subset \ol{\CW}^\la_\mu$ to the point $z^\mu$. The action via $-2\rho\colon \BC^\times \ra T$ works since conjugation by $-2\rho(t)$ contracts $U_- \subset G$ to $1 \in G$. So the desired $\BC^\times$-action on $\ol{\CW}^\la_\mu$ will be given by 
$x \mapsto (-2\rho(t))t^dx$ for $d \gg 0$, where $t^d$ acts via the loop rotation and $-2\rho(t) \in T$ acts via the $T$-action on $\ol{\CW}^\la_\mu$.
\epr

\rem{}\label{comb_cor}
{{Note that if $\mu \in \La$ is such that $\langle \al^{\vee},\mu \rangle \geqslant -1$ for any $\al^{\vee} \in \Delta_+^{\vee}$ and $\la \in \La^+$ is a dominant coweight such that $\mu \leqslant \la$ then $\mu$ is a weight of $V^\la$. This follows from Proposition~\ref{loop_contracts} together with Proposition~\ref{fixed_gen}. 
One can also give a ``representation-theoretic" proof of this fact in the same spirit as the proof of the standard fact that if $\mu,\la$ are dominant coweights and $\mu \leqslant \la$ then $V^\la_\mu \neq 0$. See Lemma~\ref{weight_rep_mu} in Appendix for more details. 
}
}
\erem 

\sec{fib}{Generalized slices vs slices in affine Grassmannian}\label{main_result_fibration}
In this section, we relate generalized transversal slices $\overline{\CW}^\la_{\mu}$ for coweights $\la \in \La^+,\, \mu \in \La$ with transversal slices in the affine Grassmannian $\on{Gr}_G$.
We prove that if $\mu$ is such that $\langle \al^{\vee},\mu \rangle \geqslant -1$ for every $\al^{\vee} \in \Delta^{\vee}_{+}$, 
then we have an isomorphism 
\begin{equation*}
\ol{\CW}^\la_\mu \simeq \ol{\CW}^\la_{\mu^+} \times \BA^{\langle2\rho^\vee,\,\mu^+-\mu\rangle}.
\end{equation*} 

Note that by~\cite[Remark~2.9]{BF1} we have 
$
\overline{\CW}^\la_{\mu^+} \simeq G[z^{-1}]_1 \cdot z^\mu \cap \overline{\on{Gr}}{}^\la_G,    
$ which is the {\em{non-generalized}} affine Grassmannian slice.
Consider also the subvariety 
$G[z^{-1}]_{U_-} \cdot z^\mu \cap \overline{\on{Gr}}{}^\la_G \subset \on{Gr}_G$, where $G[z^{-1}]_{U_-} \subset G[z^{-1}]$ is the preimage of $U_-$ with respect to the evaluation at $\infty$ morphism $G[z^{-1}] \ra G$. 
We have the natural action $G \curvearrowright \on{Gr}_G$ and denote by $\on{Stab}_{G}(z^\mu) \supset \on{Stab}_{U_-}(z^\mu)$ the corresponding stabilizers of the point $z^\mu \in \on{Gr}_G$.
Note that 
$\on{Stab}_{G}(z^\mu)=P_\mu$, where $P_\mu$ is the parabolic subgroup of $G$ such that the Lie algebra of $P_\mu$ is generated by the root spaces $\mathfrak{g}_{\al^{\vee}}$ with $\langle \al^{\vee},\mu \rangle \leqslant 0$. It is a standard fact that there exists a subgroup $U_-^\mu \subset U_-$ such that the multiplication morphism
    $
    U^{\mu}_{-} \times \on{Stab}_{U_-}(z^\mu) \ra U_-,\, (u_1,u_2) \mapsto u_1u_2
    $
  is an isomorphism of algebraic varieties. 
  Then we have an isomorphism  
\begin{equation}\label{iso_U_1}
U_-^\mu \times (G[z^{-1}]_{1} \cdot z^\mu \cap \overline{\on{Gr}}{}^\la_G) \iso (G[z^{-1}]_{U_-} \cdot z^\mu) \cap \overline{\on{Gr}}{}^\la_G
\end{equation}
given by 
$(u_-,x) \mapsto u_-x$. Note that the $LHS$ of~(\ref{iso_U_1}) identifies with $U_-^\mu \times \ol{\CW}^\la_{\mu^+}$.

Recall that we have a forgetful  morphism 
$p=p^\la_\mu\colon \ol{\CW}^\la_\mu \ra \ol{\on{Gr}}{}^\la_G$. 
We denote by $A^\la_\mu \subset \ol{\CW}^\la_\mu$ the attractor to the fixed locus with respect to the loop rotation action (do not confuse with $\mathcal A_\mu^\la$), by~\cite[Section~1.4.7]{DrGa} this is a closed affine subvariety of $\ol{\CW}^\la_\mu$.

\prop{fib}\label{fibration} The following holds.
\begin{enumerate}
\item 
We have $A^\la_\mu=p^{-1}(G[z^{-1}]_{U_-}\cdot z^\mu \cap \ol{\on{Gr}}{}^\la_G)$.
\item
The morphism $A^\la_\mu = p^{-1}(G[z^{-1}]_{U_-} \cdot z^\mu \cap \ol{\on{Gr}}{}^\la_G) \xrightarrow{p} (G[z^{-1}]_{U_-} \cdot z^\mu) \cap \ol{\on{Gr}}{}^\la_G$ is a trivial fiber bundle, i.e., we have an isomorphism 
$A^\la_\mu \simeq p^{-1}(z^\mu) \times (G[z^{-1}]_{U_-} \cdot z^\mu \cap \ol{\on{Gr}}{}^\la_G)$ such that $p$ is the projection morphism after this identification.
\item 
We have an isomorphism $A^\la_\mu \simeq   p^{-1}(z^\mu) \times U^\mu_- \times \ol{\CW}^\la_{\mu^+}$.
\end{enumerate}
\eprop
\prf
Let us prove $(1)$. Recall that by Proposition~\ref{fixed_gen} we have $(\ol{\CW}^\la_\mu)^{\BC^\times}=U_- \cdot z^\mu$. Note also that $G[z^{-1}]_{U_-}\cdot z^\mu \cap \ol{\on{Gr}}{}^\la_G$ is exactly the attractor to $U_- \cdot z^\mu$ in $\ol{\on{Gr}}{}^\la_G$ with respect to the loop rotation. As forgetful map $p$ is $\BC^\times$-equivariant it follows that $p(A^\la_\mu) \subset (G[z^{-1}]_{U_-} \cdot z^\mu) \cap \ol{\on{Gr}}{}^\la_G$. It remains to prove that $p^{-1}(G[z^{-1}]_{U_-} \cdot z^\mu \cap \ol{\on{Gr}}{}^\la_G) \subset A^\la_\mu$. Pick a point $x \in p^{-1}(G[z^{-1}]_{U_-} \cdot z^\mu \cap \ol{\on{Gr}}{}^\la_G)$. Recall that we have a locally closed embedding $\imath\colon \ol{\CW}^\la_\mu \hookrightarrow \ol{\on{Gr}}{}^\la_G \times Z^{-w_0(\la-\mu)}$ given by $o \mapsto (p(o),s^{\la}_{\mu}(o))$, here $s^{\la}_{\mu}\colon \ol{\CW}^\la_\mu \ra Z^{-w_0(\la-\mu)}$ is a certain $\BC^\times$-equivariant morphism (see~\cite[Section~2(ii)]{BFN} for the definition).
Set $(s,y):=\imath(x)$.
By Proposition~\ref{zastava_contract} the loop rotation action contracts $Z^{-w_0(\la-\mu)}$ to the unique fixed point which we will denote by $\on{pt} \in Z^{-w_0(\la-\mu)}$. Recall that $\ol{\on{Gr}}{}^\la_G$ is projective, hence, any point of $\ol{\on{Gr}}{}^\la_G$ has a limit with respect to the loop rotation action as $t$ goes to $0$. Set $s_0:=\underset{t \ra 0}{\on{lim}}\,t \cdot s \in \ol{\on{Gr}}{}^\la_G$. 
Recall that $s=p(x) \in G[z^{-1}]_{U_-}\cdot z^\mu$, so we must have $s_0 \in U_- \cdot z^\mu$ and by~\cite[Theorem~3.1~(1)]{kr} the morphism $\ol{\CW}^\la_\mu \supset U_- \cdot z^\mu \xrightarrow{p} U_- \cdot z^\mu \subset \ol{\on{Gr}}{}^\la_G$ is an isomorphism. It follows that there exists the unique $\tilde{s}_0 \in U_- \cdot z^\mu \subset \ol{\CW}^\la_\mu$ such that $p(\tilde{s}_0)=s_0$. Note now that $\tilde{s}_0 \in \ol{\CW}^{\la}_{\mu}$ is $\BC^\times$-fixed, so we must have $s^{\la}_{\mu}(\tilde{s}_0)=\on{pt}$. So we obtain $\imath(\tilde{s}_0)=(s_0,\on{pt})=\underset{t \ra 0}{\on{lim}}\,t \cdot \imath(x)$, hence, (since $\imath$ is an embedding) we must have $\underset{t \ra 0}{\on{lim}}\,t \cdot x=\tilde{s}_0$, i.e., $x \in A^\la_\mu$.  We conlude that $p^{-1}(G[z^{-1}]_{U_-}\cdot z^\mu \cap \ol{\on{Gr}}{}^\la_G) \subset A^\la_\mu$. 

Let us prove $(2)$.
We can consider the ind-scheme $\CW^{\mathrm{rat}}_\mu=\underset{\longrightarrow}{\on{lim}}\,\ol{\CW}^\la_\mu$, parametrizing triples $(\CP,\sigma,\phi)$, where we put no restrictions on the defect of $\sigma$ (in matrix description this is $G(z) \cap U[[z^{-1}]]_1z^\mu B_-[[z^{-1}]]_1$).
Note that $\ol{\CW}^\la_\mu$ is the preimage of $\ol{\on{Gr}}{}^\la_G$ under the natural morphism $p_\mu\colon \CW^{\mathrm{rat}}_\mu \ra \on{Gr}_G$, which we will simply denote by $p$.
Note that we have an action $G[z^{-1}]_{U_-} \curvearrowright \CW_\mu^{\mathrm{rat}}$ via changing the trivialization, indeed the action of $g \in G[z^{-1}]_{U_-}$ corresponding to changing the trivialization $\sigma$ does not change $\CP,\,\phi$ and the fiber at infinity $\phi_\infty$ with respect to the new trivialization is $g(\infty)B_-g(\infty)^{-1}=B_-$ since $g(\infty) \in U_-$. 
The morphism $\CW_\mu^{\mathrm{rat}} \ra \on{Gr}_G$ is clearly $G[z^{-1}]_{U_-}$-equivariant. 
We set $A^{\mathrm{rat}}_\mu:=\underset{\longrightarrow}{\on{lim}}\,A^\la_\mu$, this is the attractor in $\CW_\mu^{\mathrm{rat}}$ with respect to the loop rotation action. 
By $(1)$, the morphism $p\colon \CW_\mu^{\mathrm{rat}} \ra \on{Gr}_G$ restricts to the $G[z^{-1}]_{U_-}$-equivariant morphism $p\colon A_\mu^{\mathrm{rat}} \ra G[z^{-1}]_{U_-} \cdot z^\mu$. 
From $G[z^{-1}]_{U_-}$-equivariance it follows that the fibers of the morphism $A_\mu^{\mathrm{rat}} \ra G[z^{-1}]_{U_-} \cdot z^\mu$ are isomorphic. Our goal for now is to show that this morphism is a {\em{trivial}} fibration (this will imply $(2)$). 

Let ${\bf{Gr}}_G$ be the scheme parametrizing pairs $(\CP,\sigma_{U_\infty})$ consisting of a $G$-bundle $\CP$ on $\BP^1$ together with a trivialization $\sigma_{U_\infty}$ of $\CP$ restricted to the formal neighbourhood $U_\infty$ of $\infty \in \BP^1$. Scheme ${\bf{Gr}}_G$ is called a {\em{thick}} affine Grassmannian (see Section~\ref{thick} for more details).
We have the natural morphism ${\bf{i}}\colon \on{Gr}_G \hookrightarrow {\bf{Gr}}_G$. Recall that $\on{Gr}_G=\underset{\longrightarrow}{\on{lim}}\,\ol{\on{Gr}}{}^\la_G$. The composition $\ol{\on{Gr}}{}^\la_G \subset \on{Gr}_G \subset {\bf{Gr}}_G$ is a closed embedding 
(see
\cite[Section~10.6.2]{FM},~\cite[Proposition~1.3.2]{kt}
) so ${\bf{i}}$ is an embedding of functors, note that the morphism ${\bf{i}}$ has dense image (in Zariski topology) so is not a closed embedding. 
Let $G[[z^{-1}]]_{U_-} \subset G[[z^{-1}]]$ be the preimage of $U_-$ with respect to the natural evaluation at $\infty$ morphism $G[[z^{-1}]] \ra G$. The group $G[[z^{-1}]]$ acts on ${\bf{Gr}}_G$ via changing the trivialization.
Note that $(G[[z^{-1}]]_{U_-} \cdot z^\mu) \cap \on{Gr}_G=G[z^{-1}]_{U_-} \cdot z^\mu$  since we have decompositions
$
\on{Gr}_G=\bigsqcup_{\nu \in \La}G[z^{-1}]_{U_-} \cdot z^\nu,\,    \on{{\bf{Gr}}}_G=\bigsqcup_{\nu \in \La}G[[z^{-1}]]_{U_-} \cdot z^\nu 
$ (see~\cite[Propositions~1.3.1,~1.3.2]{kt}).
It follows that ${\bf{i}}$ restricts to the embedding (that  becomes closed after intersecting with $\ol{\on{Gr}}{}^\la_G$):
\begin{equation*}
\on{Gr}_G \supset G[z^{-1}]_{U_- } \cdot z^\mu \hookrightarrow G[[z^{-1}]]_{U_-} \cdot z^\mu \subset {\bf{Gr}}_G
\end{equation*}
to be denoted by the same symbol. 
Let $\CW_\mu$ be the moduli space of triples $(\CP,\sigma_{U_\infty},\phi)$, where $\CP$ is a $G$-bundle on $\BP^1$, $\sigma_{U_\infty}$ is the trivialization of $\CP$ at the formal neighbourhood $U_\infty$ of $\infty \in \BP^1$ and $\phi$ is a $B$-structure of degree $w_0\mu$ in $\CP$ such that the fiber of $\phi$ at $\infty$ is $B_-$ with respect to $\sigma_\infty$ (in matrix description we have $\CW_\mu=U[[z^{-1}]]_1 z^\mu B_-[[z^{-1}]]_1$). We have the natural morphism ${\bf{p}}\colon \CW_\mu \ra {\bf{Gr}}_G$ and the following diagram is cartesian (c.f.~\cite[Section~12.2]{FM}, this follows from the moduli description of $\CW^{\mathrm{rat}}_\mu,\, \CW_\mu,\,\on{Gr}_G,\, {\bf{Gr}}_G$):
\begin{equation*}
{\xymatrix{
\CW_\mu^{\mathrm{rat}} \ar[r] \ar[d]^p & \CW_\mu \ar[d]^{{\bf{p}}} \\
\on{Gr}_G \ar[r]^{{\bf{i}}} & {\bf{Gr}}_G
}
}
\end{equation*}
Set $A_\mu:={\bf{p}}^{-1}(G[[z^{-1}]]_{U_-}z^\mu)$. Note that by $(1)$ the following diagram is cartesian:
\begin{equation*}
{\xymatrix{
A_\mu^{\mathrm{rat}} \ar[r] \ar[d]^p & A_\mu \ar[d]^{{\bf{p}}} \\
G[z^{-1}]_{U_-} \cdot z^\mu \ar[r]^{{\bf{i}}} & G[[z^{-1}]]_{U_-} \cdot z^\mu
}
}
\end{equation*}


Note that $G[[z^{-1}]]_{U_-}$ is a pro-unipotent group, so it is easy to construct a subgroup $\Gamma \subset G[[z^{-1}]]_{U_-}$ such that the action of $\Gamma$ on $G[[z^{-1}]]_{U_-} \cdot z^\mu$ is free and transitive.
Indeed the Lie algebra of $G[[z^{-1}]]_{U_-}$ is $\mathfrak{u}_- \oplus z^{-1}\mathfrak{g}[[z^{-1}]]$ (here $\mathfrak{u}_-=\on{Lie}U_-$) and the Lie algebra of $\on{Stab}_{G[[z^{-1}]]_{U_-}}(z^\mu)$ is generated (topologically i.e. allowing infinite sums) by $z^{-k}\mathfrak{g}_{\al^\vee}$ such that $\langle\al,\mu\rangle+k \leqslant 0$. We can now consider the Lie subalgebra $\mathfrak{a} \subset \mathfrak{u}_- \oplus z^{-1}\mathfrak{g}[[z^{-1}]]$ generated (topologically) by $z^{-k}\mathfrak{g}_{\al^\vee}$ such that $\langle\al^\vee,\mu\rangle+k>0$, note that this is indeed the Lie subalgebra of $\mathfrak{u}_- \oplus z^{-1}\mathfrak{g}[[z^{-1}]]$ since $[z^{-k_1}\mathfrak{g}_{\al^\vee_1},z^{-k_2}\mathfrak{g}_{\al^\vee_2}] \subset z^{-k_1-k_2}\mathfrak{g}_{\al^\vee_1+\al^\vee_2}$ and $\langle\al^\vee_1,\mu\rangle+k_1>0,\, \langle\al^\vee_2,\mu\rangle+k_2>0$ implies $\langle\al^\vee_1+\al^\vee_2,\mu\rangle+k_1+k_2>0$. Now we can define $\Gamma$ as the exponent of $\mathfrak{a}$. Note that we can explicitly write  (compare with~\cite[Lemma~4.5.7 and Corollary~4.5.8]{kash})
\begin{equation*}
\Gamma
=G[[z^{-1}]]_{U_-} \cap z^\mu G[[z^{-1}]]_1z^{-\mu}.
\end{equation*}
Since ${\bf{p}}$ is $\Gamma$-equivariant we conclude that the action by $\Gamma$ induces the isomorphism $A_\mu \simeq {\bf{p}}^{-1}(z^\mu) \times (G[[z^{-1}]]_{U_-} \cdot z^\mu)$ such that the morphism ${\bf{p}}\colon A_\mu \ra G[[z^{-1}]]_{U_-} \cdot z^\mu$ is the projection onto the second factor with respect to the identification above.
So we see that we have the following {\em{cartesian}} square: 
\begin{equation*}
{\xymatrix{
A_\mu^{\mathrm{rat}} \ar[r] \ar[d]^p & {p}^{-1}(z^\mu) \times (G[[z^{-1}]]_{U_-} \cdot z^\mu) \ar[d]^{\on{proj}} \\
G[z^{-1}]_{U_-} \cdot z^\mu \ar[r]^{{\bf{i}}} & G[[z^{-1}]]_{U_-} \cdot z^\mu
}
},    
\end{equation*}
here $\on{proj}\colon {p}^{-1}(z^\mu) \times (G[[z^{-1}]]_{U_-} \cdot z^\mu) \ra G[[z^{-1}]]_{U_-} \cdot z^\mu$ is the projection onto the second factor. Since the square is cartesian we conclude that $A_\mu^{\mathrm{rat}} \simeq {p}^{-1}(z^\mu) \times (G[z^{-1}]_{U_-} \cdot z^\mu)$.
Another way to prove this is to consider the group $\Gamma^{\mathrm{rat}}:=\Gamma \cap G(z)$ and to note that its action on $G[z^{-1}]_{U_-}\cdot z^\mu \subset \on{Gr}_G$ is free and transitive. It follows that the action of $\Gamma^{\mathrm{rat}}$ induces the desired isomorphism $A^{\mathrm{rat}}_\mu \simeq p^{-1}(z^\mu) \times (G[z^{-1}]_{U_-} \cdot z^\mu)$.
Now $(2)$ follows.


Part $(3)$ follows from $(2)$ using the isomorphism~(\ref{iso_U_1}).
\epr

\rem{}
{{
Using the same approach as in the proof of Theorem~(\ref{fibration}), one can show that for any $\nu \in \La$ the morphism 
\begin{equation*}
\ol{\CW}^\la_\mu \supset p^{-1}(G[z^{-1}]_{U_-} \cdot z^\nu \cap \ol{\on{Gr}}{}^\la_G) \ra G[z^{-1}]_{U_-} \cdot z^\nu \cap \ol{\on{Gr}}{}^\la_G \subset \ol{\on{Gr}}{}^\la_G
\end{equation*}
is a trivial fibration, so 
\begin{equation*}
p^{-1}(G[z^{-1}]_{U_-} \cdot z^\nu \cap \ol{\on{Gr}}{}^\la_G) \simeq p^{-1}(z^\nu) \times (G[z^{-1}]_{U_-} \cdot z^\nu \cap \ol{\on{Gr}}{}^\la_G).
\end{equation*}
Since $\ol{\on{Gr}}{}^\la_G=\bigsqcup_{\nu \in \La} G[z^{-1}]_{U_-} \cdot z^\nu \cap \ol{\on{Gr}}{}^\la_G$, it follows that $\ol{\CW}^\la_\mu$ has the following stratification by locally closed subvarieties:
\begin{equation*}
\ol{\CW}^\la_\mu=\bigsqcup_{\nu \in \La} p^{-1}(z^\nu) \times (G[z^{-1}]_{U_-} \cdot z^\nu \cap \ol{\on{Gr}}{}^\la_G)=\bigsqcup_{\nu \in \La} p^{-1}(z^\nu) \times U^\nu_- \times \ol{\CW}^\la_{\nu^+}.
\end{equation*}
}
}
\erem

\rem{}
{{
Directly from the definitions $p^{-1}(z^\nu)$ is the moduli space of $B$-structures of degree $w_0\mu$ in the $G$-bundle $\CO(\nu)$ which are equal to $B_-$ at $\infty$ with respect to the trivialization $\sigma_{z^\nu}$. 
It would be interesting to describe this space, for example, to compute its dimension. 
}}
\erem

\th{}\label{affine_fibr}
For $\mu \in \La$ such that 
$\langle\al^{\vee},\mu\rangle \geqslant -1$ for $\al^{\vee} \in \Delta^{\vee}_+$, the image of the morphism $p$ coincides with $(G[z^{-1}]_{U_-} \cdot z^\mu) \cap \ol{\on{Gr}}{}^\la_G$
and 
\begin{equation*}
\ol{\CW}^\la_\mu \simeq 
\BA^{\langle2\rho^{\vee},\,\mu^+-\mu \rangle} \times \ol{\CW}^\la_{\mu^+}.    
\end{equation*}
\eth
\prf
It follows from Proposition~\ref{loop_contracts} that $\ol{\CW}^\la_\mu=A^\la_\mu$, so by Proposition~\ref{fibration} the image of the morphism $p$ coincides with $(G[z^{-1}]_{U_-} \cdot z^\mu) \cap \ol{\on{Gr}}{}^\la_G$ and using the isomorphism~(\ref{iso_U_1}) we obtain 
\begin{equation*}
\ol{\CW}^\la_\mu=A^\la_\mu \simeq p^{-1}(z^\mu) \times (G[z^{-1}]_{U_-} \cdot z^\mu \cap \ol{\on{Gr}}{}^\la_G) \simeq p^{-1}(z^\mu) \times U^\mu_- \times \ol{\CW}^\la_{\mu^+}.    
\end{equation*}
Since $U^\mu_-$ is isomorphic to an affine space, it remains to prove that $p^{-1}(z^\mu)$ is also isomorphic to an affine space. Let us first of all show that $p^{-1}(z^\mu)$ is smooth.

Recall that by Proposition~\ref{fibration} the morphism $p\colon \overline{\CW}^\la_\mu=A^\la_\mu \ra (G[z^{-1}]_{U_-} \cdot z^\mu) \cap \ol{\on{Gr}}{}^\la_G$ is a trivial fibration, so all its fibers are isomorphic to $p^{-1}(z^\mu)$. Note that the preimage $p^{-1}(G[z^{-1}]_{U_-} \cdot z^\mu \cap \on{Gr}^\la_G)$ is $\CW^\la_\mu$ that is smooth by~\cite{smooth_slice}. So we see that the preimage 
$p^{-1}(G[z^{-1}]_{U_-} \cdot z^\mu \cap \on{Gr}^\la_G) \simeq p^{-1}(z^{\mu}) \times (G[z^{-1}]_{U_-} \cdot z^\mu \cap \on{Gr}^\la_G)$  is smooth, hence
$p^{-1}(z^\mu)$ is also smooth. 

Recall that the loop rotation action contracts $p^{-1}(z^\mu)$ to the fixed point $z^\mu \in \ol{\CW}^\la_\mu$ and $p^{-1}(z^\mu)$ is smooth. Now  from  (very simple version of) the Bialynicki-Birula decomposition (see~\cite{BB}) we deduce that $p^{-1}(z^\mu)$ is isomorphic to an affine space.
\epr

\rem{}
{{
Let us note that from our proof of Theorem~\ref{affine_fibr} it follows that the Conjecture a) in Section~12.3.1 of \cite{FM} that $p\colon \ol{\CW}^\la_\mu \ra \ol{\on{Gr}}{}^\la_G$ is smooth onto its image is true for $\mu$ such that $\langle\al^{\vee},\mu\rangle \geqslant -1$ for $\al^{\vee} \in \Delta_+^{\vee}$. 
Let us also note that the Conjecture is false in general: indeed for $G=\on{SL}_2$ and $\mu=-\al$, $\la=\al$ (with $\al$ being the simple coroot of $\on{SL}_2$) one can check that $1 \in \on{im}p$.  Using the moduli description of $\ol{\CW}^\la_\mu$ it is easy to see that $p^{-1}(1)$ coincides with open zastava space $\oZal$, so it has dimension $2$. On the other hand, $\on{dim}\ol{\CW}^\al_{-\al}=4$ and $\on{dim}G[z^{-1}]_{U_-} \cdot z^{-\al} \cap \ol{\on{Gr}}{}^\al_G=\on{dim}U_- \cdot z^{-\al}=\on{dim}U_-=1$, i.e., $\on{dim}p^{-1}(z^{-\al})=3$. We see that $p$ has fibers of different (positive) dimensions over the points $1,\, z^{-\al} \in \on{Gr}_G$, so $p$ can not be smooth onto its image.  We are gratefull to Dinakar Muthiah for pointing out this example to us. 
}
}
\erem

\cor{}\label{more_general}
If $\mu \in W\la$ and $\langle\al^{\vee},\mu\rangle \geqslant -1$ for every $\al^{\vee} \in \Delta^{\vee}_+$, then $\ol{\CW}^\la_\mu \simeq \BA^{2\langle \rho^\vee,\,\la-\mu \rangle}$.
\ecor
\prf
Follows from Theorem~\ref{affine_fibr}.
\epr

Let us give another short proof of this fact.
\prf
Recall that by Corollary~\ref{contr_to_pt} we have a $\BC^\times$-action that contracts $\ol{\CW}^\la_\mu$ to the point $z^\mu$, note also that $z^\mu$ lies in a smooth open subset $\CW^\la_\mu \subset \ol{\CW}^\la_\mu$ (since $\mu \in W\la$), so we conclude that $\ol{\CW}^\la_\mu$ is smooth at the point $z^\mu$ (smoothness of $\CW^\la_\mu$ follows from~\cite{smooth_slice}). It follows that $\ol{\CW}^\la_\mu$ is smooth (here we use that if we have a $\BC^\times$-action on some variety $Z$ that contracts $Z$ to a point and $Z$ is smooth at this point, then $Z$ is smooth). 
Now the claim follows from the Bialynicki-Birula decomposition (see~\cite{BB}).
\epr


\rem{}\label{Dinak_exp}
{{
Note that the condition that $\langle \al^{\vee},\mu \rangle \geqslant -1$ for any positive $\al^{\vee} \in \Delta^{\vee}_+$ is equivalent to the fact that $\ol{\CW}^{\mu^+}_\mu=\CW^{\mu^+}_\mu$. 
Equivalently $\CW^{\mu^+}_{\mu} \subset \ol{\CW}^{\la}_{\mu}=\bigsqcup_{\mu \leqslant \la' \leqslant \la}\CW^{\la'}_{\mu}$ is the {\em{deepest stratum}} (i.e. for any dominant $\la' \in \La^+$ such that $\mu \leqslant \la' \leqslant \la$ we have $\CW^{\mu^+}_{\mu} \subset \ol{\CW}^{\la'}_{\mu}$) iff $\langle \al^{\vee},\mu \rangle \geqslant -1$ for any positive $\al^{\vee} \in \Delta_+^{\vee}$. This easily follows from the results of Appendix (see Corollary~\ref{equiv_deepest_comb}). We  are  grateful  to   Dinakar  Muthiah  for pointing this out to us.
}
}
\erem

\section{Generalized slices for minuscule coweights}\label{minuscule_case}
In this section we consider the example of generalized slices $\ol{\CW}^\la_\mu$ for minuscule $\la$ and $\mu \in W\la$. Since we will use the Coulomb branch description of $\ol{\CW}^\la_\mu$ we then assume in this section that $G$ is adjoint. 

\rem{}
{\em{Note that we do not loose any generality assuming that $G$ is adjoint: indeed if $G$ is any reductive group with Lie algebra $\mathfrak{g}$ and $G^{\mathrm{ad}}$ is the corresponding adjoint group then we have a surjective homomorphism $G \twoheadrightarrow G^{\mathrm{ad}}$. Let $\ol{T}$ be the image of $T$ with respect to this morphism. Note that $\ol{T}$ is the maximal torus of $G^{\mathrm{ad}}$. Then any two cocharachters $\la \in \La^{+},\,\mu \in \La$ being composed with the surjection $T \twoheadrightarrow \ol{T}$ define cocharacters of $\ol{T}$ to be denoted $\ol{\la},\,\ol{\mu}$ respectively. It is easy to see that the canonical morphism $\ol{\CW}^{\la}_{\mu} \iso \ol{\CW}^{\ol{\la}}_{\ol{\mu}}$ is an isomorphism i.e. the generalized transversal slice $\ol{\CW}^{\la}_{\mu}$ for $G$ can be realized as the generalized transversal slice $\ol{\CW}^{\ol{\la}}_{\ol{\mu}}$ for the corresponding adjoint group $G^{\mathrm{ad}}$.   
}}
\erem

Let us recall the definition of a minuscule coweight.
     \defe[Minuscule coweights]{minisc}
    A dominant coweight $\lambda \in \Lambda^+$ is called minuscule if for any coweight $\mu \in \La$ such that $V^\lambda_{\mu}\neq \{0\}$ we have $\mu \in W\la$.
    Here $V^\lambda$ is the irreducible representation of the Langlands dual group $G^{\vee}$ with the highest weight $\la$
    and $V^\la_\mu$ is the $\mu$-weight space of $V^\la$.
    \edefe

Note that since $\la$ is minuscule it follows that $\langle\al^{\vee},\la\rangle \in \{-1,0,1\}$ for any $\al^{\vee} \in \Delta^{\vee}$ (see~~\cite[VIII, \textsection 7, no. 3]{Bu}) so since $\mu \in W\la$ we conclude that $\langle\al^{\vee},\mu\rangle \geqslant -1$ for any $\al^{\vee} \in \Delta^{\vee}$ so we can apply Theorem~\ref{affine_fibr} (see also Corollary~\ref{more_general}) and conclude that $\ol{\CW}^\la_\mu \simeq \BA^{\langle 2\rho^\vee,\,\la-\mu \rangle}$. The goal of this section is to construct the isomorphism $\ol{\CW}^\la_\mu \simeq \BA^{\langle 2\rho^\vee,\, \la-\mu \rangle}$ explicitly and in particular obtain natural coordinates on $\ol{\CW}^\la_\mu$. We will also describe the Poisson structure on $\ol{\CW}^\la_\mu$ which comes from its realization as a certain Coulomb branch (see~Section~\ref{coulomb}). We will also compute a $T\times \BC^\times$-character of $\ol{\CW}^\la_\mu$ that will allow us to compute characters of tangent spaces at fixed points of convolution diagrams over slices (see Remark~\ref{characters_resolv}) and also to compute Poincar\'e polynomials of certain convolution diagrams $\widetilde{\CW}^{\ul{\la}}_{\mu}$ (see Remark~\ref{Poincare_pol}). Note that since $\mu \in W\la$ and $\la \in \La^+$ we then have $\la=\mu^+$. 

\ssec{}{Coordinates on $\ol{\CW}^{\la}_\mu=\CW^{\mu^+}_\mu$ for minuscule $\la$}
In this section we assume that $\mu \in W\la$ and $\la$ is minuscule. It follows that $\la=\mu^+$ and
 $\ol{\CW}^\la_\mu=\ol{\CW}^{\mu^+}_\mu=\CW^{\mu^+}_\mu$ (see Remark~\ref{Dinak_exp}). 

Recall that (see Definition~\ref{rep_attr})
\begin{equation*}
    \mathcal{R}^{\mu^+}_\mu:=
    \{x \in \CW^{\mu^+}_\mu \,|\,  \underset{t \ra \infty}{\on{lim}}\,2{\rho}(t)\cdot x=z^{\mu} \}.
\end{equation*}
 \lem{cfixed}\label{cfixed}
    For $\mu \in \La$ such that $\mu^+ \in \La^+$ is minuscule we have $(\mathcal W_\mu^{\mu^+})^\Cx=\CR_\mu^{\mu^+}$.
    \elem    
    \prf
    It follows from Proposition~\ref{fixed_gen} that $(\mathcal W_\mu^{\mu^+})^\Cx=U_- \cdot z^\mu \subset \CR^{\mu^+}_\mu$. It remains to prove that $\CR^{\mu^+}_\mu \subset (\mathcal {W}_\mu^{\mu^+})^\Cx$.
    Indeed recall that by Proposition~\ref{reptogr} the morphism $p|_{\CR^{\mu^+}_\mu}\colon \CR^{\mu^+}_\mu \ra \on{Gr}^{\mu^+}_\mu$ is the isomorphism onto its image. 
    Since $\mu^+$ is minuscule we have $\on{Gr}^{\mu^+}_G=G \cdot z^{\mu^+}$ so $\BC^\times$ acts trivially on $\on{Gr}^{\mu^+}_G$. Since $p|_{\CR^{\mu^+}_\mu}\colon \CR^{\mu^+}_\mu \ra \on{Gr}^{\mu^+}_\mu$ is a $\BC^\times$-equivariant embedding we conclude that $\CR^{\mu^+}_\mu \subset (\mathcal {W}_\mu^{\mu^+})^\Cx$.
    \epr
    
    \lem{image}\label{image}
    The image of the  morphism $p^{\mu^+}_\mu=p\colon\CW^{\mu^+}_\mu \ra \on{Gr}^{\mu^+}_G$
    coincides with $p(\CR^{\mu^+}_\mu) \simeq \CR^{\mu^+}_\mu$.
    \elem
    
    \prf
    Take a point $x \in \CW^{\mu^+}_\mu$. By Proposition~\ref{loop_contracts} together with Lemma~\ref{cfixed}, it flows to some point $x_0 \in \CR^{\mu^+}_\mu$ under the loop rotation action. It follows that $p(x)$ flows to $p(x_0)$ under the loop rotation action. Recall now that since $\mu^+$ is minuscule we have $\on{Gr}^{\mu^+}_G=G \cdot z^{\mu^+}$ so
    $\BC^\times$ acts trivially on $\on{Gr}^{\mu^+}_G$, hence, 
    $p(x)=p(x_0) \in p(\CR^{\mu^+}_\mu)$.
    The isomorphism 
    $p(\CR^{\mu^+}_\mu) \simeq \CR^{\mu^+}_\mu$ follows 
    from Proposition~\ref{reptogr}.
    \epr

    Recall that $\on{Stab}_{U_-}(z^\mu) \subset U_-$ is the stabilizer of the point $z^\mu \in \on{Gr}^{\mu^+}_G$
    in $U_-$. Recall that $\on{Gr}^{\mu^+}_G \simeq G/P_{\mu^+}$ is isomorphic to a parabolic flag variety and the corresponding action $U_- \curvearrowright G/P_{\mu^+}$ is given by the left multiplication. 
    There exists a canonical subgroup $U^{\mu}_{-} \subset U_-$, such that the multiplication morphism
    \begin{equation} \label{multip_U_-}
    U^{\mu}_{-} \times \on{Stab}_{z^\mu}(U_-) \ra U_-,\, (u_1,u_2) \mapsto u_1u_2
    \end{equation} 
    is an isomorphism of algebraic varieties. The group $U^{\mu}_-$ is defined as follows. Consider the nilpotent subalgebra of $\mathfrak{g}$ generated by $\mathfrak{g}_{\al^{\vee}}$ such that $\al^{\vee} \in \Delta_{-},\, \langle\al^{\vee},\mu \rangle>0$. Then $U^\mu_-$ is the exponent of this algebra. 

   \lem{free_act}\label{free_act}
    The morphism $\Phi\colon U^\mu_-\times p^{-1}(z^\mu) \iso \CW^{\mu^+}_\mu$ given by $(u_-,x) \mapsto u_- \cdot x$ is an isomorphism.
    \elem
    
    \prf
    It follows from \cref{image} and (\ref{multip_U_-}) that the group $U^\mu_-$ acts freely and transitively on the image of the morphism $p$. The desired follows.
    \epr

    \prop{aff_fib}\label{aff_fib}
    We have $p^{-1}(z^{\mu})=\CA^{\mu^+}_\mu$.
    \eprop
    \prf
    It follows from \cref{attr_contr} that $\CA^{\mu^+}_\mu \subset p^{-1}(z^{\mu})$. Note now that $\on{dim}(\CA^{\mu^+}_\mu)=\on{dim}(\CR^{\mu^+}_\mu)=\langle\rho^{\vee},\la-\mu\rangle=\on{dim}(p^{-1}(z^\mu))$ and $\CA^{\mu^+}_\mu$ is closed in $p^{-1}(z^{\mu})$. Note also that the variety $p^{-1}(z^{\mu})$ is irreducible since by \cref{irred}, $\CW^{\mu^+}_\mu$ is irreducible and by \cref{free_act}, $\CW^{\mu^+}_\mu \simeq U^\mu_- \times p^{-1}(z^\mu)$. Thus, $\CA^{\mu^+}_\mu = p^{-1}(z^\mu)$.
    \epr
    
\cor{}\label{iso_U_-_A} 
We have a natural isomorphism $\CW^{\mu^+}_\mu \simeq U^\mu_- \times \CA^{\mu^+}_\mu \simeq \CR^{\mu^+}_\mu \times \CA^{\mu^+}_\mu$.
\ecor
\prf
Follows from Lemma~\ref{free_act} and Proposition~\ref{aff_fib}.
\epr

    Recall now that the Cartan involution $\iota\colon \CW^{\mu^+}_\mu \iso \CW^{\mu^+}_\mu$  (see~\cite[Section~2(vii)]{BFN}) and Corollary~\ref{Cartan_invol_via_loop_rot}) induces an anti-$T$-equivariant isomorphism from $\CR^{\mu^+}_\mu$ to $\CA^{\mu^+}_\mu$ and that we have the $T$-equivariant isomorphism from $\CR^{\mu^+}_\mu$ to a certain Bruhat cell in $G/P_{\mu^+}$, which is naturally isomorphic to $U^{\mu}_-$. 
 Note that we have an isomorphism $\on{Lie}(U^\mu_-) \overset{\text{exp}}{\simeq}U^\mu_-$ and $\on{Lie}(U^\mu_-) \subset \on{Lie}U_-$ has the basis $\{e_{\al^{\vee}} \},\, e_{\al^{\vee}} \in \mathfrak{g}_{\al^\vee}$ parametrized by the subset $\Delta^{\vee}_{\mu,-} \subset \Delta_{-}^\vee$ defined as follows $\Delta^{\vee}_{\mu,-}:=\{\al^\vee \in \Delta_-^{\vee},\,\langle\al^\vee,\mu\rangle > 0 \}$.
 So we obtain the identifications
 \begin{equation*}
 \CA^{\mu^+}_\mu \overset{\iota}{\simeq} \CR^{\mu^+}_\mu \overset{p}{\simeq} U^\mu_- \overset{\text{log}}{\simeq} \on{Lie}(U^\mu_-) = \BA^{\langle \rho^{\vee},\,\mu^+-\mu \rangle}.
 \end{equation*}
 For $\al^{\vee} \in \Delta^{\vee}_{\mu,-}$ we denote by $y_{\al^{\vee}} \in \BC[\CR^{\mu^+}_\mu]$ the corresponding coordinate function on $\CR^{\mu^+}_\mu$. For $\al^{\vee} \in \Delta^{\vee}_{\mu,+}:=-\Delta^{\vee}_{\mu,-}$, we denote by $x_{\al^{\vee}} \in \BC[\CA^{\mu^+}_\mu]$ the function on $\CA^{\mu^+}_\mu$ which corresponds to the coordinate function $y_{-\al^{\vee}} \in \BC[\CR^{\mu^+}_\mu]$ with respect to the isomorphism $\iota|_{\CR^{\mu^+}_\mu}\colon \CR^{\mu^+}_\mu \iso \CA^{\mu^+}_\mu$.
     So we obtain natural coordinates $x_{\al^{\vee}},\, y_{-\al^{\vee}}$ on $\CA^{\mu^+}_\mu,\, \CR^{\mu^+}_\mu$ that we will call {\em{standard coordinates}}.
The following theorem follows immediately from Corollary~\ref{iso_U_-_A} and the discussion above.
\th{}\label{coordinates_and_isomorphisms}
We have natural isomorphisms 
\begin{equation*}
\CA^{\mu^+}_\mu \overset{\iota}{\simeq} \CR^{\mu^+}_\mu \overset{p}{\simeq} U^\mu_- \overset{\on{log}}{\simeq} \on{Lie}(U^\mu_-) = \BA^{\langle\rho^{\vee},\la-\mu\rangle},\, \CW^{\mu^+}_\mu \simeq \CA^{\mu^+}_\mu \times \CR^{\mu^+}_\mu \simeq \BA^{\langle\rho^{\vee},\mu^+-\mu\rangle} \times \BA^{\langle\rho^{\vee},\mu^+-\mu\rangle}
\end{equation*}
so we obtain natural coordinates $\{x_{\al^{\vee}},\,y_{-\al^{\vee}}\,|\, \al^{\vee} \in \Delta^{\vee}_{\mu,-}\}$ on $\CW^{\mu^+}_\mu$.
\eth

\rem{}\label{matrix_descr_A_R_min}
{\em{Let us now describe coordinates $x_{\al^{\vee}},\, y_{-\al^\vee},\, \al^\vee \in \Delta_{\mu,+}^\vee$ more explicitly using matrix descriptions of $\CA^{\mu^+}_\mu,\, \CR^{\mu^+}_\mu$. Let us first of all recall that the action of $U^\mu_-$ on $\CR^{\mu^+}_\mu$ is free. I claim that in matrix terms this action is given by $u_- \mapsto u_-z^\mu,\, u_- \in U^\mu_-$. To see this recall that by Remark~\ref{act_U-_matrix} the action of $g \in B_-$ on $x \in \CW^{\mu^+}_\mu \subset \CW_\mu$ is given by $x \mapsto p(g i(x) g^{-1})$, here $i\colon \CW_\mu \hookrightarrow X_\mu$ is the embedding and $p\colon X_\mu \twoheadrightarrow \CW_\mu$ is the projection. 
So we need to compute  $p(u_-z^\mu u_-^{-1})$.  We can write $u_-=\on{exp}(n_-)$ for some $n_- \in \on{Lie}U^\mu_-=\bigoplus_{\beta^\vee \in \Delta_{\mu,-}^\vee}\mathfrak{g}_{\beta^\vee}$ then since $\langle \beta^{\vee},\mu \rangle=1$ for $\beta^\vee \in \Delta_{\mu,-}^\vee$ (here we use that $\mu^+$ is minuscule) we get 
\begin{equation*}
u_-z^\mu u_-^{-1}=\on{exp}(n_-)z^{\mu} \on{exp}(-n_-)=z^{\mu}\on{exp}(z^{-1}n_-)\on{exp}(-n_-) \in X_\mu.
\end{equation*}
Projecting to $\CW_\mu$ we obtain 
\begin{equation*}
p(z^{\mu}\on{exp}(z^{-1}n_-)\on{exp}(-n_-))=z^{\mu}\on{exp}(z^{-1}n_-)=\on{exp}(n_-)z^\mu=u_-z^\mu \in \CW^{\mu^+}_\mu.
\end{equation*}
We conclude that 
\begin{equation*}
\CR^{\mu^+}_\mu=\{u_-z^{\mu}\,|\, u_- \in U^\mu_-\}=\Big\{\on{exp}\Big(\sum_{\beta^{\vee} \in \Delta_{\mu,-}^\vee} b_{\beta^\vee}e_{\beta^\vee}\Big)z^\mu,\, b_{\beta^\vee} \in \BC\Big\}
\end{equation*}
and the coordinate $y_{-\al^\vee}$ is given by 
\begin{equation*}
\on{exp}\Big(\sum_{\beta^{\vee} \in \Delta_{\mu,-}^\vee} b_{\beta^\vee}e_{\beta^\vee}\Big)z^\mu \mapsto b_{-\al^\vee}.   
\end{equation*}
Let us now describe the coordinates $x_{\al^\vee}$ in matrix terms. We denote by $U^\mu \subset B$ the exponent of the nilpotent algebra $\bigoplus_{\beta^\vee \in \Delta^\vee_{\mu,+}} \mathfrak{g}_{\beta^\vee}$.
 Recall that by Proposition~\ref{cartan_via_matrix} below the Cartan involution $\iota$ is given by $g \mapsto \mathfrak{C}(g)^{-1}$, here $\mathfrak{C}\colon G \iso G$ is the Cartan involution of the group $G$ and we denote by the same symbol the corresponding isomorphism $G((z^{-1})) \iso G((z^{-1}))$. Since $\mathfrak{C}(u_-z^\mu)^{-1}=z^{\mu}\mathfrak{C}(u_-)^{-1}$ and $\iota$ identifies $\CR^{\mu^+}_\mu$ with $\CA^{\mu^+}_\mu$ we conclude that $\CA^{\mu^+}_\mu$ in matrix terms can be described as follows 
 \begin{equation*}
\CA^{\mu^+}_\mu=\{z^{\mu}u\,|\, u \in U^\mu\}=\Big\{z^\mu\on{exp}\Big(\sum_{\al^{\vee} \in \Delta_{\mu,+}} a_{\al^\vee}e_{\al^\vee}\Big)\Big\}
\end{equation*}
and the coordinate $x_{\al^\vee}$ is given by 
\begin{equation*}
z^\mu\on{exp}\Big(\sum_{\beta^{\vee} \in \Delta_{\mu,+}^\vee} a_{\beta^\vee}e_{\beta^\vee}\Big) \mapsto -a_{\al^\vee},   
\end{equation*}
here $e_{\beta^\vee} \in \mathfrak{g}_{\beta^\vee}$ are Chevalley generators, the sign appears in the formula since the Cartan involution of $\mathfrak{g}$ sends $e_{\al^\vee}$ to $-e_{-\al^{\vee}}$. 
}}
\erem


\ssec{}{Poisson structure}
   In the paper~\cite{NW}, H.~Nakajima and A.~Weekes defined a Poisson structure on $\ol{\CW}^\la_\mu$ for any reductive group $G$ (simply laced case was treated in~\cite{BFN}, we will often refer to some results of~\cite{BFN} in the general case since their proofs can be easily generalized). Let us describe this Poisson structure.
    Let $\la,\,\mu \in \La$ be coweights of $T$ such that $\la$ is dominant and $\mu \leqslant \la$.
    We set $\la^*:=-w_0(\la),\, \mu^*:=-w_0(\mu)$.
    It is known (see~\cite[Theorem~4.1]{NW},~\cite[Theorem~3.10]{BFN}) 
    that the variety $\ol{\CW}^{\la^*}_{\mu^*}$ 
    is isomorphic to a Coulomb branch $\mathcal M(\la,\mu)$ (see Section~\ref{coulomb}) of the corresponding quiver gauge theory.
    It is known (see~\cite[Section~3(iv)]{BFN}) that the Coulomb branch  $\CM(\la,\mu)$ admits a Poisson structure. Let us briefly describe it. 
    
    Recall that the Coulomb branch $\CM(\la,\mu)$ is the spectrum of the algebra of equivariant Borel-Moore homology $H^{\on{GL}(V)[[z]]}_*(\CR_{\on{GL}(V),\,{\bf{N}}^\la_\mu})$, where the product is given by the convolution (see
    Section~\ref{coulomb}).
    The space of triples $\CR_{\on{GL}(V),\,{\bf{N}}^\la_\mu}$ is equipped with the $\BC^\times$-action via the loop rotation. We obtain a graded deformation $H^{\on{GL}(V)[[z]]\rtimes \BC^\times}_*(\CR_{\on{GL}(V),\,{\bf{N}}^\la_\mu})$ of the algebra $H^{\on{GL}(V)[[z]]}_*(\CR_{\on{GL}(V),\,{\bf{N}}^\la_\mu})$, i.e.\ $H^{\on{GL}(V)[[z]]\rtimes \BC^\times}_*(\CR_{\on{GL}(V),\,{\bf{N}}^\la_\mu})$ is a flat (graded) module over $\BC[\hbar]=H_{\BC^\times}^*(\on{pt})$, 
    such that $H^{\on{GL}(V)[[z]]\rtimes \BC^\times}_*(\CR_{\on{GL}(V),\,{\bf{N}}^\la_\mu})/(\hbar-1)=H^{\on{GL}(V)[[z]]}_*(\CR_{\on{GL}(V),\,{\bf{N}}^\la_\mu})$. We now define the Poisson bracket on $H^{\on{GL}(V)[[z]]}_*(\CR^\la_\mu)$ by the formula $\{[f],[g]\}:=[\frac{fg-gf}{\hbar}]$, where $f,\,g,\, \frac{fg-gf}{\hbar} \in H^{\on{GL}(V)[[z]]\rtimes \BC^\times}_*(\CR_{\on{GL}(V),\,{\bf{N}}^\la_\mu})$ and $[f],\,[g],\,[\frac{fg-gf}{\hbar}]$ are the corresponding elements of $H^{\on{GL}(V)[[z]]}_*(\CR_{\on{GL}(V),\,{\bf{N}}^\la_\mu})$. In this Section, we will describe the Poisson structure on $\CM(\la^*,\mu^*) \simeq \ol{\CW}^\la_\mu$ for minuscule $\la \in \La^+$ and $\mu \in W\la$.
    
    \ssec{Poiss}{Coulomb branches}(see~\cite{N},~\cite[Section~3]{BFN},~\cite{BFN2},~\cite{NW})\label{coulomb}
    Let $Q=(Q_0,Q_1)$ be the Dynkin quiver of the group $G$, here $Q_0$ stands for vertices and $Q_1$ for edges. We also fix any orientation of every element of $Q_1$. Let $(c_{ij})_{i,j \in Q_0}$ be the Cartan matrix of $G$. Let $(d_i)_{i \in Q_0} \in \BZ_{>0}$ be such that $d_ic_{ij}=d_jc_{ji}$.
    We set $f_{ij}=1$ if $i=j$ and $f_{ij}=|c_{ij}|$ otherwise.
    Consider the formal disk $D_i:=\on{Spec}\BC[[z_i]]$ and the punctured disk $\overset{\circ}{D}_i:=\on{Spec}\BC((z_i))$ for each vertex $i \in Q_0$. For each pair $(i,j)$ with $c_{ij}<0$ we take a formal disk $D=\on{Spec}\BC[[z]]$ and consider its branched coverings $\pi_{ji}\colon D_i \ra D,\, 
    \pi_{ij}\colon D_j \ra D$ given by $z_{i} \mapsto z^{f_{ij}},\, z_j \mapsto z^{f_{ji}}$ respectively.
    
    We write $\la=\sum_{i \in Q_0}l_i\omega_i,\, \la-\mu=\sum_{i \in Q_0}a_i\al_i$, where $\omega_i \in \La^+$ are fundamental coweights of $G$ and $\al_i \in \La^{\on{pos}}$ are simple coroots. For each $i \in Q_0$ we set $W_i:=\BC^{l_i},\, V_i:=\BC^{a_i}$,
    ${\bf{N}}^\la_\mu:=\bigoplus_{h \in Q_1}\on{Hom}(V_{o(h)},V_{i(h)})\oplus \bigoplus_{j \in Q_0}\on{Hom}(W_i,V_i)$
    and $\on{GL}(V):=\prod_{i\in Q_0}\on{GL}(V_i)$. The group $\on{GL}(V)$ acts naturally on ${\bf{N}}^\la_\mu$. 
    
    Following~\cite[Sections~2(ii),~5(v)]{NW} we consider 
    the moduli space $\CR=\CR_{\on{GL}(V),\,{\bf{N}}^\la_\mu}$ parametrizing 
    the following objects:
    \begin{itemize}
        \item a rank $a_i$ vector bundle $\CE_i$ on $D_i$ together with a trivialization 
        \begin{equation*}
        \sigma_i\colon \CE_i|_{\overset{\circ}{D}_i} \iso V_i \otimes \CO_{D_i}
        \end{equation*}
        for $i \in Q_0$,
        \item a homomorphism $s_i\colon W_i \otimes \CO_{D_i} \ra \CE_i$ such that $\sigma_i \circ \Big(s_i|_{\overset{\circ}{D}_i}\Big)$ extends to $D_i$ for $i \in Q_0$,
        \item a homomorphism $s_{ij} \in  \on{Hom}_{\CO_D}(\pi_{ij*}\CE_j,\pi_{ji*}\CE_i)$ such that $\Big(\pi_{ij*}\sigma_i\Big) \circ \Big(s_{ij}|_{\overset{\circ}{D}_i}\Big) \circ \Big(\pi_{ji*}\sigma_j\Big)^{-1}$ extends to $D$, where $c_{ij}<0$ and $j \ra i$ is an arrow in the quiver $Q$. 
    \end{itemize}
    Consider now the semidirect product $\on{GL}(V)[[z]] \rtimes \BC^\times$, here $\BC^\times$ acts on $\on{GL}(V)[[z]]$ via the loop rotation. 
    The group  $\on{GL}(V)[[z]] \rtimes \BC^\times$ acts on $\CR_{\on{GL}(V),\,{\bf{N}}^\la_\mu}$. We can consider the algebra of equivariant Borel-Moore homology $H^{\on{GL}(V)[[z]]}_*(\CR_{\on{GL}(V),\,{\bf{N}}^\la_\mu})$.  It follows from~\cite[Proposition~5.15]{BFN2} that this is a commutative algebra, so we can define $\CM(\la,\mu):=\on{Spec}(H^{\on{GL}(V)[[z]]}_*(\CR_{\on{GL}(V),\,{\bf{N}}^\la_\mu}))$. 
    
    \rem{}
    {\em{One can show that $\mathcal{M}(\la,\mu)$ does not depend on the choice of the orientation of the elements of $Q_1$. This was communicated to us by Alex Weekes, the proof is similar to the proof of the similar statement for simply-laced $G$ given in~\cite{BFN2}: first one considers the 
    case when $G=T$ (the proof is similar to~\cite[Section~4(v)]{BFN2}) and then deduces the general case
    using the same approach as in~\cite[Section~6(viii)]{BFN2}).    
    }}
    \erem

    The following
    proposition follows from~\cite[Theorem~3.1]{BFN},~\cite[Theorem~4.1]{NW}. 
    \prop{Comp}\label{iso_slices_coulomb_claim}
    There exists an isomorphism of algebras $\Xi\colon \BC\Big[\ol{\CW}^{\la^*}_{\mu^*}\Big] \iso H^{\on{GL}(V)[[z]]}_*(\CR_{\on{GL}(V),\,{\bf{N}}^{\la}_{\mu}})$. In particular the Poisson bracket $\{\,,\,\}$ on $H^{\on{GL}(V)[[z]]}_*(\CR_{\on{GL}(V),\,{\bf{N}}^{\la}_{\mu}})$ defines a Poisson bracket on $\BC\Big[\ol{\CW}^{\la^*}_{\mu^*}\Big]$ to be denoted by the same symbol $\{\,,\,\}$.
    \eprop

    \prop{}\label{leaf_slice}
    $\CW^{\la^*}_{\mu^*} \subset \ol{\CW}^{\la^*}_{\mu^*}$ is a symplectic leaf of the Poisson bracket $\{\,,\,\}$. We denote the corresponding symplectic from on $\CW^{\la^*}_{\mu^*}$ by $\omega$.
    \eprop
    \prf
    The claim follows from~\cite[Remark~3.19]{BFN} (see also~\cite[Proposition~6.15]{BFN2}) combined with~\cite[Theorem~1.2]{smooth_slice}. See the introduction of~\cite{smooth_slice} for more details. 
    \epr

\sssec{Poiss}{Torus action}
Recall the Cartan torus action $T \curvearrowright \ol{\CW}^{\la^*}_{\mu^*}$ of Section~\ref{actionss}. With respect to the isomorphism $\Xi$ (see Proposition~\ref{iso_slices_coulomb_claim}), we obtain the action $T \curvearrowright H^{\on{GL}(V)[[z]]}_*(\CR_{\on{GL}(V),\,{\bf{N}}^\la_\mu})$. 
Since $G$ is adjoint an action of $T$ on $H^{\on{GL}(V)[[z]]}_*(\CR_{\on{GL}(V),\,{\bf{N}}^\la_\mu})$ is the same as a $\BZ\langle\al^{\vee}_i\rangle_{i \in Q_0}$-grading on $H^{\on{GL}(V)[[z]]}_*(\CR_{\on{GL}(V),\,{\bf{N}}^\la_\mu})$.
The algebra $H^{\on{GL}(V)[[z]]}_*(\CR_{\on{GL}(V),\,{\bf{N}}^\la_\mu})$ is naturally graded by $\pi_1(\on{GL}(V))=\BZ^{Q_0}=\BZ\langle\al^{\vee}_i\rangle_{i \in Q_0}$. This is exactly the grading which corresponds to the $T$-action above (see~\cite[Remark~3.12]{BFN}).
Recall now that the Poisson structure on $H^{\on{GL}(V)[[z]]}_*(\CR^\la_\mu)$ comes from the deformation $H^{\on{GL}(V)[[z]] \rtimes \BC^\times}_*(\CR_{\on{GL}(V),\,{\bf{N}}^\la_\mu})$. 
It is easy to see that the $\BZ^{Q_0}$-grading on $H^{\on{GL}(V)[[z]]}_*(\CR^\la_\mu)$ extends to the $\BZ^{Q_0}$-grading on the associative $\BC[\hbar]$-algebra
$H^{\on{GL}(V)[[z]] \rtimes\BC^\times}_*(\CR_{\on{GL}(V),\,{\bf{N}}^\la_\mu})$, here $\BC[\hbar]=H_{\BC^\times}^*(\on{pt})$. 
As a corollary we obtain the following lemma. 

\lem{Poiss}\label{cartan_poiss}
The action $T \curvearrowright \BC\Big[\ol{\CW}^{\la^*}_{\mu^*}\Big] \simeq H^{\on{GL}(V)[[z]]}_*(\CR_{\on{GL}(V),\,{\bf{N}}^\la_\mu})$ is Poisson.
\elem

\sssec{Poiss}{Cartan involution}\label{cart_inv_prop} Recall the Cartan involution $\iota\colon \ol{\CW}^{\la^*}_{\mu^*} \iso \ol{\CW}^{\la^*}_{\mu^*}$ (see~\cite[Section~2(vii)]{BFN}). Let us describe the corresponding automorphismism $\iota^*\colon H^{\on{GL}(V)[[z]]}_*(\CR_{\on{GL}(V),\,{\bf{N}}^\la_\mu}) \iso H^{\on{GL}(V)[[z]]}_*(\CR_{\on{GL}(V),\,{\bf{N}}^\la_\mu})$ (see~\cite[Remarks~3.6,~3.16]{BFN}). 

Let $\mathfrak{i}\colon \on{Gr}_{\on{GL}(V)} \simeq 
\on{Gr}_{\on{GL}(V^*)}$
be the following automorphismism: it takes $(\CP,\sigma)$ to $(\CP^{\vee},{}^t\sigma^{-1})$. Let $\ol{Q_1}$ be the opposite orientation of our quiver. Consider the representation ${\bf{N}}^{\la}_\mu:=\bigoplus_{h \in \ol{Q_1}}\on{Hom}(V_{o(h)}^*,V_{i(h)}^*)\oplus \bigoplus_{i \in Q_0} \on{Hom}(V^*_i,W^*_i)$ 
of $\on{GL}(V^*)$.
It is easy to see that $\mathfrak{i}$ lifts to the
isomorphism 
$\mathfrak{i}^\la_\mu\colon \CR_{\on{GL}(V),\,{\bf{N}}^\la_\mu} \iso \CR_{\on{GL}(V^*),\,{\bf{N}}^{\la}_\mu}$, which together with the automorphismism $\on{GL}(V) \iso \on{GL}(V^*),\, g \mapsto {}^tg^{-1}$ induces the convolution algebra isomorphism 
\begin{equation*}
\mathfrak{i}^\la_{\mu*}\colon H^{\on{GL}(V)[[z]]}_*(\CR_{\on{GL}(V),\,{\bf{N}}^\la_\mu}) \iso H^{\on{GL}(V)[[z]]}_*(\CR_{\on{GL}(V^*),\,{\bf{N}}^{\la}_\mu}).
\end{equation*}
The composition 
$\BC\Big[\ol{\CW}^{\la^*}_{\mu^*}\Big] \simeq H^{\on{GL}(V)[[z]]}_*(\CR_{\on{GL}(V),\,{\bf{N}}^\la_\mu}) \iso H^{\on{GL}(V)[[z]]}_*(\CR_{\on{GL}(V^*),\,{\bf{N}}^{\la}_\mu}) \simeq 
\BC\Big[\ol{\CW}^{\la^*}_{\mu^*}\Big]$ is an involution of the algebra $\BC\Big[\ol{\CW}^{\la^*}_{\mu^*}\Big]$ to be denoted $\tilde{\iota}^*$. This coincides with the Cartan involution $\iota^*$ composed with the involution $\varkappa_{-1}$ of $\ol{\CW}^{\la^*}_{\mu^*}$ induced by an automorphismism $\BP^1 \iso \BP^1,\, z \mapsto -z$ and finally composed with an action of a certain element of the Cartan torus $T$.  
We are now ready to prove the following lemma.

\lem{poiss}\label{invol_poiss}
Recall the Cartan involution $\iota\colon \ol{\CW}^{\la^*}_{\mu^*} \iso \ol{\CW}^{\la^*}_{\mu^*}$.
Then $\iota^*\colon \BC\Big[\ol{\CW}^{\la^*}_{\mu^*}\Big] \iso \BC\Big[\ol{\CW}^{\la^*}_{\mu^*}\Big]$ is an antiautomorphism of the Poisson algebra $\BC\Big[\ol{\CW}^{\la^*}_{\mu^*}\Big]\simeq H^{\on{GL}(V)[[z]]}_*(\CR_{\on{GL}(V),\,{\bf{N}}^\la_\mu})$ (i.e. $\{\iota^*(f_1),\iota^*(f_2)\}=-\iota^*\{f_1,f_2\}$ for any $f_1,f_2 \in \BC\Big[\ol{\CW}^{\la^*}_{\mu^*}\Big]$).
\elem

\prf
It is easy to see that the isomorphism $\mathfrak{i}^\la_\mu$ is $\BC^\times$-equivariant (with respect to the loop rotation action). It follows that the isomorphism $\mathfrak{i}^\la_\mu$ induces the (graded) convolution algebra isomorphism 
\begin{equation*}
\mathfrak{i}^\la_{\mu*}\colon H^{\on{GL}(V)[[z]]\rtimes \BC^\times}_*(\CR_{\on{GL}(V),\,{\bf{N}}^\la_\mu}) \iso H^{\on{GL}(V)[[z]]\rtimes \BC^\times}_*(\CR_{\on{GL}(V^*),\,{\bf{N}}^{\la}_\mu}),
\end{equation*}
hence, $\tilde{\iota}^*$ is an automorphism of Poisson algebras. Note also that the automorphism $\BP^1 \iso \BP^1,\, z \mapsto -z$ induces the antiautomorphism of the graded algebra
$
H^{\on{GL}(V)[[z]] \rtimes \BC^{\times}}_*(\CR_{\on{GL}(V),\,{\bf{N}}^\la_\mu}).
$
Let us finally note that the action of the Cartan torus $T$ on $\ol{\CW}^{\la^*}_{\mu^*}$ is Poisson by Lemma~\ref{cartan_poiss}. 
\epr

\sssec{poiss}{Loop rotation}
Recall the loop rotation action $\BC^\times \curvearrowright \ol{\CW}^{\la^*}_{\mu^*}$. Consider also the morphism $s^{\la^*}_{\mu^*} \colon \ol{\CW}^{\la^*}_{\mu^*} \ra Z^{\al}$, where $\al:=\la-\mu$. It follows from~\cite{BFN},~\cite[Theorem~4.1]{NW} that this morphism is Poisson. It also follows from~\cite[Proposition~2.10]{BFN} that the restriction of $s^{\la^*}_{\mu^*}$ to $(s^{\la^*}_{\mu^*})^{-1}(\oZal)$ is an isomorphism 
$(s^{\la^*}_{\mu^*})^{-1}(\oZal) \iso \oZal$. We obtain the embedding of Poisson algebras $\BC\Big[\ol{\CW}^{\la^*}_{\mu^*}\Big] \hookrightarrow \BC\Big[\oZal\Big]$ (c.f.~\cite[Remark~3.11]{BFN}~\cite[Remark~5.14]{BFN2}). Note also that this embedding is $\BC^{\times}$-equivariant (with respect to the loop rotation action).

Recall the $\BC^\times$-action on $\ol{\CW}^{\la^*}_{\mu^*}$ via loop rotation. It induces a $\BZ$-grading on the algebra of functions $\BC\Big[\ol{\CW}^{\la^*}_{\mu^*}\Big]=\bigoplus_{n \in \BZ} \BC\Big[\ol{\CW}^{\la^*}_{\mu^*}\Big]_n$: 
\begin{equation}\label{our_grading_slice}
\BC\Big[\ol{\CW}^{\la^*}_{\mu^*}\Big]_n=\big\{f \in \BC\Big[\ol{\CW}^{\la^*}_{\mu^*}\Big],\,|\, f(t^{-1} \cdot x)=t^nf(x),\, t \in \BC^\times,\, x \in \ol{\CW}^{\la^*}_{\mu^*}\big\}.
\end{equation}

\rem{} 
{\em{Note that the loop rotation action $\BC^\times  \curvearrowright \ol{\CW}^{\la^*}_{\mu^*}$ contracts $\ol{\CW}^{\la^*}_{\mu^*}$ to its fixed points iff the grading~\ref{our_grading_slice} is nonpositive i.e. $\BC\Big[\ol{\CW}^{\la^*}_{\mu^*}\Big]_n=0$ for $n>0$.
}}
\erem

\lem{poiss}\label{deg_loop}
 The Poisson structure $\{\,,\,\}$ has degree $1$ with respect to the grading~(\ref{our_grading_slice})  (i.e. for functions $f_i,\,f_j \in \BC\Big[\ol{\CW}^{\la^*}_{\mu^*}\Big]$ of degrees $i,\,j$ respectively, the function $\{f_i,f_j\}$ has degree $i+j+1$). 
\elem

\prf
We have a $\BC^\times$-equivariant Poisson embedding $\BC\Big[\ol{\CW}^{\la^*}_{\mu^*}\Big] \hookrightarrow \BC\Big[\oZal\Big]$, so it is enough to prove our statement for $\oZal$ and for zastava space this statement follows from the computation of $\{\,,\,\}$ in \'etale coordinates made in~\cite[Section~1.2]{FKMM}. Recall the \'etale coordinates $(w_{i,r},\, y_{i,r})$, where $i \in Q_0$ and $r=1,\ldots,b_i$ ($\al=\sum_i b_i\al_i$). We have $\{w_{i,r},w_{j,s}\}=0,\,
\{w_{i,r},y_{j,s}\}=\delta_{ij}\delta_{rs}y_{j,s}$ and
$\{y_{i,r},y_{j,s}\}=(\al^{\vee}_i,\al^{\vee}_j)\frac{y_{i,r}y_{j,s}}{w_{i,r}-w_{j,s}}$ for $i \neq j$,
$\{y_{i,r},y_{i,s}\}=0$. 
Here $\al^{\vee}_i$ is a simple root of $G$, $(\,,\,)$ is the invariant scalar product on $(\on{Lie}T)^*$, such that the square length of a short root is $2$.
Note that $\on{deg}(w_{i,r})=-1,\, \on{deg}(y_{j,s})=0$ (the degree with respect to the loop rotation action).
The lemma is proven.
\epr

\rem{poiss}
{{Lemma~\ref{deg_loop} can be also deduced from~\cite[Remark~3.13]{BFN}.}}
\erem

Let $\hbar$ be the character $1$ of the loop rotating $\BC^\times$. We denote by the same symbol the character of $T \times \BC^\times$ induced by $\hbar$.
\cor{}\label{poiss_eigen}
The Poisson bivector $\Theta \in \on{Hom}_{\CO_{\ol{\CW}^{\la^*}_{\mu^*}}}(\La^2\Omega^1_{\ol{\CW}^{\la^*}_{\mu^*}},\CO_{\ol{\CW}^{\la^*}_{\mu^*}})$ corresponding to the Poisson bracket $\{\,,\,\}$ is an eigenvector of $T \times \BC^\times$-action with eigenvalue $\hbar$. In particular the symplectic form $\omega \in \Gamma(\CW^{\la^*}_{\mu^*},\La^2\Omega^1_{\CW^{\la^*}_{\mu^*}})$ is an eigenvector of $T \times \BC^\times$-action with eigenvalue $-\hbar$.
\ecor
\prf
Follows from Lemmas~\ref{cartan_poiss},~\ref{deg_loop}.
\epr

\ssec{poiss}{Minuscule case}
Now we assume that $\la=\mu^+$ is miniscule. Since for minuscule slices we have $\ol{\CW}^{\mu^+}_\mu=\CW^{\mu^+}_\mu$ it follows from Proposition~\ref{leaf_slice} and  \cref{cartan_poiss} that $\CW^{\mu^+}_\mu$ is a $T$-equivariant symplectic variety. 
    Recall the $T$-fixed point $z^\mu \in \CW^{\mu^+}_\mu$ and note that we have the $T$-equivariant Lagrangian decomposition $
    T_{z^{\mu}}\ol{\CW}^\la_\mu =
    T_{z^\mu}\CA^\la_\mu \oplus T_{z^\mu}\CR^\la_\mu 
    $
    due to the $T$-invariance of the symplectic form $\omega$.
    
   \sssec{poiss}{$T \times \BC^\times$-character of minuscule slice}
    \prop{poiss}\label{character_minis}
    Recall that $\mu^+ \in \La^+$ is minuscule.
    The character of $T \times \BC^\times$ acting on $T_{z^\mu}\CW^{\mu^+}_\mu$ equals to  $\sum_{\al^{\vee} \in \Delta^{\vee}_{\mu,-}}(e^{\al^{\vee}}+ e^{-\al^{\vee}+\hbar})$.
    \eprop 
    \prf
    Recall that we have an 
    identification $T_{z^\mu}\mathcal {W}_\mu^\lambda=T_{z^{\mu}}\CR^\la_\mu \oplus T_{z^{\mu}}\CA^\la_\mu$. Let us compute the $T\times \BC^\times$-character of  $T_{z^{\mu}}\CR^\la_\mu$. Recall that the loop rotation acts trivially on $\CR^\la_\mu$. Recall also that we have a $T$-equivariant embedding $\CR^\la_\mu \hookrightarrow G/P_{\la}$ which identifies $\CR^\la_\mu$ with the repellent $X^\la_\mu$ to the point $wP_{\la}/P_{\la} \in G/P_{\la}$ with respect to the $\BC^\times$-action via $2\rho$. It remains to compute the $T$-character of the tangent space $T_{wP_{\la}/P_{\la}}X^\la_\mu$. It is easy to see that it is equal to $\sum_{\al^{\vee} \in \Delta^{\vee}_{\mu,-}}e^{\al^{\vee}}$. We can now compute the character of $T_{z^\mu}\CA^\la_\mu$ using Corollary~\ref{poiss_eigen}. Indeed the symplectic form $\omega$ is clearly nondegenerate and has $T\times \BC^\times$-weight $-\hbar$ so it induces an isomorphism of $T\times \BC^\times$ representations  $T_{z^\mu}\CA^\la_\mu \simeq e^{\hbar} (T_{z^\mu}\CR^\la_\mu)^*$. We conclude that the $T \times \BC^\times$-character of $T_{z^\mu}\CA^\la_\mu$ is equal to $\sum_{\al^{\vee} \in \Delta^{\vee}_{\mu,-}} e^{-\al^{\vee}+ \hbar}$.
    \epr
    
\rem{}{\em{Let us give another proof of the Proposition~\ref{character_minis} that uses  only the matrix description of $\CW^{\mu^+}_\mu.$
We are gratefull to the anonymous referee for explaining this approach to us. 
    Recall that 
    \begin{equation*}
    \CW^{\mu^+}_{\mu}= \big(U[[z^{-1}]]_{1}z^{\mu}B_{-}[[z^{-1}]]_{1}\cap G[z]z^{\mu^+}G[z]\big) \subset G((z^{-1})).
    \end{equation*}
    Note that we can consider $U[[z^{-1}]]_{1}z^{\mu}B_{-}[[z^{-1}]]_{1},\, G[z]z^{\mu^+}G[z]$ as subfunctors of the functor $G((z^{-1}))$ (see~\cite{smooth_slice} for details) and by~\cite[Theorem~1.2]{smooth_slice} their intersection $\CW^{\mu^+}_{\mu}$ is smooth (considered as a subscheme of $G((z^{-1}))$). By passing to $\BC[\epsilon]/\epsilon^2$-points of our functors we conclude that
    \begin{equation*}
    T_{z^\mu}\CW^{\mu^+}_\mu=T_{z^\mu}\big(U[[z^{-1}]]_{1}z^{\mu}B_{-}[[z^{-1}]]_{1}\big) \cap T_{z^{\mu}}\big(G[z]z^{\mu^+}G[z]\big) \subset T_{z^\mu}G((z^{-1}))=\mathfrak{g}((z^{-1})).
    \end{equation*}
    Since $G[z]z^{\mu^+} G[z]=G[z]z^\mu G[z]$ we conclude  that $T_{z^\mu}\big(G[z]z^{\mu^+} G[z]\big) \subset \mathfrak{g}(z)$ is generated by vectors of the form \begin{equation*}
    v(z),\, z^{\mu}v(z) z^{-\mu},\,\text{where}~v(z) \in \mathfrak{g}[z].
    \end{equation*}
    Note also that every element of $T_{z^\mu}\big(U[[z^{-1}]]_{1}z^{\mu}B_{-}[[z^{-1}]]_{1}\big) \subset \mathfrak{g}((z^{-1}))$ can be obtained as (possibly infinite) linear combination of vectors of the form \begin{equation*}
    z^{-n-1}v_{\beta^\vee},\, z^{-n-1+\langle\al^\vee,\mu\rangle}v_{\al^\vee},\,\text{where}~v_{\beta^\vee} \in \mathfrak{g}_{\beta^{\vee}},\, v_{\al^\vee} \in \mathfrak{g}_{\al^{\vee}},\, n \geqslant 0,\, \beta^\vee \in \Delta_+,\, \al^\vee \in \Delta_-^\vee.
    \end{equation*}
    Recall now that $\mu^+$ is minuscule and $\mu \in W\mu^+$ so we have $\langle\gamma^\vee,\mu \rangle \in \{-1,0,1\}$ for every $\gamma^\vee \in \Delta^\vee$. It follows then that the tangent space to $\CW^{\mu^+}_\mu$ is generated by vectors 
    \begin{equation*}
    z^{-1}v_{\beta^{\vee}},\,  v_{\al^\vee},\,\text{such that}~\beta^\vee \in \Delta_+^\vee,\, \al^\vee \in \Delta_-^\vee,\,
    \langle\beta^\vee,\mu \rangle=-1,\,
    \langle \al^\vee,\mu \rangle=1.    
    \end{equation*}
    Note that the conditions $\beta^\vee \in \Delta_+^\vee,\, \langle\beta^\vee,\mu \rangle=-1$ are equivalent to the condition $-\beta^\vee \in \Delta^\vee_{\mu,-}$ and the conditions 
    $\al^\vee \in \Delta_-^\vee,\, \langle\al^\vee,\mu \rangle=1$ are equivalent to the condition $\al^\vee \in \Delta^\vee_{\mu,-}$. It remains to note that $v_{\al^\vee}$ has $T \times \BC^\times$-weight $e^{\al^\vee}$ and $z^{-1}v_{\beta^\vee}$ has $T \times \BC^\times$-weight $ e^{\beta^{\vee} + \hbar}$. 
    
    Another way to prove Proposition~\ref{character_minis} is to use matrix descriptions of $\CR^{\mu^+}_\mu,\, \CA^{\mu^+}_\mu$ given in Remark~\ref{matrix_descr_A_R_min}.
    }}
  \erem  

    We describe the form $\omega$ in the following theorem.
    We are gratefull to the anonymous referee for strengthening our original theorem and suggesting the possible proof of this strengthened version.
    Recall that $\{y_{\beta^{\vee}}\,|\, \beta^{\vee} \in \Delta^{\vee}_{\mu,-}\}$, $\{x_{\al^{\vee}}\,|\, \al^{\vee} \in \Delta^{\vee}_{\mu,+}\}$ are the standard coordinates on the affine spaces $ \CR^{\mu^+}_\mu$ and $\CA^{\mu^+}_\mu$. Recall the isomorphism $\CW^{\mu^+}_\mu \simeq \CA^{\mu^+}_\mu \times \CR^{\mu^+}_\mu$. 
\th{poiss}\label{poiss_structure} 
The symplectic form $\omega$ is $\sum_{\al^{\vee} \in \Delta^{\vee}_{\mu,+}} c_{\al^\vee} dx_{\al^{\vee}} \wedge dy_{-\al^{\vee}}$ for some $c_{\al^\vee} \in \BC^\times$. In other words after the rescalling $y'_{-\al^\vee}=\frac{1}{c_{\al^\vee}}y_{-\al^\vee}$ the form $\omega$ becomes standard $\omega=\sum_{\al^\vee \in \Delta^\vee_{\mu,+}}dx_{\al^\vee} \wedge dy'_{-\al^\vee}$.
\eth
    \prf
    Recall the grading~(\ref{our_grading_slice}) on $\BC\Big[\CW^{\mu^+}_\mu\Big]$ induced by the loop rotation action.
    Let us first of all compute the degrees of the coordinates $y_{-\al^\vee},\, x_{\al^\vee}$ with respect to this grading. Note that the loop rotation action $\BC^\times \curvearrowright \CR^{\mu^+}_\mu,\,\CA^{\mu^+}_\mu$ induces $\BZ$-gradings on $\BC\Big[\CR^{\mu^+}_\mu\Big],\,\BC\Big[\CA^{\mu^+}_\mu\Big]$.
    It follows from Remark~\ref{matrix_descr_A_R_min} that $\CA^\la_\mu$ has the following matrix description 
    \begin{equation*}
    \CA^{\mu^+}_\mu=   \{z^{\mu}u\,|\, u \in U^\mu\}=\Big\{z^\mu\on{exp}\Big(\sum_{\al^{\vee} \in \Delta_{\mu,+}} a_{\al^\vee}e_{\al^\vee}\Big)\Big\}= 
    \Big\{\on{exp}\Big(\sum_{\al^{\vee} \in \Delta_{\mu,+}} a_{\al^\vee} z^{-1} e_{\al^\vee}\Big)z^\mu\Big\}.
    \end{equation*}
    It then follows from Proposition~\ref{action_via_matrix} $(a)$ together with Remark~\ref{matrix_descr_A_R_min} that the restriction of $x_{\al^\vee}$ to $\CA^\la_\mu$ has degree $-1$.   Recall now that by Corollary~\ref{iso_U_-_A} we have an isomorphism $U^\mu_- \times \CA^{\mu^+}_\mu \iso \CW^{\mu^+}_\mu$ given by the action morphism $(u_-,x) \mapsto u_- \cdot x$ and by the definition function $x_{\al^\vee} \in \BC\Big[\CW^{\mu^+}_\mu\Big] \simeq \BC\Big[U^\mu_- \times \CA^{\mu^+}_\mu\Big]$ is the pull back (with respect to the projection morphism $U^\mu_- \times \CA^{\mu^+}_\mu \twoheadrightarrow \CA^{\mu^+}_\mu$) of its restriction to $\CA^{\mu^+}_\mu$. Note now that the loop rotation action on $\CW^{\mu^+}_\mu$ commutes with the $U_-$-action (follows from the descriptions of these actions in geometric terms or from Proposition~\ref{action_via_matrix} $(a)$ together with Remark~\ref{act_U-_matrix}) so we conclude that the degree of $x_{\al^\vee}$ coincides with the degree of its restriction to $\CA^{\mu^+}_\mu$ and so is equal to $-1$. Let us now compute the degree of $y_{\beta^\vee}$. Recall that we identify $\CW^{\mu^+}_\mu \simeq U^\mu_- \times \CA^{\mu^+}_\mu \simeq \CR^{\mu^+}_\mu \times \CA^{\mu^+}_\mu$ and loop rotation acts via its action on $\CA^{\mu^+}_\mu$. It follows that coordinates $y_{\beta^\vee}$ have degree $0$ being pull backs of their restrictions to $\CR^{\mu^+}_\mu$. So we have shown that 
$
    \on{deg}x_{\al^\vee}=-1,\, \on{deg}y_{\beta^\vee}=0.
$

    Recall now that $\{\,,\,\}$ has degree $1$ so it follows that $\{x_{\al^\vee},y_{\beta^{\vee}}\}$ has degree $0$ i.e. is $\BC^\times$-invariant. It follows that $\{x_{\al^\vee},y_{\beta^{\vee}}\} \in \BC\Big[\CW^{\mu^+}_\mu\Big]$ should be a polynomial on $y_{\gamma^\vee},\, \gamma^\vee \in \Delta^\vee_{\mu,-}$ (here we use that $\BC\Big[\CW^{\mu^+}_\mu\Big]=\BC[x_{-\gamma^{\vee}},y_{\gamma^{\vee}}]$ and $y_{\gamma^{\vee}}$ are $\BC^\times$-invariant while $x_{-\gamma^{\vee}}$
    has degree $-1$). So we conclude that $\{x_{\al^\vee},y_{\beta^{\vee}}\}=r(y_{\gamma^\vee})$, where $r(y_{\gamma^\vee})$ is some polynomial on all possible $y_{\gamma^\vee}$. Recall now that by Lemma~\ref{invol_poiss} the Cartan involution $\iota$ induces an antiautomorphism $\iota^*$ of the Poisson algebra $\BC\Big[\CW^{\mu^+}_\mu\Big]$. Recall also that directly from the definitions $\iota^*$ sends $y_{\gamma^\vee}$ to $x_{-\gamma^{\vee}}$. We conclude that 
    $\{x_{-\beta^{\vee}},y_{-\al^\vee}\}=r(x_{-\gamma^\vee})$. On the other hand  by the same observations as above $\{x_{-\beta^{\vee}},y_{-\al^\vee}\}$ must be a polynomial on $y_{\gamma^\vee}$. We conclude that $\{x_{-\beta^{\vee}},y_{-\al^\vee}\} \in \BC$. So we have shown that for every $\al^\vee \in \Delta^\vee_{\mu,+},\, \beta^\vee \in \Delta^\vee_{\mu,-}$ 
    we have
    $
    \{x_{\al^{\vee}},y_{\beta^\vee}\}=c_{\al^\vee,\beta^\vee} \in \BC    
    $
    for some constants $c_{\al^\vee,\beta^\vee}$.

    Let us now show that $c_{\al^\vee,\beta^\vee}=0$ if $\beta^\vee \neq -\al^\vee$. Note  that $x_{\al^{\vee}}$ has $T$-weight $-\al^\vee$ and $y_{\beta^\vee}$ has $T$-weight $-\beta^\vee$ (this follows from the fact that the isomorphism $U^-_\mu \times \CA^{\mu^+}_\mu \iso \CW^{\mu^+}_\mu,\, (u_-,x) \mapsto u_- \cdot x$ is $T$-equivariant with respect to the $T$-action on $U^-_\mu \times \CA^{\mu^+}_\mu$ via $t \cdot (u_-,x)=(tu_-t^{-1},t x)$ and our identification $U^-_\mu \simeq \CR^{\mu^+}_\mu$ is $T$-equivariant i.e. the identification $\CW^{\mu^+}_\mu \simeq \CR^{\mu^+}_\mu \times \CA^{\mu^+}_\mu$ is $T$-equivariant). Since $\{\,,\,\}$ is $T$-equivariant it follows that the bracket $ \{x_{\al^{\vee}},y_{\beta^\vee}\}$ must have $T$-weight $-\al^\vee-\beta^\vee$. On the other hand we have shown that this bracket is equal to some constant $c_{\al^\vee,\beta^\vee}$ so we conclude that $\{x_{\al^{\vee}},y_{\beta^\vee}\}=c_{\al^\vee,\beta^\vee}=0$ if $\beta^\vee \neq -\al^\vee$.
    So the possibly nonzero numbers among $c_{\al^\vee,\beta^\vee}$ are $c_{\al^\vee,-\al^\vee}$ that we will simply denote by $c_{\al^\vee}$.
    
    Let us now show that \begin{equation*}
    \{x_{\al^\vee_1},x_{\al_2^\vee}\}=0= \{y_{-\al^\vee_1},y_{-\al_2^\vee}\},\, \forall \al_1^\vee,\, \al_2^\vee \in \Delta^\vee_{\mu,+}.
    \end{equation*}
    Indeed recall that $y_{-\al_1^\vee},\, y_{-\al_2^\vee}$ are $\BC^\times$-invariant and the Poisson bracket $\{\,,\,\}$ has degree $1$. It follows that $\{y_{-\al_1^\vee}, y_{-\al_2^\vee}\}$ has degree $1$. Note now that the $\BC^\times$-grading on $\BC\Big[\CW^{\mu^+}_\mu\Big]$ is nonpositive so we must have $\{y_{-\al_1^\vee}, y_{-\al_2^\vee}\}=0$. Applying Cartan involution we conclude that $\{x_{\al_1^\vee},x_{\al_2^\vee}\}=0$.

    It follows from the above that $\omega$ is given by $\sum_{\al^\vee \in \Delta^\vee_{\mu,+}}c_{\al^\vee} dx_{\al^\vee} \wedge dy_{-\al^{\vee}}$.
    Let us now finally note that the constants $c_{\al^\vee}$ are nonzero since $\omega$ is nondegenerate.
    The Theorem is proven.
    \epr
    
    \rem{}
    {\em{Note that the constants $c_{\al^\vee}$ depend on the choice of $\{e_{\beta^\vee}\}_{\beta^{\vee} \in \Delta}$. Assume that $\mathfrak{g}$ is semisimple and let $(\,,\,)$ be the Killing form on $\mathfrak{g}$. We can then normalize elements $e_{\beta^{\vee}}$ in such a way that $(e_{\beta^{\vee}},e_{-\beta^{\vee}})=1$. It is easy to see that after this normalization constants $c_{\al^\vee}$  from Theorem~\ref{poiss_structure} are determined uniquely.  
    It would be interesting to compute them.}}
    \erem

	\sec{ap}{Covering of convolution diagrams over slices}\label{application}
\ssec{ap}{Convolution diagrams over generalized slices}\label{conv_diagr_sl}
Fix $\la \in \La^+,\, \mu \in \La,\, \mu \leqslant \la$. 
	Fix $N \in \BZ_{\geqslant 1}$ and pick $N$-tuples $\ul{\la}=(\la_1,\ldots,\la_N),\, \la_i \in \La^+$ such that $\la_1+\ldots+\la_N=\la$. For any other $N$-tuple $\ul{\nu}$ we say that $\ul{\nu} \leqslant 
	\ul{\la}$ if $\nu_i \leqslant \la_i$ for $i=1,\ldots,N$. We also set $|\ul{\nu}|=\nu_1+\ldots+\nu_N$.
	\defe[Convolution diagrams over slices]{conv_slice}
	Following~\cite[Section~3.4]{F},~\cite[Section~5(i)]{BFN3} we define $\widetilde{\CW}^{\ul{\la}}_{\mu}$ as the moduli space of the data $(\CP^{triv}=\CP_0,\CP_1,\ldots,\CP_N,\sigma_1,\ldots,\sigma_N,\phi)$, where
	
	(a) $\CP_i$ is a $G$-bundle on $\BP^1$;
	
	(b) $\sigma_i\colon \CP_{i-1}|_{\BP^1 \setminus \{0\}} \iso \CP_i|_{\BP^1 \setminus \{0\}}$ is an isomorphism having a pole of degree $\leqslant \la_i$ at zero;
	
	(c) $\phi$ is a $B$-structure on $\CP_N$ of degree $w_0\mu$, having no defect at $\infty$ and having fiber $B_-$ at $\infty$ with respect to $\sigma_N \circ \sigma_{N-1} \circ \ldots \circ \sigma_1$.
	
We denote by $\CW^{\ul{\la}}_\mu$ the open subscheme of $\widetilde{\CW}^{\ul{\la}}_{\mu}$ consisting of $(\CP_0,\CP_1,\ldots,\CP_N,\sigma_1,\ldots,\sigma_N,\phi)$ such that for $i=1,\ldots,N$ the trivialization $\sigma_i$ has pole of degree equal to $\la_i$ at zero.	
	\edefe

\defe[Convolution affine Grassmannian]{}	Consider the set 
\begin{equation}\label{def_matrix_conv_gr}
\underbrace{G({\CK}) \times_{G({\CO})} \times \ldots \times_{G(\CO)} G({\CK}) \times_{G({\CO})} \on{Gr}_G}_{\text{N}}.
\end{equation}
This set of points  has  a  natural  structure  of  an  ind-scheme  over $\BC$ and  this  object  is  called  a convolution affine Grassmannian and will be denoted by  $\widetilde{\on{Gr}}_{G,N}$. 
It is the moduli space of the following data:

(a) $G$-bundles $\CP^{triv}=\CP_0,\CP_1,\ldots,\CP_N$ on $\BP^1$.
	
(b) Isomorphisms $\sigma_i\colon \CP_{i-1}|_{\BP^1 \setminus \{0\}} \iso \CP_i|_{\BP^1 \setminus \{0\}}$, $i=1,\ldots,N$.
\edefe

Let us denote by $\widetilde{\on{Gr}}_{G,N}^{\ul{\la}} \subset \widetilde{\on{Gr}}_{G,N}$ the closed {\em{reduced}} subscheme of $\widetilde{\on{Gr}}_G$ consisting of the following data:

(a) $G$-bundles $\CP^{triv}=\CP_0,\CP_1,\ldots,\CP_N$ on $\BP^1$.
	
(b) Isomorphisms $\sigma_i\colon \CP_{i-1}|_{\BP^1 \setminus \{0\}} \iso \CP_i|_{\BP^1 \setminus \{0\}}$ having a pole of degree $\leqslant \la_i$ at zero.

We denote by $\on{Gr}^{\ul{\la}}_{G,N}$ the open subscheme of $\widetilde{\on{Gr}}^{\ul{\la}}_{G,N}$ consisting of $(\CP_0,\ldots,\CP_N,\sigma_1,\ldots,\sigma_N)$ such that $\sigma_i$ has pole of degree equal to $\la_i$ at zero.
Directly from the definitions we have 
\begin{equation*}
\widetilde{\CW}^{\ul{\la}}_\mu=\widetilde{\on{Gr}}_{G,N}^{\underline{\la}} \times_{'\on{Bun}_G(\BP^1)} \on{Bun}^{w_0\mu}_B(\BP^1),\, \CW^{\ul{\la}}_\mu=\on{Gr}_{G,N}^{\underline{\la}} \times_{'\on{Bun}_G(\BP^1)} \on{Bun}^{w_0\mu}_B(\BP^1), 
\end{equation*}
where $'\on{Bun}_G(\BP^1)$ is the stack of $G$-bundles on $\BP^1$ with a $B$-structure at $\infty$ and $\on{Bun}^{w_0\mu}_B(\BP^1)$ is the stack of $B$-bundles on $\BP^1$ of degree $w_0\mu$.
	
We have the natural morphism $\varpi^{\ul{\la}}_\mu\colon \widetilde{\CW}^{\ul{\la}}_\mu \ra \ol{\CW}^\la_\mu$ given by 
	\begin{equation}\label{morp_compose}
	(\CP_0,\CP_1,\ldots,\CP_N,\sigma_1,\ldots,\sigma_N,\phi) \mapsto (\CP_N, \sigma_N\circ \sigma_{N-1}\circ\ldots\circ\sigma_1,\phi).
	\end{equation}
	This morphism is a partial resolution of singularities, it is a resolution of singularities when all $\la_i$ are {\em{minuscule}} coweights (see Proposition~\ref{conv_prop} bellow).

The following proposition is well-known to the experts. 	
	\prop{ap}\label{conv_prop} The following holds.
	\begin{enumerate}
	\item 
	The morphism $\varpi^{\ul{\la}}_\mu\colon \widetilde{\CW}^{\ul{\la}}_\mu \ra \ol{\CW}^\la_\mu$ is projective and stratified semismall with respect to the stratifications $\ol{\CW}^\la_\mu=\bigsqcup_{\mu \leqslant \nu \leqslant \la}\CW^{\nu}_\mu,\, \widetilde{\CW}^{\ul{\la}}_\mu=\bigsqcup_{\ul{\nu} \leqslant \ul{\la},\,\mu \leqslant |\ul{\nu}|}\CW^{\ul{\nu}}_\mu$, in particular, $\varpi^{\ul{\la}}_\mu$ is birational. 
	\item 
	The variety $\widetilde{\mathcal{W}}^{\ul{\la}}_\mu$ is Cohen-Macaulay and normal.
    \item 
    Assume that $\la_i$ are minuscule for $i=1,\ldots,N$, then the variety $\widetilde{\mathcal{W}}^{\ul{\la}}_\mu$ is smooth.
	\end{enumerate}
	\eprop
	\prf
	We have 
	$\widetilde{\CW}^{\ul{\la}}_\mu=\widetilde{\on{Gr}}^{\ul{\la}}_{G,N} \times_{'\on{Bun}_G(\BP^1)} \on{Bun}^{w_0\mu}_B(\BP^1)$,
	$\overline{\CW}^{{\la}}_\mu=\ol{\on{Gr}}{}^{\la}_{G} \times_{'\on{Bun}_G(\BP^1)} \on{Bun}^{w_0\mu}_B(\BP^1)$ and the morphism $\varpi^{\ul{\la}}_\mu$ is the base change of the morphism 
	$\widetilde{\on{Gr}}^{\ul\la}_{G,N} \ra \ol{\on{Gr}}{}^{\la}_{G}$ that is projective and stratified semismall with respect to the stratifications 
	$\overline{\on{Gr}}{}^\la_G=\bigsqcup_{\nu \leqslant \la} \on{Gr}_G^\nu$, $\widetilde{\on{Gr}}^{\ul{\la}}_{G,N}=\bigsqcup_{\ul{\nu} \leqslant \ul{\la}, \mu \leqslant |\ul{\nu}|}\on{Gr}^{\ul{\la}}_{G,N}$ (see~\cite[Lemma~4.4]{MV}).
	Now part $(1)$ follows from the definitions and the fact that $\on{dim}{\CW}^{\ul{\nu}}_\mu=\langle2 \rho^{\vee},|\ul{\nu}|-\mu\rangle,\,
	\on{dim}{\CW}^{\nu}_\mu=\langle 2\rho^{\vee},\nu-\mu\rangle$, $\on{dim}\on{Gr}^{\ul{\nu}}_{G,N}=\langle 2\rho^\vee,|\ul{\nu}|\rangle,\, \on{dim}\on{Gr}^{\nu}_G=\langle 2\rho^\vee,\nu \rangle$.
	Let us prove part $(2)$.
    The proof of the fact that $\widetilde{\CW}^{\ul{\la}}_\mu$ is Cohen-Macaulay is the same as the proof of~\cite[Lemma~2.16]{BFN}. Since $\widetilde{\CW}^{\ul{\la}}_\mu$ is Cohen-Macaulay we conclude that to prove that $\widetilde{\CW}^{\ul{\la}}_\mu$ is normal it is enough to show that it is regular in codimension $1$. This immediately follows from the stratified semismallness of $\varpi^{\ul{\la}}_\mu $ and the fact that the open stratum $\CW^\la_\mu \subset \overline{\CW}^\la_\mu$ is smooth with the complement of codimension $\geqslant 2$. 
    Part $(3)$ can be proved using the same methods as in the proof of the main Theorem of~\cite{smooth_slice}, see~\cite[Section~3.2.1]{w} and Remark~\ref{smooth_transv_crit} below.
	\epr

\rem{}\label{smooth_transv_crit}
Part $(3)$ of Proposition~\ref{conv_prop} can be strengthened in the following way.
Recall that $\CW^{\ul{\la}}_\mu \subset \widetilde{\CW}^{\ul{\la}}_\mu$ is the open subscheme consisting of tuples $(\CP_0,\ldots,\CP_N,\sigma_1,\ldots,\sigma_N,\phi)$ such that the pole of $\sigma_i$ at zero has degree $\la_i$. Then the variety $\CW^{\ul{\la}}_\mu$ is smooth without any restrictions on $\ul{\la},\, \mu$.  The proof can be obtained in the same way as the proof of the main theorem of~\cite{smooth_slice}, for the another proof  see~\cite[Section~$2$]{y}.
Note now that for minuscule $\la_i$  we have $\CW^{\ul{\la}}_\mu=\widetilde{\CW}^{\ul{\la}}_\mu$. Moreover we claim that in general $\CW^{\ul{\la}}_\mu \subset \widetilde{\CW}^{\ul{\la}}_\mu$ is precisely the smooth locus. We will give two possible approaches to see this. There is also a third approach that is described in~\cite[Section~$2$]{y}, see~\cite[Proposition $2.4$]{y}. First approach was communicated to us by Dinakar Muthiah.
This approach is very short and nice but uses the Poisson structure on $\widetilde{\CW}^{\ul{\la}}_\mu$ and does not have a chance to work over arbitrary ground field. Another one is a combination of ideas from~\cite[Theorem~2.9]{kwwy},~\cite{smooth_slice} and has an advantage that it (potentially) works in any characteristic. 

Let us explain the first approach. Our goal is to show that the smooth locus $(\widetilde{\CW}^{\ul{\la}}_\mu)^{\mathrm{sm}}$ coincides with $\CW^{\ul{\la}}_\mu$. We already know that  $\CW^{\ul{\la}}_\mu \subset (\widetilde{\CW}^{\ul{\la}}_\mu)^{\mathrm{sm}}$ so we only need to show that $(\widetilde{\CW}^{\ul{\la}}_\mu)^{\mathrm{sm}} \subset \CW^{\ul{\la}}_\mu$. 
Indeed recall that we have a stratification $\widetilde{\CW}^{\ul{\la}}_\mu=\bigsqcup_{\ul{\nu} \leqslant \ul{\la},\mu \leqslant |\ul{\nu}|}\CW^{\ul{\nu}}_\mu$  by Poisson subvarieties. It follows from~\cite[Theorem~5]{w} that $(\widetilde{\CW}^{\ul{\la}}_\mu)^{\mathrm{sm}} \subset \widetilde{\CW}^{\ul{\la}}_\mu$ is symplectic. It follows that $\CW^{\ul{\nu}}_\mu \cap (\widetilde{\CW}^{\ul{\la}}_\mu)^{\mathrm{sm}} \subset (\widetilde{\CW}^{\ul{\la}}_\mu)^{\mathrm{sm}}$ is Poisson. From~\cite[Lemma~1.4]{kal} we conclude that  $\CW^{\ul{\nu}}_\mu \cap (\widetilde{\CW}^{\ul{\la}}_\mu)^{\mathrm{sm}}=\varnothing$ for $\ul{\nu} \neq \ul{\la}$. It follows that $(\widetilde{\CW}^{\ul{\la}}_\mu)^{\mathrm{sm}} \subset \CW^{\ul{\la}}_\mu$.



Let us now sketch the second approach.
Recall (see~\cite{smooth_slice}) that
we can define subfunctors $\mathcal{X}^\la=G[z]z^\la G[z],\, \ol{\mathcal{X}}^\la=\ol{G[z]z^\la G[z]}$ of the functor $G((z^{-1}))$. More generally
we can define functors $\mathcal{X}^{\ul{\la}},\, \widetilde{\mathcal{X}}^{\ul{\la}}$ as follows. 
Consider the convolution thick affine Grassmanian 
\begin{equation*}
\widetilde{\on{\bf{Gr}}}_{G,N}:=\underbrace{G((z^{-1})) \times_{G[z]} \times \ldots \times_{G[z]} G((z^{-1})) \times_{G[z]} \on{{\bf{Gr}}}_G}_{\text{N}}.    \end{equation*}
We have the natural morphism 
\begin{equation*}
\underbrace{G((z^{-1})) \times_{G[z]} \times \ldots \times_{G[z]} G((z^{-1})) \times_{G[z]} G((z^{-1}))}_{\text{N}} \xrightarrow{{\bf{a}}} \widetilde{\on{{\bf{Gr}}}}_{G,N}
\end{equation*}
that is a Zariski locally trivial principal bundle for the group $G[z]$.

Note that we have embeddings 
$
\on{Gr}^{\ul{\la}}_{G,N},\widetilde{\on{Gr}}^{\ul{\la}}_{G,N} \subset \widetilde{\on{Gr}}_{G,N} \subset \widetilde{\on{{\bf{Gr}}}}_{G,N}.   
$

We then define
\begin{equation*}
\mathcal{X}^{\ul{\la}}={\bf{a}}^{-1}({\bf{\on{Gr}}}^{\ul{\la}}_{G,N}),\, \widetilde{\mathcal{X}}^{\ul{\la}}={\bf{a}}^{-1}(\widetilde{{\bf{\on{Gr}}}}^{\ul{\la}}_{G,N}).
\end{equation*}
As in~\cite[Equation~(3.1)]{smooth_slice} we can also define the functor
\begin{equation*}
\mathcal{X}_\mu=U((z^{-1})) z^\mu T[[z^{-1}]]_1 U_-((z^{-1}))=U((z^{-1}))T[[z^{-1}]]_1U_-((z^{-1}))z^\mu \subset G((z^{-1})).
\end{equation*}

We have the multiplication morphism 
\begin{equation*}
\underbrace{G((z^{-1})) \times_{G[z]} \times \ldots \times_{G[z]} G((z^{-1})) \times_{G[z]} G((z^{-1}))}_{\text{N}} \ra G((z^{-1}))
\end{equation*}
and define $\widetilde{\CW}_{\mu}$ as the preimage of $\CW_\mu=U[[z^{-1}]]_1 z^\mu B_-[[z^{-1}]]_1$ with respect to this morphism. We also define $\widetilde{\mathcal{X}}_\mu$ to be the preimage of $\mathcal{X}_\mu$ with respect to this morphism.

We also define functors
\begin{equation*}
\mathcal{X}^{\ul{\la}}_\mu=\mathcal{X}^{\ul{\la}} \cap \widetilde{\mathcal{X}}_{\mu},\, \widetilde{\mathcal{X}}^{\ul{\la}}_\mu=\widetilde{\mathcal{X}}^{\ul{\la}} \cap \widetilde{\mathcal{X}}_{\mu}.
\end{equation*}
Note that $\widetilde{\CW}^{\ul{\la}}_\mu = \widetilde{\CW}_{\mu} \cap \widetilde{\mathcal{X}}^{\ul{\la}} ,\, \CW^{\ul{\la}}_\mu = \widetilde{\CW}_\mu \cap \mathcal{X}^{\ul{\la}}$ (this is the matrix description of convolution diagrams over slices).

By the same reasons as in~\cite[Proposition~3.8]{smooth_slice} we have isomorphisms of functors $\widetilde{\mathcal{X}}_\mu \simeq U[z] \times \widetilde{\CW}_\mu \times U_-[z]$ that induce isomorphisms 
\begin{equation}\label{tilde_X_via_W}
\mathcal{X}^{\ul{\la}}_\mu \simeq U[z] \times \CW^{\ul{\la}}_\mu \times U_-[z],\, \widetilde{\mathcal{X}}^{\ul{\la}}_{\mu} \simeq U[z] \times \widetilde{\CW}^{\ul{\la}}_\mu \times U_-[z].
\end{equation}

Pick now a closed point $p_0 \in \widetilde{\CW}^{\ul{\la}}_\mu(\BC) \setminus \CW^{\ul{\la}}_\mu(\BC)$. Our goal is to show that $\widetilde{\CW}^{\ul{\la}}_\mu$ is not smooth at this point. Since $\widetilde{\CW}^{\ul{\la}}_\mu$ is of finite type then it is enough to show that  $\widetilde{\CW}^{\ul{\la}}_\mu$ is not formally smooth at $p_0$ (i.e. every  open neighbourhood of $p_0 \in \widetilde{\CW}^{\ul{\la}}_\mu$ is not formally smooth). Recall the isomorphism $\widetilde{\mathcal{X}}^{\ul{\la}}_\mu \simeq U[z] \times \widetilde{\CW}^{\ul{\la}}_\mu \times U_-[z]$ from~(\ref{tilde_X_via_W}) and let us denote the point $\{1\} \times \{p_0\} \times \{1\} \in \widetilde{\mathcal{X}}^{\ul{\la}}_\mu$ by $p$.
Recall that we have an embedding $\widetilde{\mathcal{X}}^{\ul{\la}}_\mu \subset \widetilde{\mathcal{X}}^{\ul{\la}}$
and the Zariski locally-trivial principal bundle $\widetilde{\mathcal{X}}^{\ul{\la}} \twoheadrightarrow \widetilde{\on{Gr}}^{\ul{\la}}_{G,N}$. We denote by $p'$ the image of the point $p$ under the composition $\widetilde{\mathcal{X}}^{\ul{\la}}_\mu \subset \widetilde{\mathcal{X}}^{\ul{\la}} \twoheadrightarrow  \widetilde{\on{Gr}}^{\ul{\la}}_{G,N}$. Note now that the smooth locus of $\widetilde{\on{Gr}}^{\ul{\la}}_{G,N}$ is $\on{Gr}^{\ul{\la}}_{G,N}$ (this can be deduced from~\cite[Theorem~0.1(b)]{em}, see also~\cite[Corollary~B]{mov} and~\cite[Proof of Theorem~2.9]{kwwy}) so $p'$ is not a formally smooth point of $\widetilde{\on{Gr}}^{\ul{\la}}_{G,N}$.
So we can find a local Artinian $\BC$-algebra $(A,\mathfrak{m})$ and an element $x \in \widetilde{\on{Gr}}^{\ul{\la}}_G(A)$ that sends $\mathfrak{m} \in \on{Spec}A$ to $p'$ and such that $x$ can not be lifted to an element of $\widetilde{\on{Gr}}^{\ul{\la}}_{G,N}(\tilde{A})$ for some nilpotent extension $\tilde{A} \twoheadrightarrow A$ (here we use~\cite{stack}). 
Since $\widetilde{\mathcal{X}}^{\ul{\la}} \twoheadrightarrow \widetilde{\on{Gr}}^{\ul{\la}}_{G,N}$ is Zariski locally-trivial $G[z]$-torsor, $G[z]$ is formally smooth and $\mathfrak{m}$ is nilpotent we can lift point $x \in \widetilde{\on{Gr}}^{\ul{\la}}_{G,N}(A)$ to some point $y \in \widetilde{\mathcal{X}}^{\ul{\la}}(A)$ that maps $\mathfrak{m}$ to $p \in \widetilde{\mathcal{X}}^{\ul{\la}}(\BC)$.
Recall that the maximal ideal $\mathfrak{m} \subset A$ is nilpotent and the image of $\mathfrak{m}$ in  $\widetilde{\mathcal{X}}^{\ul{\la}}(\BC)$ is equal to $p$ so the image of $\mathfrak{m}$ lies in the intersection $\widetilde{\mathcal{X}}^{\ul{\la}}(\BC) \cap \widetilde{\mathcal{X}}_{\mu}(\BC)$. It can then be deduced from the proof of~\cite[Theorem~3.14]{smooth_slice} that there exists an element $t \in T[A[z]]$ such that $ty \in \widetilde{\mathcal{X}}^{\ul{\la}}(A) \cap \widetilde{\mathcal{X}}_{\mu}(A)$. Note that the image of $t$ in $T[\BC[z]]=T(\BC)$ is equal to $1 \in T$ since if $t_0$ is this image then $p,\, t_0p \in \widetilde{\mathcal{X}}_{\mu}$ that implies $t_0=1$.
Recall now that $\widetilde{\mathcal{X}}^{\ul{\la}} \cap \widetilde{\mathcal{X}}_{\mu} \simeq U[z] \times \widetilde{\CW}^{\ul{\la}}_\mu \times U_-[z]$ so the point $ty$ corresponds to some triple $(u,w,u_-) \in U(A[z]) \times \widetilde{\CW}^{\ul{\la}}_\mu(A) \times U_-(A[z])$ which image in $U(\BC[z]) \times \widetilde{\CW}^{\ul{\la}}_\mu(\BC) \times U_-(\BC[z])$ is $p=\{1\} \times \{p_0\} \times \{1\}$. So the image of $w \in \widetilde{\CW}^{\ul{\la}}_\mu(A)$ in $\widetilde{\CW}^{\ul{\la}}_\mu(\BC)$ is equal to $p_0$. Assume for the sake of contradiction that $\widetilde{\CW}^{\ul{\la}}_\mu$ is smooth at the point $p_0$. Then it follows that there exists a lift of $w$ to some element $\tilde{w} \in \widetilde{\CW}^{\ul{\la}}_\mu(\tilde{A})$. Since $U[z],\, U_-[z]$ are formally smooth we can then lift $u,u_-$ to elements $\tilde{u} \in U[\tilde{A}[z]],\, \tilde{u}_- \in U_-[\tilde{A}[z]]$. Then the triple $(\tilde{u},\tilde{w},\tilde{u}_-)$ corresponds to some lift of the element $ty$ that we denote by $k$. 
Let us now finally
note that $T[z]$ is formally smooth, so we can find $\tilde{t} \in T[\tilde{A}[z]]$ lifting $t$. We can now consider the  product $\tilde{y}:=\tilde{t}^{-1}k \in \widetilde{\mathcal{X}}^{\ul{\la}}(\tilde{A})$ and note that this element lifts $y \in \widetilde{\mathcal{X}}^{\ul{\la}}(A)$. It remains to note that the image of $\tilde{y}$ in $\widetilde{\on{Gr}}^{\ul{\la}}_{G,N}(\tilde{A})$ lifts $x$. This contradiction finishes the proof.

Let us now note that if $\mu$ is dominant then by~\cite[Theorem~2.9]{kwwy} $\ol{\CW}^{\ul{\la}}_{\mu}$ is smooth iff for every decomposition $\mu=\mu_1+\ldots+\mu_N$ with $\mu_i$ being weights of $V^{\la_i}$ we have $\mu_i \in W\la_i$. The implication $\Rightarrow$ is  true in general (without assuming that $\mu$ is dominant). 
The implication $\Leftarrow$ is not true in general. For example we can take $N=1$ and consider any dominant $\la$ that is not minuscule and take $\mu=w_0(\la)$. Then since $\la$ is not minuscule we then can find dominant $\nu < \la$ that is a weight of $V^\la$. It follows then that $\nu \geqslant w_0(\la)$. We conclude that $\ol{\CW}^{\nu}_{w_0(\la)} \neq \varnothing$ so $\ol{\CW}^{\la}_{\mu} \setminus \CW^\la_\mu$ is nonempty since it contains $\ol{\CW}^{\nu}_{w_0(\la)}$. It follows that $\ol{\CW}^{\la}_{w_0(\la)}$ is not smooth. 

By the above results  we see  that $\ol{\CW}^\la_\mu$ is smooth iff $\ol{\CW}^\la_\mu=\CW^\la_\mu$.
So we conclude that $\ol{\CW}^{\ul{\la}}_{\mu}$ is smooth iff for every $N$-tuple of dominant coweights such that $\nu_i \leqslant \la_i$ and $\mu \leqslant |\ul{\nu}|$ we must have $\ul{\nu} = \ul{\la}$.
\erem


\ssec{ap}{Multiplication morphisms for convolution diagrams over slices}	The goal for now is to define multiplication morphisms between convolution diagrams over slices (c.f.~\cite[Section~2.(vi)]{BFN}). Following the approach of~\cite{BFN} we start with a symmetric definition of the variety $\widetilde{\mathcal{W}}^{\ul{\la}}_\mu$ (c.f.~\cite[Section~2(v)]{BFN}). 

\sssec{ap}{Symmetric definition of convolution diagrams over slices}\label{symmetric_slices} Fix a decomposition $\mu=\mu_- + \mu_+$ and define $\widetilde{\mathcal{W}}^{\ul{\la}}_{\mu_-,\mu_+}$ as the moduli space of the following data:
	
	(a) $G$-bundles $\CP_0,\CP_1,\ldots,\CP_N$.
	
	(b) Isomorphisms $\sigma_i\colon \CP_{i-1}|_{\BP^1 \setminus \{0\}} \iso \CP_i|_{\BP^1 \setminus \{0\}}$ having a pole of degree $\leqslant \la_i$ at zero.
	
	(c) A trivialization $s=s_0$ of ${\CP_0}|_\infty$ (we denote by $s_i\colon G \iso {\CP_i}|_\infty$ the trivialization of ${\CP_i}|_\infty$ equal to ${\sigma_i}|_{\infty}\circ\dots\circ {\sigma_1}|_{\infty} \circ s$).
	
	(d) A $B_-$-structure $\phi_-$ on  $\CP_0$ such that the induced $T$-bundle has degree $-w_0\mu_-$ and the fiber of $\phi_-$ at $\infty$ (with respect to our fixed trivialization) is $B \subset G$.
	
	(e) A $B$-structure $\phi_+$ on  $\CP_N$ such that the induced $T$-bundle has degree $w_0\mu_+$ and the fiber of $\phi_+$ at $\infty$ (with respect to our fixed trivialization) is $B_- \subset G$.
	
	\prop{}\label{iso_symm_ord}
	There exists an isomorphism $\widetilde{\CW}^{\ul{\la}}_{\mu_-,\mu_+} \iso \widetilde{\CW}^{\ul{\la}}_\mu$.
	\eprop
	\prf
	The proof repeats the one in~\cite[Section~2(v)]{BFN}. Let us briefly recall how to construct a morphism $\widetilde{\CW}^{\ul{\la}}_{\mu_-,\mu_+} \ra \widetilde{\CW}^{\ul{\la}}_\mu$. Pick $(\CP_i,\sigma_i,\phi_-,\phi_+, s) \in \widetilde{\CW}^{\ul{\la}}_{\mu_-,\mu_+}$. Trivializations $\sigma_j$ define us $B$-structure $\phi^{i}_+$ and $B_-$-structure $\phi^{i}_-$ on $\CP_i|_{\BP^1 \setminus \{0\}}$ for any $i=0,\ldots,N$. 
	Recall that $U_{\infty}$ is a formal neighbourhood of $\infty \in \BP^1$.
	Note that $\phi^i_+$ and $\phi^i_-$ are transversal being restricted to $U_\infty$.
	Let $\CP^T_i$ be the corresponding $T$-bundles on $U_{\infty}$. Consider now the $T$-bundles $'\CP^T_i:=\CP^T_i(w_0\mu_-\cdot \infty)$ and note that the bundles $'\CP^T_i$ and $\CP^T_i$ are isomorphic off $\infty \in U_\infty$. We define $'\CP_i =\,\CP_i(w_0\mu_-)$ as the result of gluing $\CP_i$ and the induced $G$-bundle $'\CP^T_i \times^T G$ in the punctured neighbourhood of $\infty \in U_\infty$. 
	Bundle $\CP_i(w_0\mu_-)$ is called the {\em{Hecke transform}} of $\CP_i$ (w.r.t. $w_0\mu \in \La$).
	The isomorphism $\sigma_i\colon '\CP_{i-1}|_{\BP^1 \setminus \{0,\infty\}} \iso '\CP_i|_{\BP^1 \setminus \{0,\infty\}}$ extends to $\BP^1 \setminus \{0\}$ to some isomorphism $'\sigma_i$ 
	and $\phi_-$ (resp. $\phi_+$) extends from $\BP^1 \setminus \{\infty\}$ to a $B_-$-structure of degree $0$ on $'\CP_0$ (resp. $B$-structure of degree $\mu$ on $'\CP_n$ to be denoted $\phi$). We now send 
	\begin{equation*}
	 \widetilde{\CW}^{\ul{\la}}_{\mu_-,\mu_+}\ni (\CP_i,\sigma_i,\phi_-,\phi_+, s) \mapsto ('\CP_i,'\sigma_i,\phi) \in \widetilde{\CW}^{\ul{\la}}_\mu.	\end{equation*}
    \epr	
    
\sssec{ap}{Multiplication morphism}	
Let us now recall the multiplication morphisms between convolution diagrams over generalized slices. Pick $N,\,N' \in \BZ_{\geqslant 1}$ and fix $N$-tuple of coweights $\ul{\la}=(\la_1,\ldots,\la_N)$, $N'$-tuple of coweights $\ul{\la}'=(\la'_1,\ldots,\la'_{N'})$ and two coweights $\mu,\,\mu' \in \La$. We define the multiplication morphism $\tilde{{\bf{m}}}^{\ul{\la},\ul{\la}'}_{\mu,\mu'}\colon\widetilde{\CW}^{\ul{\la}}_\mu \times \widetilde{\CW}^{\ul{\la}'}_{\mu'} \ra \widetilde{\CW}^{(\ul{\la},\ul{\la}')}_{\mu+\mu'}$ as the composition of the following morphisms: 
	\begin{equation}\label{comp_multt}
	\widetilde{\CW}^{\ul{\la}}_\mu \times \widetilde{\CW}^{\ul{\la}'}_{\mu'} \iso \widetilde{\CW}^{\ul{\la}}_{\mu,0} \times \widetilde{\CW}^{\ul{\la}}_{0,\mu'} \xrightarrow{{\bf{c}}} \widetilde{\CW}^{(\ul{\la},\ul{\la}')}_{\mu,\mu'} \iso \widetilde{\CW}^{(\ul{\la},\ul{\la}')}_{\mu+\mu'},
	\end{equation}
	where the morphism ${\bf{c}}$ is given by
	\begin{equation*}
	\widetilde{\CW}^{\ul{\la}}_{\mu,0} \times \widetilde{\CW}^{\ul{\la}'}_{0,\mu'} \ni (\CP_i,\sigma_i,\phi, s,\CP'_i,\sigma'_i,\phi', s') \mapsto (\CP_0,\ldots,\CP_N=\CP^{triv}=\CP'_0,\ldots,\CP'_{N'},\sigma_i,\sigma'_i,\phi,\phi', s).
	\end{equation*}
In the case when $N=N'=1$, then $\ul{\la}=\la_1,\, \ul{\la}'=\la_1'$ and $\widetilde{\CW}^{\ul{\la}}_{\mu}=\overline{\CW}^{\la_1}_{\mu_1},\, \widetilde{\CW}^{\ul{\la}'}_{\mu'}=\overline{\CW}^{\la'_1}_{\mu'_1}$ we obtain multiplication morphisms 
\begin{equation*}
\tilde{{\bf{m}}}^{\la_1,\la'_1}_{\mu_1,\mu'_1}\colon \overline{\CW}^{\la_1}_{\mu_1} \times \overline{\CW}^{\la_2}_{\mu_2} \ra 
\widetilde{\CW}^{(\la_1,\la_2)}_{\mu_1+\mu_2}.
\end{equation*}
We can consider compositions
\begin{equation}\label{ordin_multipl}
\overline{\CW}^{\la_1}_{\mu_1} \times \overline{\CW}^{\la_2}_{\mu_2} \xrightarrow{\tilde{{\bf{m}}}^{\la_1,\la'_1}_{\mu_1,\mu'_1}}  \widetilde{\CW}^{(\la_1,\la_2)}_{\mu_1+\mu_2} \xrightarrow{\varpi^{(\la_1,\la_2)}_{\mu_1+\mu_2}} \overline{\CW}^{\la_1+\la_2}_{\mu_1+\mu_2}
\end{equation}
to be denoted ${\bf{m}}^{\la_1,\la'_1}_{\mu_1,\mu'_1}$. These morphisms were defined in~\cite[Section~2(vi)]{BFN}.

	We can also define multiplication morphisms 
	\begin{equation*}
	\tilde{{\bf{m}}}^{\ul{\la}^1,\ldots,\ul{\la}^k}_{\mu^1,\ldots,\mu^k}\colon \widetilde{\CW}^{\ul{\la}^1}_{\mu^1} \times \ldots \times \widetilde{\CW}^{\ul{\la}^k}_{\mu^k} \ra \widetilde{\CW}^{(\ul{\la}^1,\ul{\la}^2,\ldots,\ul{\la}^k)}_{\mu^1+\ldots+\mu^k},
	\end{equation*}
	where $\ul{\la}^i=(\la_1^i,\ldots,\la^i_{N_i}),\, \mu^i \in \La,\, N_i \in \BZ_{\geqslant 1}$ as the following composition:
	\begin{multline*}
	 \widetilde{\CW}^{\ul{\la}^1}_{\mu^1}\times\ldots \times \widetilde{\CW}^{\ul{\la}^k}_{\mu^k}  \xrightarrow{\on{Id} \times \tilde{{\bf{m}}}^{\ul{\la}^{k-1},\ul{\la}^k}_{\mu^{k-1},\mu^k}} \widetilde{\CW}^{\ul{\la}^1}_{\mu^1}\times\ldots \times \widetilde{\CW}^{\ul{\la}^{k-2}}_{\mu^{k-2}} \times \widetilde{\CW}^{(\ul{\la}^{k-1},\ul{\la}^k)}_{\mu^{k-1}+\mu^k}  \xrightarrow{\on{Id} \times \tilde{{\bf{m}}}^{\ul{\la}^{k-2},(\ul{\la}^{k-1},\ul{\la}^k)}_{\mu^{k-2},\mu^{k-1}+\mu^{k}}} \ldots \\ 
	 \ldots \xrightarrow{\on{Id}\times \tilde{{\bf{m}}}^{\ul{\la}^2,(\ul{\la}^3,\ldots,\ul{\la}^k)}_{\mu^2,\mu^3+\ldots+\mu^k}} \widetilde{\CW}^{\ul{\la}^1}_{\mu^1} \times \widetilde{\CW}^{(\ul{\la}^2,\ldots,\ul{\la}^{k})}_{\mu^2+\ldots+\mu^k} \xrightarrow{\on{Id}\times \tilde{{\bf{m}}}^{\ul{\la}^1,(\ul{\la}^2,\ldots,\ul{\la}^k)}_{\mu^1,\mu^2+\ldots+\mu^k}} 
	 \widetilde{\CW}^{(\ul{\la}^1,\ul{\la}^2,\ldots,\ul{\la}^k)}_{\mu^1+\ldots+\mu^k}.
	\end{multline*}
	
\rem{}	
{{Let us point out that morphisms 
$\tilde{{\bf{m}}}^{\ul{\la}^1,\ldots,\ul{\la}^k}_{\mu^1,\ldots,\mu^k}$ are not associative in general.
}}
\erem

\prop{}\label{mult_prop} 
Morphism $\tilde{{\bf{m}}}^{\ul{\la},\ul{\la}'}_{\mu,\mu'}$ or more generally $\tilde{{\bf{m}}}^{\ul{\la}^1,\ldots,\ul{\la}^k}_{\mu^1,\ldots,\mu^k}$ is an open embedding.
	\eprop
	\prf
    It is enough to show that the morphism ${\bf{c}}$ in~(\ref{comp_multt}) is an open embedding. Note that the image of ${\bf{c}}$ consists of tuples 
    \begin{equation*}
    (\CP_0,\ldots,\CP_N,\CP'_1,\ldots,\CP'_{N'},\sigma_i,\sigma'_i,\phi,\phi', s, s') \in  \widetilde{\CW}^{\ul{\la},\ul{\la}'}_{\mu,\mu'},  
    \end{equation*}
    such that $\CP_N$ is trivial. This is an open condition, so the image of ${\bf{c}}$ is open.
    Note also that we can construct the inverse morphism 
    $
    \on{im}{\bf{c}} \ra \widetilde{\CW}^{\ul{\la}}_{\mu} \times \widetilde{\CW}^{\ul{\la}'}_{\mu'}
    $ as follows: 
    there exists the unique trivialization $\CP_N \simeq \CP^{triv}$ extending our given trivialization at infinity.
    Now we send 
    $(\CP_0,\ldots,\CP_N,\CP'_1,\ldots,\CP_{N'},\sigma_i,\sigma'_i,\phi,\phi', s, s')$ to 
    \begin{equation*}
    \Big((\CP_0,\ldots,\CP_N\simeq \CP^{triv},\sigma_i,\phi, s),(\CP^{triv} \simeq \CP_N,\CP'_1\ldots,\CP'_N,\sigma'_i,\phi', s')\Big) \in \widetilde{\CW}^{\ul{\la}}_\mu \times \widetilde{\CW}^{\ul{\la}'}_{\mu'}.
    \end{equation*}
    
	\epr
	
\rem{}
{{
Note that it follows from Propositions~\ref{conv_prop},~\ref{mult_prop} that the multiplication morphisms for generalized slices (see~(\ref{ordin_multipl}),~\cite[Section~2(vi)]{BFN}) are birational.}}
\erem

\ssec{ap}{Cartan torus fixed points of convolution slices}	
Note that we have a natural action $T \curvearrowright \widetilde{\CW}^{\ul{\la}}_\mu$.
Let us describe the fixed points of this action.
\prop{ap}\label{torus_fixed_pts}
The set $(\widetilde{\CW}^{\ul{\la}}_\mu)^T$ is parametrized by $N$-tuples $\ul{\mu}=(\mu_1,\ldots,\mu_N)$ such that $\mu=\mu_1+\ldots+\mu_N$ and  $\mu_i$ appears as a weight of $V^{\la_i}$. The point which corresponds to the tuple $\ul{\mu}$ will be denoted by $z^{\ul{\mu}}$ and can be obtained as the image of the point $z^{\mu_1}\times \ldots \times z^{\mu_N} \in \ol{\CW}^{\la_1}_{\mu_1}\times \ldots \times \ol{\CW}^{\la_N}_{\mu_N}$ under the morphism $\tilde{{\bf{m}}}^{\ul{\la}}_{\ul{\mu}}=\tilde{{\bf{m}}}^{\la_1,\ldots,\la_N}_{\mu_1,\ldots,\mu_N}\colon \ol{\CW}^{\la_1}_{\mu_1}\times \ldots \times \ol{\CW}^{\la_N}_{\mu_N} \ra \widetilde{\CW}^{\ul{\la}}_\mu$.
\eprop
\prf
It follows from the definitions that the morphism $\tilde{{\bf{m}}}^{\ul{\la}}_{\ul{\mu}}$ is $T$-equivariant.
Recall now that by Proposition~\ref{fixedpoints} we have $(\overline{\CW}^{\la_i}_{\mu_i})^T=z^{\mu_i}$, so we conclude that the point $(z^{\mu_1},\ldots,z^{\mu_N}) \in \ol{\CW}^{\la_1}_{\mu_1}\times \ldots \times \ol{\CW}^{\la_N}_{\mu_N}$ is $T$-fixed and its image $z^{\ul{\mu}}$ is the $T$-fixed point. It remains to prove that any $T$-fixed point $x \in (\widetilde{\CW}^{\ul{\la}}_{\mu})^T$ is $z^{\ul{\mu}}$ for some $N$-tuple $\ul{\mu}$ as above. 

Recall the convolution affine Grassmannian
$\widetilde{\on{Gr}}_{G,N}$ (see Section~\ref{conv_diagr_sl}).
Note that we have the action $G \curvearrowright \widetilde{\on{Gr}}_G$ induced by the action of $G$ on $\CP^{triv}$ via authomorphisms. In ``matrix" terms (see~(\ref{def_matrix_conv_gr})) the action is induced by the left multiplication $G \curvearrowright G(\CK)$.
We have a forgetful  morphism 
$\tilde{p}^{\ul{\la}}_\mu\colon \widetilde{\CW}^{\ul{\la}}_\mu \ra \widetilde{\on{Gr}}^{\ul{\la}}_{G,N}$ that is clearly $T$-equivariant.

Recall also that we have the natural morphism 
$
\varpi^{\ul{\la}}_\mu\colon \widetilde{\CW}^{\ul{\la}}_\mu \ra \overline{\CW}^\la_\mu
$ (see~(\ref{morp_compose}))
that is also $T$-equivariant. 
We have morphisms 
\begin{equation*}
\pi^{\ul{\la}} \colon \widetilde{\on{Gr}}^{\ul{\la}}_{N,G} \ra \overline{\on{Gr}}{}^\la_G,\, 
(\CP_0,\CP_1,\ldots,\CP_N,\sigma_1,\ldots,\sigma_N) \mapsto (\CP_N,\sigma_N\circ \sigma_{N-1}\circ\ldots\circ\sigma_1),
\end{equation*}
\begin{equation*}p^\la_\mu\colon\overline{\CW}^{\la}_\mu \ra \overline{\on{Gr}}{}^{\la}_G,\,
(\CP,\sigma,\phi) \mapsto (\CP,\sigma)
\end{equation*}
and it follows from the definitions that the compositions $\pi^{\ul{\la}} \circ \tilde{p}^{\ul{\la}}_\mu,\, p^{\la}_\mu \circ \varpi^{\ul{\la}}_\mu$ coincide with the morphism
\begin{equation*}
\widetilde{\CW}^{\ul{\la}}_\mu \ra \ol{\on{Gr}}{}^\la_G,\, 	(\CP_0,\CP_1,\ldots,\CP_N,\sigma_1,\ldots,\sigma_N,\phi) \mapsto (\CP_N, \sigma_n\circ \sigma_{N-1}\circ\ldots\circ\sigma_1).
\end{equation*}
So we obtain a morphism 
\begin{equation*}
\widetilde{\CW}^{\ul{\la}}_\mu \xrightarrow{(\tilde{p}^{\ul{\la}}_\mu,\varpi^{\ul{\la}}_\mu)} \widetilde{\on{Gr}}^{\ul{\la}}_{G,N} \times_{\ol{\on{Gr}}{}^\la_G} \overline{\CW}^\la_\mu
\end{equation*}
that is clearly an embedding.

It follows from the identification~(\ref{def_matrix_conv_gr}) and the definitions that the set $(\widetilde{\on{Gr}}^{\ul{\la}}_{G,N})^T$ consists of the points of the form $[(z^{\mu_1},\ldots,z^{\mu_N})]$ such that $z^{\mu_i} \in \ol{\on{Gr}}{}^{\la_i}_G$ so $\mu_i$ is a weight of $V^{\la_i}$.
It follows from Proposition~\ref{fixedpoints} that the set $(\overline{\CW}^{\la}_\mu)^T$ coincides with $\{z^\mu\}$ if $\mu$ is a weight of $V^\la$ and is empty otherwise. 
We conclude that if $x \in (\widetilde{\CW}^{\ul{\la}}_\mu)^T$, then we must have $\big(\tilde{p}^{\ul{\la}}_\mu,\varpi^{\ul{\la}}_\mu\big)(x)=([(z^{\mu_1},\ldots,z^{\mu_N})],z^\mu)$ and we must have $\mu=\mu_1+\ldots+\mu_N$ since the images of $[(z^{\mu_1},\ldots,z^{\mu_N})],\, z^\mu$ in $\ol{\on{Gr}}{}^\la_G$ should coincide. Since $(\tilde{p}^{\ul{\la}}_\mu,\varpi^{\ul{\la}}_\mu)$ is an embedding it follows that $z^{\ul{\mu}}=x$.
\epr


\ssec{ap}{Loop rotation action and covering}
We are gratefull to Hiraku Nakajima for lots of explanations on the results of this Section.

Let us consider the following $\BC^\times$-action on $\ol{\CW}^{\la_1}_{\mu_1} \times \ldots \times \ol{\CW}^{\la_N}_{\mu_N}$: 
\begin{multline}\label{twisted_action}
(x_1,\ldots,x_N) \mapsto \\
\mapsto \big(tx_1,\mu_1(t^{-1})tx_2,\ldots, (\mu_1+\ldots+\mu_{k-1})(t^{-1})tx_{k},\ldots,(\mu_1+\ldots+\mu_{N-1})(t^{-1})tx_N\big),    
\end{multline}
where $t$ acts via the loop rotation and $(\mu_1+\ldots+\mu_{k-1})(t^{-1}) \in T$ acts via the natural action of $T$.
\lem{}\label{shifted_loop_equiv_mult}
The morphism $\tilde{{\bf{m}}}^{\la_1,\ldots,\la_N}_{\mu_1,\ldots,\mu_N}\colon \ol{\CW}^{\la_1}_{\mu_1} \times \ldots \times \ol{\CW}^{\la_N}_{\mu_N} \ra \widetilde{\CW}^{\ul{\la}}_{\mu}$ is $\BC^\times$-equivariant where $\BC^\times$ acts on 
$\ol{\CW}^{\la_1}_{\mu_1} \times \ldots \times \ol{\CW}^{\la_N}_{\mu_N}$ via~(\ref{twisted_action}) and $\BC^\times$ acts on $\widetilde{\CW}^{\ul{\la}}_{\mu}$ via the loop rotation. 
\elem
\prf
For any $\la'_i \in \La^{+},\, \mu' \leqslant \la'_1+\ldots+\la'_{N'},\, N' \in \BZ_{\geqslant 1}$  we have $\BC^\times \times T$-equivariant birational (by Proposition~\ref{conv_prop} $(1)$)  morphisms $\varpi^{\la'_1,\ldots,\la'_{N'}}_{\mu'}\colon \widetilde{\CW}^{\la'_1,\ldots,\la'_{N'}}_{\mu'} \ra \ol{\CW}^{\la'_1+\ldots+\la'_{N'}}_{\mu'}$, so it follows from the continuity argument that it is enough to show that the multiplication morphism ${\bf{m}}^{\la_1,\ldots,\la_N}_{\mu_1,\ldots,\mu_N}\colon \ol{\CW}^{\la_1}_{\mu_1} \times \ldots \times \ol{\CW}^{\la_N}_{\mu_N} \ra \ol{\CW}^{\la_1+\ldots+\la_N}_{\mu_1+\ldots+\mu_N}$ is $\BC^\times$-equivariant with respect to the twisted action~(\ref{twisted_action}) on $\ol{\CW}^{\la_1}_{\mu_1} \times \ldots \times \ol{\CW}^{\la_N}_{\mu_N}$ and the loop rotation action on $\ol{\CW}^{\la_1+\ldots+\la_N}_{\mu_1+\ldots+\mu_N}$. Let us consider the case $N=2$, the general case then follows by the induction. Let us use the ``matrix" descriptions of slices $\ol{\CW}^{\la_1}_{\mu_1},\, \ol{\CW}^{\la_2}_{\mu_2},\, \ol{\CW}^{\la_1+\la_2}_{\mu_1+\mu_2}$ and of the morphism 
${\bf{m}}^{\la_1,\la_2}_{\mu_1,\mu_2}\colon \ol{\CW}^{\la_1}_{\mu_1} \times \ol{\CW}^{\la_2}_{\mu_2} \ra \ol{\CW}^{\la_1+\la_2}_{\mu_1+\mu_2}$. Recall that 
    \begin{equation*}
   \ol{\CW}^{\la_i}_{\mu_i} \simeq \big(U[[z^{-1}]]_{1}z^{\mu_i}T[[z^{-1}]_1U_{-}[[z^{-1}]]_{1}\cap \ol{G[z]z^{\la_i}G[z]}\big),\, i=1,2,
    \end{equation*}
    \begin{equation*}
   \ol{\CW}^{\la_1+\la_2}_{\mu_1+\mu_2} \simeq \big(U[[z^{-1}]]_{1}z^{\mu_1+\mu_2}T[[z^{-1}]]_1U_{-}[[z^{-1}]]_{1}\cap \ol{G[z]z^{\la_1+\la_2}G[z]}\big)
    \end{equation*}
and by Lemma~\ref{action_via_matrix} the loop rotation actions on slices $\ol{\CW}^{\la_1}_{\mu_1},\, \ol{\CW}^{\la_2}_{\mu_2},\, \ol{\CW}^{\la_1+\la_2}_{\mu_1+\mu_2}$ in ``matrix" terms are given by 
\begin{equation*}
x_1(z) \mapsto x_1(t^{-1}z)t^{\mu_1},\, x_2(z) \mapsto x_2(t^{-1}z)t^{\mu_2},\, 
x(z) \mapsto x(t^{-1}z)t^{\mu_1+\mu_2}
\end{equation*}
respectively
and $T$-actions are given by conjugation, here $x_i \in \ol{\CW}^{\la_i}_{\mu_i},\, i=1,2,\, \, x \in \ol{\CW}^{\la_1+\la_2}_{\mu_1+\mu_2}$.
Recall also that the morphism    ${\bf{m}}^{\la_1,\la_2}_{\mu_1,\mu_2}$ sends $(x_1,x_2) \in \ol{\CW}^{\la_1}_{\mu_1} \times \ol{\CW}^{\la_2}_{\mu_2}$ to 
$\psi([x_1x_2]) \in \ol{\CW}^{\la_1+\la_2}_{\mu_1+\mu_2}$, where $\psi$ is the isomorphism 
\begin{equation*}
\psi\colon U[z] \backslash U((z^{-1}))z^{\mu_1+\mu_2}T[[z^{-1}]]_1U_{-}((z^{-1}))/U_-[z]  \iso U[[z]]_{1}z^{\mu_1+\mu_2}T[[z^{-1}]]_1U_{-}[[z^{-1}]]_{1}.    
\end{equation*}
The following chain of equalities finishes the proof
\begin{multline*}
{\bf{m}}^{\la_1,\la_2}_{\mu_1,\mu_2}(t\cdot (x_1,x_2))={\bf{m}}^{\la_1,\la_2}_{\mu_1,\mu_2}(x_1(t^{-1}z)t^{\mu_1},t^{-\mu_1}x_2(t^{-1}z)t^{\mu_1+\mu_2})=\\
=
\psi([x_1(t^{-1}z)x_2(t^{-1}z)t^{\mu_1+\mu_2}])=\psi([x_1x_2])(t^{-1}z)t^{\mu_1+\mu_2}=t \cdot {\bf{m}}^{\la_1,\la_2}_{\mu_1,\mu_2}(x_1,x_2).
\end{multline*}
\epr

\rem{}
{{
In order to check the equality $\psi([x_1(t^{-1}z)x_2(t^{-1}z)t^{\mu_1+\mu_2}])=\psi([x_1x_2])(t^{-1}z)t^{\mu_1+\mu_2}$ 
we decompose 
\begin{equation*}
x_1x_2t^{\mu_1+\mu_2}=u_1(z)u_2(z)z^{\mu_1+\mu_2}h(z)u_{-,2}(z)u_{-,1}(z)
\end{equation*}
with $u_1(z) \in U[z],\, u_2(z) \in U[[z^{-1}]]_1,\, h(z) \in T[[z^{-1}]]_1,\, u_{-,1}(z) \in U_-[z],\, u_{-,2}(z) \in U_-[[z^{-1}]]_1$ and note that 
\begin{multline*}
x_1(t^{-1}z)x_2(t^{-1}z)t^{\mu_1+\mu_2}=   \\
=
u_1(t^{-1}z)u_2(t^{-1}z)z^{\mu_1+\mu_2}t^{-\mu_1-\mu_2}h(t^{-1}z)u_{-,2}(t^{-1}z)u_{-,1}(t^{-1}z)t^{\mu_1+\mu_2}=\\
=
u_1(t^{-1}z)u_2(t^{-1}z)z^{\mu_1+\mu_2}h(t^{-1}z)(t^{-\mu_1-\mu_2}u_{-,2}(t^{-1}z)t^{\mu_1+\mu_2})(t^{-\mu_1-\mu_2}u_{-,1}(t^{-1}z)t^{\mu_1+\mu_2}),
\end{multline*}
so 
\begin{multline*}
\psi([x_1(t^{-1}z)x_2(t^{-1}z)t^{\mu_1+\mu_2}])=  
u_2(t^{-1}z)z^{\mu_1+\mu_2}h(t^{-1}z)(t^{-\mu_1-\mu_2}u_{-,2}(t^{-1}z)t^{\mu_1+\mu_2})=\\
\psi([x_1x_2])(t^{-1}z)t^{\mu_1+\mu_2}.
\end{multline*}
}}
\erem

\cor{}\label{global_shifted}
The morphism $\tilde{{\bf{m}}}^{\ul{\la}}_{\ul{\mu}}\colon \ol{\CW}^{\la_1}_{\mu_1} \times \ldots \times \ol{\CW}^{\la_N}_{\mu_N} \ra  \widetilde{\CW}^{\ul{\la}}_\mu$ is $T \times \BC^\times$-equivariant, where $T \times \BC^\times$ acts on $\widetilde{\CW}^{\ul{\la}}_\mu$ standardly and the action of $(h,t) \in T \times \BC^\times$ on $\ol{\CW}^{\la_1}_{\mu_1} \times \ldots \times \ol{\CW}^{\la_N}_{\mu_N}$ is given by the following formula:
\begin{multline*}
(h,t) \cdot (x_1,\ldots,x_N) = \\
=\big(htx_1,h\mu_1(t^{-1})tx_2,\ldots, h(\mu_1+\ldots+\mu_{k-1})(t^{-1})tx_{k},\ldots,h(\mu_1+\ldots+\mu_{N-1})(t^{-1})tx_N\big).   
\end{multline*}
\ecor
\prf
Follows from Lemma~\ref{shifted_loop_equiv_mult} together with the fact that the morphism $\tilde{m}^{\ul{\la}}_{\ul{\mu}}$ is $T$-equivariant.
\epr



\th{}\label{covering_figovering}
Let $\mu \in \La$ be such that $\langle\al^{\vee},\mu \rangle \geqslant -1$ for any positive root $\al^{\vee} \in \Delta_+^{\vee}$ and $\la \in \La^+$.
Let $\ul{\la}=(\la_1,\ldots,\la_N)$ be an $N$-tuple of dominant coweights such that $\la_1+\ldots+\la_N=\la$.
Then the variety $\widetilde{\CW}^{\underline{\la}}_\mu$ can be covered by the images of the open emeddings $\tilde{{\bf{m}}}^{\ul{\la}}_{\ul{\mu}}\colon \overline{\CW}^{\la_1}_{\mu_1} \times 
\overline{\CW}^{\la_2}_{\mu_2} \times \ldots \times \overline{\CW}^{\la_N}_{\mu_N} \hookrightarrow \widetilde{\CW}^{\ul{\la}}_\mu$ with $\mu_i$ being weights of $V^{\la_i}, i=1,2,\ldots,N$.
\eth
\prf
Recall that by Proposition~\ref{torus_fixed_pts} we have
\begin{equation*}
\big(\widetilde{\CW}^{\ul{\la}}_\mu\big)^T=\{z^{\ul{\mu}}\,|\,\ul{\mu}=(\mu_1,\ldots,\mu_N),\, \mu_1+\ldots+\mu_N=\mu,\, V^{\la_i}_{\mu_i} \neq 0\},
\end{equation*}
where $V^{\la_i}_{\mu_i}$ is the $\mu_i$ weight space of $V^{\la_i}$.
Recall that by Proposition~\ref{loop_contracts} the loop rotation action contracts  $\ol{\CW}^{\la}_\mu$  to the orbit $U_- \cdot z^\mu$ and the $\BC^\times$-action via $(-2\rho)\colon \BC^\times \ra T$ contracts $U_- \cdot z^\mu$ to the point $z^\mu$. Pick $d \gg 0$ and consider the following action 
\begin{equation}\label{actactact}
\BC^\times \curvearrowright \widetilde{\CW}^{\ul{\la}}_\mu,\, t \cdot x:=(-2\rho)(t)t^dx.  
\end{equation}
Using that the morphism $\varpi^{\ul{\la}}_\mu \colon \widetilde{\CW}^{\ul{\la}}_\mu \ra \ol{\CW}^\la_\mu$ is projective (hence, proper) we conclude that for $d$ large enough we have 
\begin{equation*}
\big(\widetilde{\CW}^{\ul{\la}}_\mu\big)^{\BC^\times}=\{z^{\ul{\mu}}\,|\,\ul{\mu}=(\mu_1,\ldots,\mu_N),\, \mu_1+\ldots+\mu_N=N,\, V^{\la_i}_{\mu_i} \neq 0\}
\end{equation*}
and the $\BC^\times$-action contracts $\widetilde{\CW}^{\ul{\la}}_\mu$ to $\big(\widetilde{\CW}^{\ul{\la}}_\mu\big)^{\BC^\times}$. 

Pick now a point $x \in \widetilde{\CW}^{\ul{\la}}_\mu$. Let $z^{\ul{\mu}} \in \big(\widetilde{\CW}^{\ul{\la}}_\mu\big)^{\BC^\times}$ be the point to which $x$ flows when $t \ra 0$. Recall the multiplication morphism 
$\tilde{{\bf{m}}}^{\ul{\la}}_{\ul{\mu}}\colon \ol{\CW}^{\la_1}_{\mu_1} \times \ldots \times \ol{\CW}^{\la_N}_{\mu_N} \hookrightarrow \widetilde{\CW}^{\ul{\la}}_\mu$. 
If follows from Corollary~\ref{global_shifted} that there exists a $\BC^\times$-action on $\ol{\CW}^{\la_1}_{\mu_1} \times \ldots \times \ol{\CW}^{\la_N}_{\mu_N}$ which makes a morphism $\tilde{{\bf{m}}}^{\ul{\la}}_{\ul{\mu}}$ equivariant (the $\BC^\times$-action on $\widetilde{\CW}^{\ul{\la}}_\mu$ is given by~(\ref{actactact})). 
Since $\tilde{{\bf{m}}}^{\ul{\la}}_{\ul{\mu}}$ is an open embedding and $\underset{t \ra 0}{\on{lim}}\,t\cdot x = z^{\ul{\mu}} \in \on{im}\big(\tilde{{\bf{m}}}^{\ul{\la}}_{\ul{\mu}}\big)$ we conclude that there exists $t_0 \in \BC^\times$ such that $t_0 \cdot x \in \on{im}\big(\tilde{{\bf{m}}}^{\ul{\la}}_{\ul{\mu}}\big)$. The $\BC^\times$-equivariance of $\tilde{{\bf{m}}}^{\ul{\la}}_{\ul{\mu}}$ now implies that  $x \in \on{im}\big(\tilde{{\bf{m}}}^{\ul{\la}}_{\ul{\mu}}\big)$.
\epr

Let us assume now that $\la_i$ are minuscule. It follows from Proposition~\ref{conv_prop} that the variety $\widetilde{\CW}^{\ul{\la}}_\mu$ is smooth.
As a direct corollary of the results above we obtain the following Theorem.
\th{ap}\label{theo_covering_min}
Assume that  $\la_i,i=1,\ldots,N$ are minuscule and $\mu$ is such that $\langle \al^{\vee},\mu \rangle \geqslant -1$ for every $\al^{\vee} \in \Delta_+^{\vee}$,
$\la:=\la_1+\ldots+\la_N$. Then there is an open cover of $\widetilde{\CW}^{\ul{\la}}_\mu$ by open subsets $O_{\ul{\mu}}$ parametrized by $N$-tuples $\ul{\mu}=(\mu_1,\ldots,\mu_N),\, \mu=\mu_1+\ldots+\mu_N,\, \mu_i \leqslant \la_i$. Each $O_{\ul{\mu}}$ is isomorphic to the affine space $\BA^{\langle 2\rho^{\vee},\la-\mu\rangle}$ and contains exactly one $T$-fixed point $z^{\ul{\mu}} \in \widetilde{\CW}^{\ul{\la}}_\mu$. Variety $O_{\ul{\mu}}$ is the image of the multiplication morphism $\tilde{{\bf{m}}}^{\ul{\la}}_{\ul{\mu}}\colon\CW^{\la_1}_{\mu_1}\times \ldots \times \CW^{\la_N}_{\mu_N} \hookrightarrow \widetilde{\CW}^{\ul{\la}}_\mu$.
\eth

\begin{Rem}[Characters of tangent spaces at torus fixed points of $\widetilde{\CW}^{\ul{\la}}_\mu$]\label{characters_resolv}
{{
Note that assuming that $\la_i$ are minuscule and without any restrictions on $\mu$ we can compute $T \times \BC^\times$-characters of $T_{z^{\ul{\mu}}}\widetilde{\CW}^{\ul{\la}}_{\mu}$ of tangent spaces to $T \times \BC^\times$-fixed points $z^{\ul{\mu}} \in (\widetilde{\CW}^{\ul{\la}}_{\mu})^{T\times \BC^\times}$. 
Indeed it follows from Proposition~\ref{character_minis} together with Corollary~\ref{global_shifted} that
\begin{equation*}
\on{ch}_{T}T_{z^{\ul{\mu}}}\widetilde{\CW}^{\ul{\la}}_{\mu}=\sum_{i=1}^N\sum_{\al^{\vee} \in \Delta^{\vee}_-,\, \langle \al^{\vee},\mu_i \rangle>0}\big(e^{\al^{\vee}}+e^{-\al^{\vee}}\big)    
\end{equation*}
\begin{equation}\label{tangent_charchar}
\on{ch}_{T \times \BC^\times}T_{z^{\ul{\mu}}}\widetilde{\CW}^{\ul{\la}}_{\mu}=    
\sum_{i=1}^N\sum_{\al^{\vee} \in \Delta^{\vee}_-,\, \langle \al^{\vee},\mu_i \rangle>0}\big(\hbar^{-\langle\al^{\vee},\mu_1+\ldots+\mu_{i-1}\rangle} e^{\al^{\vee}}+\hbar^{1+\langle\al^{\vee},\mu_1+\ldots+\mu_{i-1}\rangle}e^{-\al^{\vee}}\big) 
\end{equation}

Let us consider the following example. Let $G=\on{GL}_2$ and $\ul{\la}=\underbrace{(\omega_1,\ldots,\omega_1)}_{N},\, \mu=N\omega_1-k\al_1$, here $\omega_1=(1,0)$ is the fundamental coweight, $\al_1=(1,-1)$ is the simple positive coroot and $N \in \BZ_{\geqslant 0},\, 0 \leqslant k \leqslant N$ (in this case $\widetilde{\CW}^{\ul{\la}}_{\mu}$ is symplectically dual to $T^*\on{Gr}(k,N)$). Then the set of fixed points $\big(\widetilde{\CW}^{\ul{\la}}_{\mu}\big)^T$ identifies with the set of subsets 
$\{i_1,\ldots,i_k\} \subset \{1,2,\ldots,N\}$ via the map sending $\{i_1,\ldots,i_k\}$ to the fixed point $z^{\ul{\mu}}$ (see Proposition~\ref{torus_fixed_pts}), where $\mu_{i}=(0,1)$ for $i \in \{i_1,\ldots,i_k\}$ and $\mu_i=(1,0)$ otherwise. It then easily follows that 
\begin{equation}\label{example_char}
{\on{ch}_{T\times \BC^\times}T_{z^{\ul{\mu}}}\widetilde{\CW}^{\ul{\la}}_{\mu}}=\sum_{l=1}^k\big(\hbar^{i_l-2l+1}e^{-\al_1^{\vee}}+\hbar^{2l-i_l}e^{\al_1^{\vee}}\big).
\end{equation}
Note that for $k=1$ the variety $\widetilde{\CW}^{\ul{\la}}_\mu$ is the resolution of $A_{N-1}$-singularity that is a toric variety and the character ${\on{ch}_{T\times \BC^\times}T_{z^{\ul{\mu}}}\widetilde{\CW}^{\ul{\la}}_{\mu}}$ can be computed via the corresponding toric diagram.

This example shows us that it is not true in general that every open subset $O_{\ul{\mu}} \subset \widetilde{\CW}^{\ul{\la}}_\mu$ of Theorem~\ref{theo_covering_min} can be obtained as the attracting locus to the fixed point $z^{\ul{\mu}} \in \ol{\CW}^{\ul{\la}}_\mu$ with respect to the $\BC^\times$-action via some cocharacter $\eta\colon \BC^\times \ra T \times \BC^\times$ such that $\big(\widetilde{\CW}^{\ul{\la}}_{\mu}\big)^{\eta(\BC^\times)}=\{z^{\ul{\mu}} \in \widetilde{\CW}^{\ul{\la}}_{\mu}\}$. Indeed take $\ul{\la},\, \mu$ as in the example. The condition that $\langle \al^\vee,\mu \rangle \geqslant -1$ for every $\al^\vee \in \Delta^\vee_+$ is equivalent to $N \geqslant 2k-1$. Taking $N=3,\,k=2$ and $\mu_1=(1,0),\,\mu_2=\mu_3=(0,1)$ we get
\begin{equation}\label{char_exam_contr}
{\on{ch}_{T\times \BC^\times}T_{z^{\ul{\mu}}}\widetilde{\CW}^{\ul{\la}}_{\mu}}=\hbar e^{-\al_1^{\vee}}+e^{\al_1^{\vee}}+e^{-\al_1^\vee}+\hbar e^{\al_1^\vee}. 
\end{equation}
We conclude that if $\eta$ as above exists then it must be trivial along $T$ (since both $e^{\al_1^\vee}$ and $e^{-\al^{\vee}_1}$ appear in the character) but this is impossible since the fixed point set of $\widetilde{\CW}^{\ul{\la}}_\mu$ with respect to the loop rotation action is non discrete because by~(\ref{char_exam_contr}) 
$T_{z^{\ul{\mu}}}\widetilde{\CW}^{\ul{\la}}_{\mu}$ has nonzero weight zero component (with respect to the loop rotation action).
}}
\end{Rem}

\begin{Rem}[Poincar\'e polynomials of $\widetilde{\CW}^{\ul{\la}}_{\mu}$.]{}\label{Poincare_pol}
We can compute Poincar\'e polynomials of the varieties $\widetilde{\CW}^{\ul{\la}}_{\mu}$ for $\mu$ such that $\langle\al^{\vee},\mu\rangle \geqslant -1$ for $\al^{\vee} \in \Delta^{\vee}_+$. We pick a cocharacter of $T \times \BC^\times$ given by $t \mapsto ((-2\rho(t)),t^d)$ for $d \gg 0$ and consider the corresponding action $\BC^\times \curvearrowright \widetilde{\CW}^{\ul{\la}}_{\mu}$. 
It follows from~(\ref{tangent_charchar}) that the dimension of the attracting part of $T_{z^{\ul{\mu}}}\widetilde{\CW}^{\ul{\la}}_{\mu}$ at the fixed point $z^{\ul{\mu}} \in (\widetilde{\CW}^{\ul{\la}}_{\mu})^{\BC^\times}$ equals to 
\begin{equation*}
\sum_{i=1}^n|\Delta_{\mu_i,-}^{\vee}|+\sum_{i=1}^{n}\big|\{\al^{\vee} \in \Delta_{\mu_i,-}^{\vee},\, \langle \al^{\vee},\mu_1+\ldots+\mu_{i-1}\rangle=-1\}\big|   
\end{equation*}
so using Bialynicki-Birula decomposition we obtain 
\begin{multline*}
\sum_{n \in \BZ_{\geqslant 0}}\on{dim}H^n_c(\widetilde{\CW}^{\ul{\la}}_{\mu})q^n=\\
=\sum_{\ul{\mu},\,|\ul{\mu}|=\mu,\, \mu_i \in W\la_i}q^{2\big(\sum_{i=1}^n|\Delta_{\mu_i,-}^{\vee}|+\sum_{i=1}^{n}\big|\{\al^{\vee} \in \Delta_{\mu_i,-}^{\vee},\, \langle \al^{\vee},\mu_1+\ldots+\mu_{i-1}\rangle=-1\}\big|\big)}.
\end{multline*}
\end{Rem}

\ssec{ap}{Thick affine Grassmannian 
}\label{thick}
In the rest of this Section we prove that without any assumptions on $\la_i$ and $\mu$ the images of $\tilde{{\bf{m}}}^{\ul{\la}}_{\ul{\mu}}$ such that $\mu_i$ are weights of $V^{\la_i}$ for all $i=1,2,\ldots,N-1$ cover $\widetilde{\CW}^{\ul{\la}}_\mu$ (note that there are no conditions on $\mu_N$). We start from some recollections.




Let
$\Gr_G:=G((z^{-1}))/G[z]$ be the {\em{thick}} affine 
Grassmannian 
of $G$.
For $\nu \in \La$, we denote by the same symbol $z^\nu$ the corresponding point of 
$\Gr_G$. 

We have the left action $G((z^{-1})) \curvearrowright 
\Gr_G$. 
We set $I_-:=\on{ev}_{\infty}^{-1}(U)$, where 
$\on{ev}_\infty\colon G[[z^{-1}]] \ra G$ is the evaluation at infinity. Note that $I_-$ is the pro-unipotent group.
The following proposition holds by~\cite[Proposition~1.3.1]{kt}.
\prop{}\label{decomp_thick}   
\begin{enumerate}
	\item The $I_-$-orbits on $\Gr_G$ are in bijection with
	$
	\La$ via 
	\begin{equation*}
	\la \mapsto I_- \cdot z^\nu=:{\bf{\Omega}}_\nu.
	\end{equation*}
	\item
	The $G[[z^{-1}]]$-orbits on $\Gr_G$ are in bijection with  $\La^+$ via 
	\begin{equation*}
	z^\nu \mapsto G[[z^{-1}]] \cdot z^\nu=:{\bf{Gr}}_\nu.
	\end{equation*}
	\end{enumerate}
\eprop

We will also need the following well known lemma 
about the (thin) affine Grassmannian $\on{Gr}_G$. For $\mu \in \La^+$ we set $\on{Gr}_{G,\mu}:=G[z^{-1}] \cdot z^{\mu} \subset \on{Gr}_G$.
\lem{}\label{type_G(O)}
For $\la,\,\mu \in \La^+$ we have $\overline{\on{Gr}}{}^\la_G \cap \on{Gr}_{G,\mu} \neq \varnothing$ iff $\mu \leqslant \la$.
\elem
\prf
Let us prove the implication $\Rightarrow$.
Assume that $\overline{\on{Gr}}{}^\la_G \cap \on{Gr}_{G,\mu} \neq \varnothing$. Pick $x \in \overline{\on{Gr}}{}^\la_G \cap \on{Gr}_{G,\mu}$. Recall that the loop rotation action contracts $\on{Gr}_{G,\mu}$ to $Gz^\mu$. Note now that $\overline{\on{Gr}}{}^\la_G \subset \on{Gr}_G$ is closed, so we must have $gz^\mu \in \overline{\on{Gr}}{}^\la_G$ for some $g \in G$, i.e., $z^\mu \in \overline{\on{Gr}}{}^\la_G$. It follows that $\mu \leqslant \la$.
The implication $\Leftarrow$ is clear since for dominant $\mu  \leqslant \la$ we must have $z^\mu \in \overline{\on{Gr}}{}^\la_G$. 
\epr

Recall that any $G$-bundle $\CP$ on $\BP^1$ is isomorphic to $z^\mu \in \on{Bun}_G$ for some $\mu \in \La^+$. We will say that $\CP$ has {\em{type}} $\mu$.
\cor{}\label{type_orbit_cor}
If $x \in \overline{\on{Gr}}{}^\la_G$ and $\CP$ is the corresponding $G$-bundle, then the type $\mu \in \La^+$ of $\CP$ is $\leqslant \la$.
\ecor
\prf
Follows from Lemma~\ref{type_G(O)} together with the fact that $\on{Gr}_{G,\mu} \subset \on{Gr}_G$ precisely consists of pairs $('\CP,\,'\sigma) \in \on{Gr}_G$ such that $'\CP$ has type $\mu$.
\epr

Let us now prove the following lemma. Note that directly from the definitions we have ${\bf{\Omega}}_0={\bf{Gr}}_0$.

\lem{}\label{mult_grassm} The following holds.
\begin{enumerate}
\item
For every $\nu \in \La$ we have  $z^{-\nu}\cdot {\bf{\Omega}}_{\nu} \subset {\bf{\Omega}}_0={\bf{Gr}}_0$.
\item
For every $x \in \Gr_G$, there exists  $\nu \in \La$ such that $z^{-\nu} x \in {\bf{Gr}}_0$.
\end{enumerate}
\elem
\prf
Note that the group $I_- \cap z^{\nu}G[[z^{-1}]]_1z^{-\nu}$ acts freely and transitively on ${\bf{\Omega}}_\nu$ (this is the standard result, it follows from the fact that $I_-$ is pro-unipotent and so the claim can be checked at the level of Lie algebras, compare with the proof of part $(2)$ of Proposition~\ref{fibration} and~\cite[Lemma~4.5.7 and Corollary~4.5.8]{kash}).
It follows that for every $x \in {\bf{\Omega}}_\nu$ there exists $x_- \in G[[z^{-1}]]_1 G[z]/G[z]$ such that $x=z^{\nu} x_-$.
We conclude that $z^{-\nu} x = x_- \in {\bf{\Omega}}_0$ and part $(1)$ follows. 

Part $(2)$ follows from part $(1)$ together with Proposition~\ref{decomp_thick}: indeed by  Proposition~\ref{decomp_thick} for every $x \in {\bf{Gr}}_G$ there exists $\nu \in \La$ such that $x \in {\bf{\Omega}}_{\nu}$. Then from part $(1)$ we get $z^{-\nu}x \in {\bf{Gr}}_0$.
\epr

\rem{}
{\em{Lemma~\ref{mult_grassm} should be compared to \cite[Corollary~4.5.5 and Definition~4.5.6]{kash} in the paper of Kashiwara where he defined the (thick) affine flag variety. Lemma~\ref{mult_grassm} tells us that ${\bf{Gr}}_G=\bigcup_{\nu \in \La}z^{\nu}{\bf{Gr}}_0$ i.e. that ${\bf{Gr}}_G$ can be covered by open (infinite dimensional) cells $z^{\nu}{\bf{Gr}}_0$.
}}
\erem


\ssec{ap}{Hecke transformations and  covering} 
Pick a point 
\begin{equation*}
(\CP^{triv}=\CP_0,\CP_1,\ldots,\CP_N,\sigma_1,\ldots,\sigma_N,\phi, s) \in \widetilde{\CW}^{\underline{\la}}_\mu.
\end{equation*}
Recall that the bundles $\CP_0,\CP_1,\ldots,\CP_N$ are identified on $\BP^1 \setminus \{0\}$. Let us denote the corresponding bundle by $\CP$. Recall also that $\CP$ has $B$ and $B_-$-structures that are transversal around $\infty \in \BP^1$, so they define a reduction of $\CP$ to a $T$-bundle after restricting it to the formal neighbourhood of $\infty \in \BP^1$ which we denote by $U_\infty$.
We denote by $\CP^T$ the corresponding $T$-bundle on $U_{\infty}$. Note that we have the fixed identification $\CP^T_{\infty} \simeq T$, i.e., a trivialization of $\CP^T$ at $\infty$. Consider now the triple $(\CP_1, \CP^T, \sigma_{U_\infty})$, where by $\sigma_{U_\infty}$ we denote the isomorphism between $\CP^G:=\on{Ind}_T^G\CP^T$ and $\CP_1|_{U_\infty}$. Let us now fix any trivialization of $\CP^T$ which is $\on{Id}$ at $\infty$ and note that now the triple $(\CP_1, \CP^T, \sigma_{U_\infty})$ defines us a point of the {\em{thick}} affine Grassmannian ${\bf{Gr}}_G=G((z^{-1}))/G[z]$ to be denoted $x \in {\bf{Gr}}_G$. 

Our goal (see Lemma~\ref{Hecke_to_trans} bellow) is to 
show that there exists a Hecke transformation of $x$ such that the vector bundle $\CP_1$ becomes trivial, 
i.e., the corresponding point of ${\bf{Gr}}_G$ lies in $G[[z^{-1}]]\cdot 1=\Gr_{0}$. Recall that Hecke transformations are parametrized by $\nu \in \La$ and we denote by $\CP_1(\nu)$ the corresponding Hecke transform (see for example~\cite[Section~3.1]{bg} for some discussion of Hecke transforms). Note that Hecke transformations in this language are nothing else but the left multiplications by $z^{\nu},\, \nu \in \La$.

\lem{}\label{Hecke_to_trans}
Pick a point $(\CP^{triv}=\CP_0,\CP_1,\ldots,\CP_N,\sigma_1,\ldots,\sigma_N,\phi, s) \in \widetilde{\CW}^{\underline{\la}}_\mu$. Then there exists a weight $\mu_1 \in \La$ of $V^{\la_1}$ such that the Hecke transform  $\CP_1(-\mu_1)$ is trivial.
\elem
\prf
Consider a point $(\CP_1,\CP^T,\sigma_{U_\infty}) \in \Gr_G$ as above to be denoted $x$. 
Recall that $(\CP_1,\sigma_1) \in \overline{\on{Gr}}{}^{\la_1}_G$, so the vector bundle $\CP_1$ has type $\mu_1^+ \in \La^+$ such that $\mu_1^+$ is a weight of $V^{\la_1}$ (follows from Corollary~\ref{type_orbit_cor}). 
It follows from Lemma~\ref{mult_grassm} that there exists $\mu_1 \in \La$ such that $z^{-\mu_1} x \in \Gr_0$ and $x \in {\bf{\Omega}}_{\mu_1}$. We conclude that $W\mu_1=W\mu_1^+$ (since $\CP_1$ has type $\mu^+$ and $(\CP_1,\CP^T,\sigma_{U_{\infty}})=x \in {\bf{\Omega}}_{\mu_1}$). It follows that  $\mu_1$ is a weight of $V^{\la_1}$.  It remains to note that the condition $z^{-\mu_1}x \in {\bf{Gr}}_0$ exactly means that the vector bundle $\CP_1(-\mu_1)$ is trivial.
\epr

We can finally formulate and prove the last Theorem of this Section. 
\th{}\label{covering_general}
Variety $\widetilde{\CW}^{\underline{\la}}_\mu$ can be covered by $\overline{\CW}^{\la_1}_{\mu_1} \times 
\overline{\CW}^{\la_2}_{\mu_2} \times \ldots \times \overline{\CW}^{\la_N}_{\mu_N}$ with $\mu_i$ being weights of $V^{\la_i}$ for $i=1,2,\ldots,N-1$. 
\eth
\prf
Pick a point $x \in \widetilde{\CW}^{\ul{\la}}_\mu$. It follows from Lemma~\ref{Hecke_to_trans} that there exists a weight $w_0\mu_1 \in \La$ of $V^{\la_1}$ such that the Hecke transform $\CP_1(-w_0\mu_1)$ is trivial. There is a unique trivialization of $\CP_1(-w_0\mu_1)$ compatible with our trivialization of $\CP_1$ at $\infty$ (here we use that the automorphism group of $\CP_1(-w_0\mu_1) \simeq \CP^{triv}$ is isomorphic to $G$). After this identification we can assume that $\CP_1(-w_0\mu_1)=\CP^{triv}$. It now follows from the definitions that the point $x \in \widetilde{\CW}^{\ul{\la}}_\mu$ is the image of some point $(x_1,y) \in \overline{\CW}^{\la_1}_{\mu_1} \times \overline{\CW}^{(\la_2,\ldots,\la_N)}_{\mu-\mu_1}$ under the multiplication morphism 
$\overline{\CW}^{\la_1}_{\mu_1} \times \overline{\CW}^{(\la_2,\ldots,\la_N)}_{\mu-\mu_1} \ra \widetilde{\CW}^{\underline{\la}}_\mu$. Continuing applying Lemma~\ref{Hecke_to_trans} we end up with a collection of weights $\mu_1,\ldots,\mu_{N-1} \in \La$ such that $\mu_i$ is a weight of $V^{\la_i}$, and points $x_i \in \overline{\CW}^{\la_i}_{\mu_i}$ such that $\tilde{{\bf{m}}}^{\ul{\la}}_{\ul{\mu}}(x_1,\ldots,x_N)=x$, here $\mu_N:=\mu-\mu_1-\ldots-\mu_{N-1}$. 
\epr

\section{Appendix: some representation-theoretic statements}\label{appendix}
Here we prove some 
facts from representation theory which we use in Remarks~\ref{comb_cor},~\ref{Dinak_exp}.

\lem[See Remark~\ref{comb_cor}]{}\label{weight_rep_mu}
Let $\la \in \La^+$ be a dominant weight (of the Langlands dual group $G^{\vee}$) and $\mu \leqslant \la$ be a weight such that $\langle \al^{\vee},\mu \rangle \geqslant -1$ for any $\al^{\vee} \in \Delta_+^{\vee}$. Then $\mu$ is a weight of $V^{\la}$.
\elem
\prf
Let us decompose $\la-\mu=\beta_1+\ldots+\beta_k$, $\beta_i \in \Delta_+$ with $k$ minimal possible. Set $\mu_i:=\la-\beta_1-\ldots-\beta_i$ ($\mu_k=\mu$). Let us prove by the decreasing induction on $k$ that $\langle\beta_j^{\vee}, \mu_i \rangle \geqslant -1$ for $j \leqslant i$. Indeed for $i=k$ the claim follows from the fact that $\mu_k=\mu$ and our assumptions on $\mu$. Let us prove the induction step. Assume that there exists $j \leqslant i$ such that $\langle\beta_j^{\vee},\mu_i\rangle \leqslant -2$. By the induction hypothesis we have
\begin{equation*}
\langle\beta_j^{\vee},\mu_{i+1}\rangle \geqslant -1 \Rightarrow 
\langle\beta_j^{\vee},\mu_i-\mu_{i+1}\rangle <0 \Rightarrow \langle \beta_j^{\vee},\beta_{i+1} \rangle < 0 \Rightarrow \beta_{i+1}+\beta_j \in \Delta_+
\end{equation*}
so we arrive to the contradiction with the minimality of $k$. 

For any coroot $\beta_i$ consider the corresponding $\mathfrak{sl}_2$-triple $e_{\beta_i},h_{\beta_i},e_{-\beta_i}$ in the Langlands dual Lie algebra $\mathfrak{g}^{\vee}$. Note that $\langle \beta_j^{\vee},\mu_i\rangle =\mu_i(h_{\beta_i})$ so we have shown that $\mu_i(h_{\beta_{i}}) \geqslant -1$ for $j \leqslant i$.
It remains to check that $e_{-\beta_k}\ldots e_{-\beta_1}(v_\la) \neq 0$, where $v_\la \in V^\la$ is a highest weight vector.
We check this by the (increasing) induction on $k$. Case $k=0$ is clear. Let us now assume that $v^{(l-1)}_\la:=e_{-\beta^{\vee}_{l-1}}\ldots e_{-\beta^{\vee}_1}(v_\la) \neq 0$. Consider the subalgebra $\on{Span}_{\BC}(e_{\beta_l},h_{\beta_l},e_{-\beta_l}) \subset \mathfrak{g}^{\vee}$. To check that $e_{-\beta_l}v^{(l-1)}_\la \neq 0$ it is enough to show that $\mu_{l-1}(h_{\beta_l}) > 0$ (here we use that $v^{(l-1)}_\la$ has weight $\mu_{l-1}$). Recall now that 
\begin{equation*}
\mu_{l-1}(h_{\beta_l})=\langle \beta_l^{\vee},\mu_{l-1} \rangle=\langle \beta_l^{\vee},\mu_l\rangle + \langle \beta_l^{\vee}, \beta_l\rangle = 2+\langle \mu_l,\beta_l \rangle \geqslant 1 > 0
\end{equation*}
so the claim follows.
\epr

\lem{}\label{no_dom_est}
Let $\mu \in \La$ be such that there is no dominant $\la' \in \La^+$ such that $\mu \leqslant \la' < \mu^+$. Then $\langle\al^{\vee},\mu\rangle \geqslant -1$ for every $\al^{\vee} \in \Delta^{\vee}_+$.
\elem
\prf
The proof is the modification of~\cite[VIII,~\S 7,~no.~3,~Proposition~6]{Bu} for the minuscule case.
Set $\la:=\mu^+$.
Assume that there exists $\al^{\vee} \in \Delta^{\vee}_+$ such that $\langle\al^{\vee},\mu\rangle \leqslant -2$. Consider the weight $\mu':=s_{\al}(\mu)=\mu-\langle\al^{\vee},\mu\rangle\al$ and note that  $\la \geqslant \mu' > \mu$ (since $\langle\al^{\vee},\mu\rangle<0$) and $\langle\al^{\vee},\mu'\rangle=-\langle\al^{\vee},\mu\rangle \geqslant 2$. Consider now the weight $\nu:=\mu'-\al$. Again we have $\la \geqslant \mu'>\nu>\mu$. We claim that 
\begin{equation}\label{ineq_haha}
(\mu',\mu')>(\nu,\nu),   
\end{equation}
where $(\,,\,)$ is a nondegenerate symmetric bilinear form on $\La_{\BC}:=\La \otimes_{\BZ} \BC$ such that 
\begin{equation*}
\langle \al^{\vee},\beta \rangle=2\frac{(\al,\beta)}{(\al,\al)}~\text{for}~\al \in \Delta,\, \beta \in \La_{\BC},\, (\al,\al)>0~\text{for}~\al \in \Delta.
\end{equation*}
To see~(\ref{ineq_haha}) note that
\begin{equation*}
(\mu',\mu')-(\nu,\nu)=2(\al,\mu')-(\al,\al)=\frac{(\al,\al)}{2} (2\langle \al^{\vee},\mu' \rangle - 2)>0,
\end{equation*}
where the last inequality holds since $\langle\al^{\vee},\mu'\rangle \geqslant 2 >1,\, (\al,\al)>0$. Let $\nu^+ \in W\nu$ be the dominant representative of $\nu$. Since $\nu$ is a weight of $V^\la$ we must have $\la \geqslant \nu^+$. We also have $\nu^+ \geqslant \nu > \mu$ so we conclude that $\mu^+=\la \geqslant \nu^+ > \mu$. It remains to note that 
\begin{equation*}
(\nu^+,\nu^+)=(\nu,\nu)<(\mu',\mu')=(s_{\al}(\mu),s_{\al}(\mu))=(\mu,\mu)=(\la,\la)
\end{equation*}
so $\la \neq \nu^+$ and contradiction finishes the proof. 
\epr

\prop{}\label{conditions_pairing_orbit}
The following conditions on $\mu \in \La$ are equivalent. 

(1) For any positive $\al \in \Delta_+$ we have $\langle \al,\mu \rangle \geqslant -1$, 

(2) For any dominant $\la' \in \La^+$ the inequality $\mu \leqslant \la' \leqslant \mu^+$ implies $\la'=\mu^+$.
\eprop
\prf
Let us prove the implication $(1) \Rightarrow (2)$. Let $\la' \in \La^+$ be a dominant coweight such that $\mu \leqslant \la' \leqslant \mu^+$. Then by Lemma~\ref{weight_rep_mu} $\mu$ is a weight of $V^{\la'}$ so $\mu^+$ is also a weight of $V^{\la'}$, hence, $\mu^+ \leqslant \la'$. It follows that $\mu^+=\la'$.

Let us now prove the implication $(2) \Rightarrow (1)$. This is exactly the content of Lemma~\ref{no_dom_est}.
\epr

\rem{}
{{
Note that for antidominant $\mu$ the conditions $(1),\, (2)$ of Proposition~\ref{conditions_pairing_orbit} are equivalent to the fact that $\mu^+$ is minuscule.}}
\erem

\begin{Cor}[see Remark~\ref{Dinak_exp}]{}\label{equiv_deepest_comb}
Let $\la \in \La^+$ be dominant and $\mu \leqslant \la$. The following  are equivalent. 

(1) $\CW^{\mu^+}_{\mu} \subset \ol{\CW}^{\la}_\mu$ is the deepest stratum i.e. that if we write $\ol{\CW}^{\la}_\mu=\bigsqcup_{\mu \leqslant \la' \leqslant \la}\CW^{\la'}_\mu$ then $\CW^{\mu^+}_\mu \subset \ol{\CW}^{\la'}_\mu$.

(2) We have $\ol{\CW}^{\mu^+}_\mu=\CW^{\mu^+}_{\mu}$.

(3) We have $\langle \al^{\vee},\mu \rangle \geqslant -1$ for any $\al^{\vee} \in \Delta^{\vee}_+$.
\end{Cor}
\prf
The implication $(1) \Rightarrow (2)$ is clear.
The implication $(2) \Rightarrow (3)$ follows from Lemma~\ref{no_dom_est}. 

Let us now prove the implication $(3) \Rightarrow (1)$. Pick $\la' \in \La^+$ such that $\mu \leqslant \la' \leqslant \la$. Our goal is to show that $\CW^{\mu^+}_{\mu} \subset \ol{\CW}^{\la'}_{\mu}$ i.e. that $\la' \geqslant \mu^+$. Lemma~\ref{weight_rep_mu} implies that $\mu$ is a weight of $V^{\la'}$, hence, $\mu^{+}$ is a weight of $V^{\la'}$ so $\mu^+ \leqslant \la'$ and the claim follows.
\epr

\end{document}